\documentclass[a4paper,10pt]{article}
\usepackage{latexsym}

\newcommand{\thm}{\sc{Theorem}\ }
\newcommand{\pro}{\sc{Proposition}\ }

\newcommand{\lem}{\sc{Lemma}\ }
\newcommand{\proof}{\sc{Proof.}}
\usepackage{color}
\begin{document}
\begin{center}
  {\bf \Large Generic conformally flat hypersurfaces and surfaces in 3-sphere\footnote{accepted for publication in SCIENCE CHINA Mathematics}}
\end{center}

\vspace*{2mm}
\begin{center}
  {\it Dedicated to Professor Qing-Ming Cheng for his 60th birthday}
\end{center} 
 
\vspace{2mm}
\begin{center}
{\Large Yoshihiko Suyama\footnote{Department of Applied Mathematics, Fukuoka University, Fukuoka 814-0180 (Japan), \\ 
\hspace*{1cm}e-mail: suyama@fukuoka-u.ac.jp}}
\end{center}

\vspace*{2mm}
{\bf Abstract} \ The aim of this paper is to verify that the study of generic conformally flat hypersurfaces in 4-dimensional space forms is reduced to a surface theory in the standard 3-sphere. 
The conformal structure of generic conformally flat (local-)hypersurfaces is characterized as conformally flat (local-)3-metrics with the Guichard condition.
Then, there is a certain class of orthogonal analytic (local-)Riemannian 2-metrics with constant Gauss curvature -1 such that any 2-metric of the class gives rise to a one-parameter family of conformally flat 3-metrics with the Guichard condition.  
In this paper, we firstly relate 2-metrics of the class to surfaces in the 3-sphere:  
for a 2-metric of the class, a 5-dimensional set of (non-isometric) analytic surfaces in the 3-sphere is determined such that 
any surface of the set gives rise to an evolution of surfaces in the 3-sphere issuing from the surface 
and the evolution is the Gauss map of a generic conformally flat hypersurface in the Euclidean 4-space. Secondly, we characterize analytic surfaces in the 3-sphere which give rise to generic conformally flat hypersurfaces. \\



\vspace{2mm}
\noindent
{\large\bf 1. Introduction}

The aim of this paper is to verify that the study of generic conformally flat hypersurfaces in 4-dimensional space forms 
is reduced to a surface theory in the standard 3-sphere $S^3$. 
Here, we say that a hypersurface is generic if it has distinct three principal curvatures at each point. 

\vspace{2mm}
The classification of conformally flat hypersurfaces in $(n+1)$-dimensional space forms was determined by Cartan in paper (\cite{ca}) for $n\geq 4$: 
a hypersurface in an $(n+1)$-dimensional space form is conformally flat if and only if it is a branched channel hypersurface. 
Now, although 3-dimensional branched channel hypersurfaces are conformally flat as well, in the case of $n=3$ there are also generic 3-dimensional conformally flat hypersurfaces. 
Our theme of this paper is to study these generic 3-dimensional conformally flat hypersurfaces: 
we relate generic conformally flat hypersurfaces to analytic surfaces in $S^3$.

\vspace{2mm}
Any generic conformally flat hypersurfaces in a 4-dimensional space form has a special principal curvature line coordinate system $(x,y,z)$: 
its first fundamental form $I$ is expressed as \ $I=l_1^2(dx)^2+l_2^2(dy)^2+l_3^2(dz)^2$ \ 
and, for the functions $l_i^2$ ($i=1,2,3$), there is a permutation $\{i,j,k\}$ of $\{1,2,3\}$ such that a Guichard condition $l_i^2+l_j^2=l_k^2$ is satisfied. 
Such a coordinate system is called a principal Guichard net of a generic conformally flat hypersurface, and  
a pair $\{(x,y,z),[g]\}$ is called a Guichard net, where $(x,y,z)$ is a coordinate system on a simply connected  domain $U$ in the Euclidean 3-space $R^3$ and $[g]$ is the conformal class of a conformally flat metric $g$ satisfying the Guichard condition with respect to the coordinate system.

Conversely, up to M$\ddot{\rm o}$bius transformation, there uniquely exists, for a given Guichard net $\{(x,y,z),[g]\}$, a generic conformally flat hypersurface with its {\it canonical} principal Guichard net in a 4-dimensional space form (for the canonical principal Guichard net, see Remark 1 in \S2.2). 
Hence, we describe the canonical principal Guichard net for a hypersurface as only ``the Guichard net''.   
This existence theorem was induced by the integrability condition on a generic conformally flat hypersurface with 
the Guichard net in the conformal 4-sphere (cf. \cite{he1}, also see \cite{hs1}). 

Now, let $f=f(x,y,z)$ be a generic conformally flat hypersurface in $R^4$. 
Let $\kappa_i$ ($i=1,2,3$) be the principal curvatures of $f$ corresponding to the coordinate lines $x$, $y$ and $z$, respectively. 
For the sake of simplicity suppose that $\kappa_3$ 
is the middle principal curvature, i.e.,\ $\kappa_1>\kappa_3>\kappa_2$\ or\ $\kappa_1<\kappa_3<\kappa_2$.\ 
Then, by the Guichard condition there is a conformally flat metric $g$ such that $g$ is expressed as
$$g=\cos^2\varphi(dx)^2+\sin^2\varphi(dy)^2+(dz)^2  \eqno{(1.1)}$$
with a function $\varphi=\varphi(x,y,z)$ and  the pair of $g$ and  
the coordinate system $(x,y,z)$ is a representative of the Guichard net determined by $f$, furthermore, 
the first fundamental form $I_f$ of $f$ is given by 
$$I_f=e^{2P}g:=e^{2P}(\cos^2\varphi(dx)^2+\sin^2\varphi(dy)^2+(dz)^2)  \eqno{(1.2)}$$
with a function $P=P(x,y,z)$ on $U$. 
Therefore, the existence problem of generic conformally flat hypersurfaces is reduced to that of conformally flat metrics $g$ (or functions $\varphi$) given by (1.1). In the case that $\varphi$ satisfies \ $\varphi_{zx}=\varphi_{zy}\equiv 0$, \ all hypersurfaces with such Guichard nets were classified and constructed in 4-dimensional space forms by a series of papers (\cite{su2}, \cite{su3}, \cite{hs1}), 
where $\varphi_{zx}$ is the second derivative of $\varphi$ with respect to $z$ and $x$. 
In this paper, we mainly study the case that $\varphi$ satisfies $\varphi_{zx}\varphi_{zy}\neq 0$. 

For the sake of simplicity for the description later, we assume that the domain $U$, where $g$ is defined, is given by $U:=V\times I\subset R^2\times R$, where $0\in I$, in particular, $U$ intersects the plane $z=0$.

Let $g$ be a conformally flat metric given by (1.1) from a function $\varphi(x,y,z)$. 
If $\varphi(x,y,z)$ satisfies $\varphi_z\varphi_{zx}\varphi_{zy}\neq 0$, \ then $\varphi(x,y,z)$ (or $g$) gives rise to an evolution $\hat g(z), \ z\in I,$ of orthogonal 2-metrics on $V$ with constant Gauss curvature $-1$ (\cite{bhs}). 
The converse problem was also studied in \cite{bhs}, 
and then a class $Met^0$ of orthogonal {\it analytic} 2-metrics on $V$ with constant Gauss curvature $-1$ was defined: 
for any $\hat g\in Met^0$, 
a one-parameter family $g^t$ (with parameter $t\in R\setminus \{0\}$) of conformally flat metrics on $U$ given by (1.1) is determined as evolutions $\hat g^t(z)$ of 2-metrics on $V$ issuing from $\hat g^t(0)=\hat g$. 
In other words, any $\hat g\in Met^0$ determines a one-parameter family $(\varphi(x,y),\varphi_z^t(x,y))$ of two real-analytic functions on $V$ and 
a solution $\varphi^t(x,y,z)$ to certain evolution equation in $z$ under these initial data $\varphi^t(x,y,0)=\varphi(x,y)$ and $\varphi_z^t(x,y,0)=\varphi_z^t(x,y)$ determines $g^t$ by (1.1) (more precisely, see \S2.1). Note that, in the above argument, the evolution is taking in $z$-direction (not  in $t$-direction).

By use of the notation $\varphi^t_z(x,y)$ above, we may anticipate that it would be the $z$-derivative $\varphi^t_z(x,y,0)$ of some function 
$\varphi^t(x,y,z)$. 
The Cauchy-Kovalevskaya theorem on analytic evolution equations ensures the existence and the uniqueness of $\varphi^t(x,y,z)$ under real-analytic initial data. The theorem is also fundamental for our consideration in this paper.

\vspace*{2mm}
Now, let $g$ be a conformally flat metric given by (1.1) from $\varphi(x,y,z)$ and $f(x,y,z)$ be a generic conformally flat hypersurface in $R^4$ with the Guichard net $g$. Let $I_f=e^{2P}g$ be the first fundamental form of $f$. 
Let \ $\{X_{\alpha},X_{\beta},X_{\gamma},N\}$ \ be an orthonormal frame field of $R^4$ along $f(x,y,z)$, where $X_{\alpha},$ $X_{\beta}$ and $X_{\gamma}$ are orthonormal principal vector fields in the direction of $f_x,$ $f_y$ and $f_z$, respectively, and $N$ is a unit normal vector field. 
Suppose that $(\kappa_1\kappa_2)(x,y,0)\neq 0$ is satisfied. 
Then, a mapping $\phi(x,y):=N(x,y,0)$ for $(x,y)\in V$ defines a surface in $S^3$ and $(x,y)$ is a principal curvature line coordinate system of $\phi$. 
Let $X^0_{\alpha}(x,y):=X_{\alpha}(x,y,0)$, $X^0_{\beta}(x,y):=X_{\beta}(x,y,0)$ and $\xi(x,y):=X_{\gamma}(x,y,0)$. 
Then, $X^0_{\alpha}(x,y)$ and $X^0_{\beta}(x,y)$ are principal curvature vector fields and $\xi(x,y)$ is a unit normal vector field (in $S^3$) of $\phi$. 
We further define several functions on $V$ from ones on $U$ as follows: 
$$\varphi(x,y):=\varphi(x,y,0), \ \ \varphi_z(x,y):=\varphi_z(x,y,0), \ \ 
e^{\bar P}(x,y):=e^P(x,y,0),$$
$$(e^{\bar P})_z(x,y):=(e^P)_z(x,y,0), \ \ \ \bar\kappa_i(x,y):=\kappa_i(x,y,0).$$ 
Then, the surface $\phi$ satisfies the following conditions (a), (b) and (c) (see \S2.2), where $D'$ is the canonical Riemannian connection on $S^3$:

$$\phi_x=-{\bar \kappa}_1e^{\bar P}\cos\varphi\ X^0_{\alpha}, \ \ \ \ 
\phi_y=-{\bar \kappa}_2e^{\bar P}\sin\varphi\ X^0_{\beta}.  \leqno{ \ \ \ \ \ (a)}$$

$$D'_{\partial/\partial y}X^0_{\alpha}={(e^{\bar P}\sin\varphi)_x}/{(e^{\bar P}\cos\varphi)}~X_{\beta}^0, \ \ \ \ 
D'_{\partial/\partial x}X^0_{\beta}={(e^{\bar P}\cos\varphi)_y}/{(e^{\bar P}\sin\varphi)}~X_{\alpha}^0. \leqno{ \ \ \ \ \ (b)}$$

$$D'_{\partial/\partial x}\xi=e^{-\bar P}(e^{\bar P}\cos\varphi)_z\ X_{\alpha}^0, \ \ \ \ 
D'_{\partial/\partial y}\xi=e^{-\bar P}(e^{\bar P}\sin\varphi)_z\ X_{\beta}^0,  \leqno{ \ \ \ \ \ (c)}$$
where \ $(e^{\bar P}\cos\varphi)_z:=(e^{\bar P})_z\cos\varphi-e^{\bar P}\varphi_z\sin\varphi, \ \  
(e^{\bar P}\sin\varphi)_z:=(e^{\bar P})_z\sin\varphi+e^{\bar P}\varphi_z\cos\varphi.$

\vspace*{2mm}
Our first aim is to study the converse problem by starting from $\hat g\in Met^0$. 
In Main Theorem 1 below, ${\bar P}(x,y)$ and ${\bar P}_z(x,y)$ are analytic functions on $V$ independent of each other, and we denote $(e^{\bar P})_z(x,y):=(e^{\bar P}\bar P_z)(x,y)$. 
Let $(\varphi(x,y),\varphi_z^t(x,y))$ be 
a pair of functions on $V$ determined by $\hat g\in Met^0$ and arbitrarily fixed $t\neq 0$. Then, $\varphi^t(x,y,z)$ is determined as an evolution from the initial data $(\varphi(x,y),\varphi_z^t(x,y))$, as mentioned above. 

For given $(\varphi(x,y),\varphi_z^t(x,y))$ and $({\bar P}(x,y),{\bar P}_z(x,y))$, three functions $\bar\kappa_i(x,y)$ $(i=1,2,3)$ on $V$, respectively, are uniquely determined (see Definition and Notation 3 in \S3). 
A family \ $\{e^P(x,y,z), \kappa_i(x,y,z) | \ i=1,2,3\}$ \ of analytic functions on $U$ is called an extension of 
\ $\{e^{\bar P}(x,y), (e^{\bar P})_z(x,y),\bar\kappa_i(x,y)| \ i=1,2,3 \}$, \ if   
$$e^P(x,y,0)=e^{\bar P}(x,y), \ \ \ (e^P)_z(x,y,0)=(e^{\bar P})_z(x,y), \ \ \  \kappa_i(x,y,0)=\bar\kappa_i(x,y)$$ 
are satisfied on $V$. Hence, $P(x,y,z)$ is also an analytic function on $U$. Let us denote by $f(x,y,z)$ (resp. $\phi(x,y)$) a hypersurface in $R^4$ (resp. a surface in $S^3$). \\ 

\vspace*{2mm}
{\sc Main Theorem 1.} \ \ {\it Let $\hat g\in Met^0$, and $\Phi^t=(\varphi(x,y),\varphi_z^t(x,y))$ be 
a pair of functions on $V$ determined by ${\hat g}$ and arbitrarily fixed $t\neq 0$.   
For $\hat g$ and $\Phi^t$, we have the following facts (1)$\sim$(3): \\[-3mm]

(1) A 5-dimensional set of pairs $({\bar P}(x,y), {\bar P}_z(x,y) )$ consisting of analytic functions  on $V$ is \\
\hspace*{0.8cm}determined: for any pair $({\bar P}(x,y), {\bar P}_z(x,y) )$ of the set, there is an analytic surface $\phi(x,y)$ in $S^3$ \\
\hspace*{0.8cm}defined for $(x,y)\in V$ such that $\phi$ satisfies the above conditions (a), (b) and (c), by replacing \\
\hspace*{0.8cm}$\varphi_z(x,y)$ in (c) by $\varphi^t_z(x,y)$, with respect to an orthonormal frame field $\{X^0_{\alpha}(x,y),X^0_{\beta}(x,y),\xi(x,y)\}$.\\[-3mm]

(2) For any class $\{e^{\bar P}(x,y), (e^{\bar P})_z(x,y), \bar\kappa_i(x,y) \}$ in (1), its extension $\{e^P(x,y,z), \kappa_i(x,y,z)\}$ to $U$ \\  
\hspace*{0.8cm}is determined such that there is an evolution $\phi^z(x,y)$, $z\in I$, of surfaces in $S^3$ issuing from \\ 
\hspace*{0.8cm}$\phi^0(x,y)=\phi(x,y)$ and each surface $\phi^z$ with $z$ also satisfies (a), (b) and (c): describing (a) \\
\hspace*{0.8cm}explicitly, the differential of each surface $\phi^z(x,y)$ with $z$ is expressed as 
$$(d\phi^z)(x,y)=-(\kappa_1e^P\cos\varphi^t)(x,y,z)X_{\alpha}(x,y,z)dx-(\kappa_2e^P\sin\varphi^t)(x,y,z)X_{\beta}(x,y,z)dy$$
\hspace*{0.8cm}with respect to an orthonormal frame field $\{X_{\alpha}(x,y,z), X_{\beta}(x,y,z)\}$ satisfying \\
\hspace*{0.8cm}$X_{\alpha}(x,y,0)=X_{\alpha}^0(x,y)$ and $X_{\beta}(x,y,0)=X_{\beta}^0(x,y)$.\\[-3mm]

(3) Each evolution $\phi^z(x,y)$ in (2) gives rise to an evolution $f^z(x,y)$ of analytic surfaces in \\
\hspace*{0.8cm}$R^4$ such that the differential of each surface $f^z(x,y)$ with $z$ satisfies  
$$ (df^z)(x,y)=e^{P(x,y,z)}\{\cos\varphi^t(x,y,z)X_{\alpha}(x,y,z)dx+\sin\varphi^t(x,y,z)X_{\beta}(x,y,z)dy\}$$ 
\hspace*{0.8cm}and $f(x,y,z):=f^z(x,y)$ is a generic conformally flat hypersurface in $R^4$ defined on $U$. Then, \\
\hspace*{0.8cm}$f(x,y,z)$ has the Guichard net $g^t$ defined by (1.1) from $\varphi^t(x,y,z)$, the first fundamental form \\
\hspace*{0.8cm}$I_f=e^{2P}g^t$ and principal curvatures $\kappa_i$. In particular, $N(x,y,z):=\phi^z(x,y)$ is the Gauss map \\ 
\hspace*{0.8cm}of $f(x,y,z)$.}\\

\vspace*{2mm}
In order to verify Main Theorem 1, we firstly review, in \S2.3, the result of \cite{hs3}: 
for a conformally flat metric $g$ given by (1.1), the problem of determining the first and the second fundamental forms for a generic conformally flat 
hypersurface $f$ realized in $R^4$ is reduced to that of finding only one function $P(x,y,z)$ on $U$, which determines the conformal element $e^{2P}$ of $I_f$ in (1.2).

Next, we show, in \S3, that the functions $P(x,y,z)$ are 
solutions to a certain evolution equation in $z$ under suitable initial conditions $({\bar P}(x,y),\bar P_z(x,y))$ on $V\times \{0\}.$
These initial data $({\bar P}(x,y),{\bar P}_z(x,y))$ are given as solutions to a system of differential equations on $V$ defined by $\Phi^t$ (Theorem 1, and also see Definition and Notation 1 in \S2.1). 
Then, in order that $P(x,y,z)$ is actually determined as such a solution, it is necessary that $\bar\kappa_1(x,y)$ and $\bar\kappa_2(x,y)$ determined from $({\bar P}(x,y),{\bar P}_z(x,y))$ satisfy $(\bar\kappa_1\bar\kappa_2)(x,y)\neq 0$ on $V$. 
However, when $P(x,y,z)$ is attendant on a generic conformally flat hypersurface with a Guichard net $g^t$ arising from $\hat g\in Met^0$, it always satisfies the condition $(\kappa_1\kappa_2)(x,y,0)=(\bar\kappa_1\bar\kappa_2)(x,y)\neq 0$ on $V$, by choosing a suitable subdomain $V'$ of $V$ if necessary (Proposition 3.3).   
Through examples, for a given $\Phi^t$ we show that the system for $({\bar P}(x,y),{\bar P}_z(x,y))$ can be completely solved and all 5-dimensional solutions are obtained. 

The existence of surfaces $\phi$ (resp. evolutions $\phi^z$) in $S^3$ follows from a geometrical interpretation of the system for $({\bar P}(x,y),{\bar P}_z(x,y))$, which is verified in \S4 by Theorem 2 (resp. Theorem 3). 
If two surfaces $\phi$ and $\bar\phi$ in $S^3$ are isometric, then we regard $\phi$ and $\bar\phi$ as the same surface. 
Then, the set of all surfaces $\phi(x,y)$ for any $\Phi^t$ is 4-dimensional. Hence, all surfaces $\phi$ for any $\hat g\in Met^0$ generate a 5-dimensional set (Remark 4-(3)).

In \S5, we characterize analytic surfaces $\phi(x,y)$ in $S^3$ such that they give rise to generic conformally flat hypersurfaces $f$ in $R^4$. 
Let us express the differential of $\phi$ as \ $d\phi=-(a_1dx)X^0_{\alpha}-(a_2dy)X^0_{\beta}$, \ where $\{X^0_{\alpha},X^0_{\beta}\}$ is defined as above. Then, we have the following Main Theorem 2: \\

\vspace*{2mm}
{\sc Main Theorem 2.} \ \ {\it Let $\phi:V\ni (x,y)\mapsto \phi(x,y)\in S^3$ be a generic analytic surface and $(x,y)$ be a principal curvature line coordinate system.   
Suppose that, for $\phi$, there is a principal curvature line coordinate system $(x,y)$ and analytic functions $(\varphi( x, y),\bar P( x, y))$ on $V$ such that they satisfy \ 
$$(1) \ a_1\sin\varphi- a_2\cos\varphi=1 \ (resp. -1), \ \ \ \ 
(2) \ (a_1)_{ x}=-\varphi_{ x} a_2+\bar P_{ x} a_1, \ \  ( a_2)_{ y}=\varphi_{ y} a_1+\bar P_{ y} a_2.$$  
Then, three objects: an evolution $\varphi(x,y,z)$ on $U$ of $\varphi(x,y)$, an evolution $\phi^z(x,y)$ of surfaces issuing from $\phi$ and a generic conformally flat hypersurface $f(x,y,z)$ in $R^4$, are determined such that $\varphi(x,y,z)$ gives the Guichard net $g$ of $f$ by (1.1) and $\phi^z(x,y)$ is the Gauss map of $f$.  
The surface $\phi(x,y)$ further gives rise to the dual $f^*(x,y,z)$ of $f$, hence there is another pair $(\varphi^*( x, y),\bar P^*( x, y))$ leading to $f^*$ such that it satisfies (1) and (2), where we replace only the equation (1) by \ $a_1\sin\varphi^*- a_2\cos\varphi^*=-1 \ (resp.~ 1)$. 
In particular, the principal curvatures in the direction of $z$ for $f$ and $f^*$ are middle among three principal curvatures of each hypersurface. \\ }

We verify Main Theorem 2 by Theorem 4 and  Corollaries 5.1 and 5.2 together with reviews of dual generic conformally flat hypersurfaces in (\cite{hsuy} and \cite{hs3}).  Owing to Main Theorem 2, we can judge by only the first fundamental form whether a surface $\phi$ leads to a generic conformally flat hypersurface or not, and then it is possible that the Guichard net $g$ of the hypersurface is not the one arising from $\hat g\in Met^0$ (Remark 5-(2) and (3)).

\vspace*{2mm}
At the end of the introduction, we refer to some recent results on generic conformally flat hypersurfaces: 
S.Canevari and R.Tojeiro (\cite{ct}, \cite{ct2}) provided another characterization of such hypersurfaces in $R^4$ (different from the existence of the Guichard net) (see Remark 2 in \S3); 
Z.X.~Xie, C.P.~Wang and X.Z.~Wang (\cite{wxcwxw}, \cite{xzcwxw}, \cite{xwcwwx}) are studying conformally flat Lorentzian hypersurfaces in the Lorenzian space $R^4_1$. \\

\newpage
\vspace{2mm}
\noindent
{\large\bf 2. Preliminaries}

In this section, we summarize known results on generic conformally flat hypersurfaces in addition to them in the introduction, and we fix our notations for the description later.

\vspace{2mm}
\noindent
{\bf 2.1. Guichard nets and evolutions of 2-metrics with constant curvature $-1$} 

In this subsection, we summarize results in \cite{bhs}. Let $g$ be a conformally flat metric expressed as 
$$g=\cos^2\varphi(dx)^2+\sin^2\varphi(dy)^2+(dz)^2.  \eqno{(2.1.1)}$$ 
Here, $\varphi=\varphi(x,y,z)$ is a function on a simply connected domain \ 
$U=V\times I\subset R^2\times R$, where $0\in I$. 
When a conformally flat metric $g$ is expressed as (2.1.1) with respect to a coordinate system $(x,y,z)$, $g$ is called a 
conformally flat metric with the Guichard condition. 
As mentioned in the introduction, for a given conformally flat metric $g$ with the Guichard condition, there is a generic conformally 
flat hypersurface $f(x,y,z)$ in $R^4$ with the canonical principal Guichard net $g$, uniquely up to a conformal transformation.
Note that two conformally flat metrics $g$ determined from $\varphi(x,y,z)$ and $\tilde\varphi(x,y,z)$ define the same Guichard net if and only if there are three constants $a_1$, $a_2$ and $a_3$ such that \ $\tilde\varphi(x,y,z)=\varphi(\pm x+a_1,\pm y+a_2,\pm z+a_3)$. 

Now, we have the following proposition (cf. \cite{bhs}, Proposition 1.1 and Theorem 1).

\vspace{2mm}
{\pro 2.1}. \ \ {\it The following facts (1), (2) and (3) are equivalent to each other. 

(1) A metric $g$ on $U$ given by (2.1.1) is conformally flat.

(2) The function $\varphi(x,y,z)$ on $U$ determining $g$ in (2.1.1) satisfies the following equations (i)-(iv): 

$$ \varphi_{xyz}+\varphi_x\varphi_{yz}\tan\varphi-\varphi_y\varphi_{xz}\cot\varphi=0,  \leqno{(i)}$$
$$\frac{\varphi_{xxx}-\varphi_{yyx}+\varphi_{zzx}}{2}-\frac{(\varphi_{xx}-\varphi_{yy})\cos2\varphi-\varphi_{zz}}
{\sin2\varphi}\varphi_x-\varphi_{xz}\varphi_z\cot\varphi=0,  \leqno{(ii)}$$
$$\frac{\varphi_{xxy}-\varphi_{yyy}-\varphi_{zzy}}{2}-\frac{(\varphi_{xx}-\varphi_{yy})\cos2\varphi
-\varphi_{zz}}{\sin2\varphi}\varphi_y-\varphi_{yz}\varphi_z\tan\varphi=0,  \leqno{(iii)}$$
$$\frac{\varphi_{xxz}+\varphi_{yyz}+\varphi_{zzz}}{2}+\frac{\varphi_{xx}-\varphi_{yy}-\varphi_{zz}\cos2\varphi}
{\sin2\varphi}\varphi_z-\varphi_x\varphi_{xz}\cot\varphi+\varphi_y\varphi_{yz}\tan\varphi=0.  \leqno{(iv)}$$ 

(3) Let $\varphi(x,y,z)$ and $\psi(x,y,z)$ be functions on $U$ satisfying the following equations (i)-(iv): 

$$(i)\hspace{2mm} \psi_{xz}=-\varphi_{xz}\cot\varphi, \hspace{1.2cm} (ii)\hspace{2mm}  \psi_{yz}=\varphi_{yz}\tan\varphi,$$
$$ (iii)\hspace{2mm} \psi_{zz}=(\varphi_{xx}-\varphi_{yy})\sin2\varphi-(\psi_{xx}-\psi_{yy})\cos2\varphi,$$
$$ (iv)\hspace{2mm} \varphi_{zz}=(\varphi_{xx}-\varphi_{yy})\cos2\varphi+(\psi_{xx}-\psi_{yy})\sin2\varphi.$$ 
Then, $g$ is defined from $\varphi$ by (2.1.1).
In particular, we can assume that $\psi$ does not have any linear term 
for $x$, $y$ and $z$.  }\\

\vspace{2mm}
So far as we look at Proposition 2.1-(3), 
$\psi(x,y,z)$ is not unique even if it does not have any linear term for $x$, $y$ and $z$.
However, we can impose the following additional constraints on $\psi(x,y,z)$ (cf. \cite{bhs}, Theorem 4 in \S3.1);   
$$\psi_{xy}=\varphi_x\varphi_y, \ \ \ \ (\varphi_{xx}-\varphi_{yy})\sin2\varphi-(\psi_{xx}-\psi_{yy})\cos2\varphi
=-\Delta\psi+\varphi_x^2+\varphi_y^2+\varphi{_z^2},  \eqno{(2.1.2)}$$
where \ $\Delta \psi:=(\frac{\partial^2}{\partial x^2}+\frac{\partial^2}{\partial y^2})\psi$. 
Then, $\psi(x,y,z)$ is uniquely determined for $\varphi(x,y,z)$ (\cite{bhs}, Proposition 4.1).

Now, let $g$ be a conformally flat metric defined by (2.1.1) from $\varphi(x,y,z)$. 
Suppose that $\varphi$ satisfies the condition \ $(\varphi_z\varphi_{xz}\varphi_{yz})(x,y,z)\neq 0$ \ 
on $U$. 
For $g$, we define functions $\hat{A}(x,y,z)$ and $\hat{B}(x,y,z)$ on $U$ by $$\hat{A}:=-\frac{\varphi_{zx}}{\varphi_z\sin\varphi}=\frac{\psi_{zx}}{\varphi_z\cos\varphi}, \hspace*{0.5cm}
\hat{B}:=\frac{\varphi_{zy}}{\varphi_z\cos\varphi}=\frac{\psi_{zy}}{\varphi_z\sin\varphi},  \eqno{(2.1.3)}$$
respectively, which are well-defined by (i) and (ii) of Proposition 2.1-(3). Then, each 2-metric \ $\hat{g}(z):=\hat{A}^2(x,y,z)(dx)^2+\hat{B}^2(x,y,z)(dy)^2$ \ on $V$ with $z\in I$ has constant Gauss curvature $-1$. That is, for the metric $g$, 
an evolution $\hat{g}(z)$ of 2-metrics on $V$ with constant Gauss curvature $-1$ is determined. 

In \cite{bhs}, the inverse problem was also studied, and then a class $Met^0$ of orthogonal {\it analytic} 2-metrics on $V$ with constant Gauss curvature $-1$ is defined:    
for a 2-metric $\hat g\in Met^0$, a one-parameter family $g^t$ (with parameter $t\in R\setminus\{0\}$) of conformally flat metrics on $U$ with the Guichard condition is determined as evolutions ${\hat g}^t(z)$ of 2-metrics on $V$ with constant Gauss curvature $-1$ issuing from ${\hat g}^t(0)={\hat g}$.  
We review the method of determining the family $g^t$ from ${\hat g}\in Met^0$ in the following.  

We interpret (iii) and (iv) of Proposition 2.1-(3) as a system of evolution equations in $z$ for $\psi(x,y,z)$ and $\varphi(x,y,z)$. 
Hence, we must determine its suitable initial data on $z=0$; four functions $\varphi(x,y)$, $\psi(x,y)$, $\varphi_z(x,y)$ and $\psi_z(x,y)$ on $V$ for $\hat g\in Met^0$. In this process, the definition of $Met^0$ is determined. 
To see it, let \ $\hat{g}=\hat{A}^2(x,y)(dx)^2+\hat{B}^2(x,y)(dy)^2$ \ be an arbitrary analytic 2-metric on $V$ with constant Gauss curvature $-1$ and let $(x_0,y_0)$ be a point of $V$. When we regard $\hat{A}$ and $\hat{B}$ for $\hat g$ as them in (2.1.3) on $z=0$, \ three functions $\varphi(x,y)$, $\varphi_z(x,y)$ and $\psi_z(x,y)$ are determined from $\hat{g}$ as follows: $\varphi$ is uniquely determined by giving $\varphi(x_0,y_0)=\lambda$, but both $\varphi_z$ and $\psi_z$ are only determined up to the same constant multiple $t\neq 0$ even if we assume $\psi_z(x_0,y_0)=0$, \  that is, 
\ $\varphi_z=\varphi_z^t:=t\varphi_z^1$ \ and \ $\psi_z=\psi_z^t:=t\psi_z^1$ \ for $t\in \ R\setminus\{0\}$. In general, $\psi(x,y)$ is not determined from $\hat g$ itself. 

Suppose that, for $\hat g$, there is an evolution $\hat g(z)$ of 2-metrics issuing from $\hat g$ such that $\hat g(z)$ leads to a conformally flat metric $g$ with the Guichard condition. 
Let $\varphi(x,y,z)$ and $\psi(x,y,z)$ be functions on $U$ in Proposition 2.1-(3) determined for $g$. 
Then, $\varphi(x,y,z)$ and $\psi(x,y,z)$ satisfy 
$$[\psi_{xz}+\varphi_{xz}\cot\varphi]_z(x,y,0)=0, \hspace*{1cm}
[\psi_{yz}-\varphi_{yz}\tan\varphi]_z(x,y,0)=0. \eqno{(2.1.4)}$$
In (2.1.4), we interpret \ $\psi_{xzz}(x,y,0):=\psi_{zzx}(x,y,0)$, \ $\varphi_{xzz}(x,y,0):=\varphi_{zzx}(x,y,0)$, \ $(\cot\varphi)_z(x,y,0):=-[\varphi_z/\sin^2\varphi](x,y,0)$ \ and so on. Then, considering (iii) and (iv) of Proposition 2.1-(3), we obtain a certain condition on $\psi(x,y,0)$ by (2.1.4), which is determined from $\varphi(x,y,0)$ and $\varphi_z(x,y,0)$. 
Now, we say that $\hat g$ belongs to $Met^0$ if all classes $(\varphi(x,y), \varphi^t_z(x,y),\psi^t_z(x,y))$ with $t\in R\setminus \{0\}$ determined for $\hat{g}$ satisfy (2.1.4). Then, for $\hat g\in Met^0$, each class $(\varphi(x,y), \varphi^t_z(x,y),\psi^t(x,y),\psi^t_z(x,y))$ with $t\in R\setminus \{0\}$ is determined.  

Actually, for $\hat g\in Met^0$, a one-parameter family $\{(\varphi(x,y),\varphi_z^t(x,y))\}_{t\in R\setminus \{0\}}$ leading to $\hat g$ is firstly determined, 
next each $(\psi^t(x,y),\psi_z^t(x,y))$ with $t\neq 0$ is determined from $(\varphi(x,y),\varphi_z^t(x,y))$ with the same $t$ (\cite{bhs}, Theorems 5 and 6 in \S3.2). Furthermore, we can replace $\psi^t(x,y)$ with a new one such that $(\varphi,\varphi_z:=\varphi^t_z,\psi:=\psi^t)$ satisfies (2.1.2) on $z=0$.

Let $\psi^t(x,y,z):=\psi(x,y,z)$ and $\varphi^t(x,y,z):=\varphi(x,y,z)$ be solutions to the system (iii) and (iv) of Proposition 2.1-(3) under the above initial condition determined by $\hat g\in Met^0$ and $t\neq 0$.  
Then, $\psi^t(x,y,z)$ and $\varphi^t(x,y,z)$ also satisfy (i) and (ii) of Proposition 2.1-(3) (\cite{bhs}, Proposition 4.2 and Theorem 7). 
Therefore, $\varphi^t(x,y,z)$ gives a conformally flat metric $g^t$ by (2.1.1). 

\vspace*{2mm}
Under the facts above, we fix several notations: \\ 

\vspace*{2mm}
{\sc Definition and Notation 1}.  

\vspace*{2mm}
(1) We define two subsets $CFM$ and ${Met}$ of 3-metrics $g$ on $U=V\times I$ given by (2.1.1), respectively, as follows: 
$$CFM:=\{g |\ g \ {\rm is \ a \ conformally \ flat \ metric \ with \ the \ Guichard \ 
condition}\},  $$
$${Met}:=\{ g \ \in CFM\  | \ \varphi \ {\rm \ is \ analytic \ and \ satisfies \ } \varphi_z\varphi_{xz}\varphi_{yz}\neq 0 
 \},  $$
where $\varphi:=\varphi(x,y,z)$ is a function on $U$ determining $g$ in (2.1.1).

(2) Let ${Met}^0$ be a class of orthogonal analytic 2-metrics on $V$ with constant Gauss curvature -1, as above. That is, 
for any ${\hat g}\in Met^0$, a one-parameter family $g^t\subset Met$ is 
determined as evolutions ${\hat g}^t(z)$ of 2-metrics issuing from ${\hat g^t}(0)={\hat g}$.

\vspace*{2mm}
For a 2-metric ${\hat g}\in Met^0$ and a fixed $t\neq 0$, four functions $\varphi(x,y)$, $\varphi^t_z(x,y)$, $\psi^t(x,y)$ and $\psi^t_z(x,y)$ on $V$ have been determined such that they lead to a metric $g^t\in Met$. Among these four functions, $\psi^t(x,y)$ and $\psi^t_z(x,y)$ are determined from $(\varphi(x,y),\varphi^t_z(x,y))$, as mentioned above.  
Furthermore, $\psi^t_{zz}(x,y):=\psi_{zz}^t(x,y,0)$ and $\varphi^t_{zz}(x,y):=\varphi_{zz}^t(x,y,0)$ are also determined on $V$ from $(\varphi(x,y),\varphi^t_z(x,y))$ by (iii) and (iv) of Proposition 2.1-(3), respectively.

Now, in the definitions of (3) and (4) below, we write $g$ for $g^t$.

(3) For $\hat g\in Met^0$, we simply say that $(\varphi(x,y),\varphi_z(x,y))$ is a pair of functions on $V$ arising from $\hat g$.  
Then, this $(\varphi,\varphi_z)$ implies the above class $(\varphi,\varphi^t_z, \varphi_{zz}^t,\ \psi^t, \psi^t_z,\psi_{zz}^t)$ consisting 
of six functions on $V$ determined by $\hat g$ and some $t$.

(4) By the notation\ $g={\hat g}\in(Met^0\subset)Met$, \ we imply a metric\ $g:=g^t\in Met$\ 
determined by $\hat g\in Met^0$ and $t$.  
Furthermore, for\ $g={\hat g}\in(Met^0\subset)Met$,\ we say that $(\varphi(x,y),\varphi_z(x,y))$ 
is a pair of functions on $V$ arising from $\hat g$ and $\varphi(x,y,z)$ is a function on $U$ arising from $g$ by (2.1.1). 
Then, $\varphi(x,y,z)$ implies a pair $(\varphi^t(x,y,z),\psi^t(x,y,z))$ obtained as an evolution from $(\varphi(x,y), \psi^t(x,y), \varphi^t_z(x,y), \psi^t_z(x,y))$ with the same $t$. 
\hspace{\fill}$\Box$\\

\vspace*{2mm}
Now, let $S^3$ be the standard 3-sphere. 
Our first main result is stated as follows:  
for ${\hat g}\in Met^0$, a 5-dimensional set of (non-isometric) analytic surfaces in $S^3$ is determined such that any surface $\phi$ of the set gives rise to an evolution $\phi^z$ of surfaces in $S^3$ issuing from $\phi$ and a generic conformally flat hypersurface with the canonical principal Guichard net $g={\hat g}(\in Met^0\subset)Met$ is realized in $R^4$ via the evolution $\phi^z$.  \\

\noindent
{\bf 2.2. Covariant derivatives and curvature tensors of conformally flat hypersurfaces}

In this subsection, we summarize formulae on the covariant derivative and the curvature tensor for a generic conformally flat hypersurface in $R^4$ 
(\cite{hs3}, \S2, and \cite{su3}, \S2.1). 

Let $g\in CFM$, and $\varphi(x,y,z)$ be a function on $U$ arising from $g$ by (2.1.1). 
Let $f(x,y,z)$ be a generic conformally flat hypersurface in $R^4$ with the canonical principal Guichard net $g$. 
Then, the first fundamental form $I_f$ of $f$ is expressed as $I_f=e^{2P}g$ with a function $P=P(x,y,z)$. 
Let $\kappa_i \ (i=1,2,3)$ be the principal curvatures of $f$ corresponding to coordinate lines $x$, $y$ and $z$, respectively. 
Then, functions $\varphi$, $e^P$ and $\kappa_i$ satisfy the following relations:   
$$\cos^2\varphi=\frac{\kappa_3-\kappa_2}{\kappa_1-\kappa_2}, \hspace*{1cm} 
\sin^2\varphi=\frac{\kappa_1-\kappa_3}{\kappa_1-\kappa_2},  \eqno{(2.2.1)}$$
$$\kappa_1=e^{-P}\tan\varphi+\kappa_3, \hspace{1cm} \kappa_2=-e^{-P}\cot\varphi+\kappa_3.  \eqno{(2.2.2)}$$
In (2.2.2), we have assumed \ $(\kappa_1-\kappa_2)\cos\varphi\sin\varphi(=e^{-P})>0$, \ because we can replace a normal vector field $N$ of $f$ by $-N$ if necessary.
 
\vspace{2mm}
{\sc Remark 1}. \ \ The equations in (2.2.2) are slightly different from the ones in \cite{hs3} and \cite{su3}. 
The reason is as follows. 
For any $g\in CFM$, there is an associated family of generic conformally 
flat hypersurfaces such that each member of the family has the same principal Guichard net $g$. 
We can distinguish between these two members by constant $c^2(\in R\setminus \{0\})$ 
defined by $c^2:=e^{2P}(\kappa_1-\kappa_3)(\kappa_3-\kappa_2)$. 
If $c^2$ is same (resp. different), then two members are equivalent (resp. non-equivalent). 
We say that $f$ has the canonical principal Guichard net $g$ if $c^2=1$ holds.  
The results in \cite{hs3} and \cite{su3} are stated simultaneously for all members of one associated family. 

Through this paper, for $g\in CFM$, we only consider a hypersurface with its canonical principal Guichard net $g$, and further we shall fix $c=1$ for the sake of simplicity. However, we note that, if $c=-1$, then $N$ and (2.2.2) change to $-N$ and \ $\kappa_1=-e^{-P}\tan\varphi+\kappa_3$, \ $\kappa_2=e^{-P}\cot\varphi+\kappa_3$, respectively, \ where $\kappa_3$ is the third principal curvature for $-N$.  
For each member of one associated family, we can replace its principal Guichard net with a canonical one by modifying the original function $\varphi$ (cf. \cite{hs3}, \S2). 

From now on, we shall describe the canonical principal Guichard net for a hypersurface as only ``the Guichard net". 
\hspace{\fill}$\Box$\\

\vspace*{2mm}
Let $\nabla'$ be the standard connection of $R^4$. Let $X_{\alpha}:=\frac{e^{-P}}{\cos \varphi}f_x$, $X_{\beta}:=
\frac{e^{-P}}{\sin\varphi}f_y$ and $X_{\gamma}:=e^{-P}f_z$, which are orthonormal principal vector fields of $f$. 
For a function $h(x,y,z)$, we denote by $h_{\alpha}, \ h_{\beta}, \ h_{\gamma}$ derivatives $X_{\alpha}h, \ X_{\beta}h, \ X_{\gamma}h$, respectively. Then, the following formulae hold: 

\hspace*{10mm} $\left\{\begin{array}{lllrl}
\nabla'_{X_\alpha}X_\alpha \ = &   & \displaystyle \frac{(\kappa_1)_\beta}
{\kappa_1-\kappa_2}X_\beta 
& + \ \  \displaystyle \frac{(\kappa_1)_\gamma}{\kappa_1-\kappa_3}X_\gamma & +\ \ 
\kappa_1 N, \\[3mm]
\nabla'_{X_\beta}X_\beta \ = & \displaystyle \frac{(\kappa_2)_\alpha}{\kappa_2-\kappa_1}
X_\alpha & 
& + \ \ \displaystyle \frac{(\kappa_2)_\gamma}{\kappa_2-\kappa_3}X_\gamma & + \ \ \kappa_2 N, \\[3mm]
\nabla'_{X_\gamma}X_\gamma \ = & \displaystyle \frac{(\kappa_3)_\alpha}{\kappa_3 - \kappa_1}
X_\alpha \ \ + & \displaystyle \frac{(\kappa_3)_\beta}{\kappa_3 - \kappa_2}X_\beta 
&  & + \ \ \kappa_3 N, \\[3mm]
\nabla'_{X_\alpha}X_\beta \ = & \displaystyle -\frac{(\kappa_1)_\beta}
{\kappa_1-\kappa_2}X_\alpha, &  
&  \nabla'_{X_\alpha}X_\gamma \ = & \displaystyle -\frac{(\kappa_1)_\gamma}
{\kappa_1-\kappa_3}X_\alpha, \\[3mm]
\nabla'_{X_\beta}X_\alpha \ = & \displaystyle -\frac{(\kappa_2)_\alpha}{\kappa_2 - \kappa_1}
X_\beta, &  
 & \nabla'_{X_\beta}X_\gamma \ = & \displaystyle -\frac{(\kappa_2)_\gamma}
{\kappa_2 - \kappa_3}X_\beta, \\[3mm]
\nabla'_{X_\gamma}X_\alpha \ = & \displaystyle -\frac{(\kappa_3)_\alpha}
{\kappa_3 - \kappa_1}X_\gamma, &  
& \nabla'_{X_\gamma}X_\beta \ = & \displaystyle -\frac{(\kappa_3)_\beta}
{\kappa_3 - \kappa_2}X_\gamma. 
\end{array}\right.  \hfill{ (2.2.3)}$ \\[3mm]
Note that the covariant derivative $\nabla$ with respect to the metric $I_f$ of $f$ is also determined by (2.2.3).     

By comparing Christoffel's symbols of $I_f$ with equations in (2.2.3), we have\\[3mm]
\hspace*{10mm} $\left\{\begin{array}{lllrl}
\displaystyle \frac{(\kappa_1)_\beta}{\kappa_1-\kappa_2} \ = & -\displaystyle\frac{e^{-2P}}{\sin\varphi\cos\varphi}
(e^P\cos\varphi)_y, &  
&  \displaystyle \frac{(\kappa_1)_\gamma}{\kappa_1-\kappa_3} \ = & 
-\displaystyle\frac{e^{-2P}}{\cos\varphi}(e^P\cos\varphi)_z,\\[3mm]
\displaystyle \frac{(\kappa_2)_\alpha}{\kappa_2-\kappa_1} \ = &  
-\displaystyle\frac{e^{-2P}}{\sin\varphi\cos\varphi}(e^P\sin\varphi)_x, &  
&  \displaystyle \frac{(\kappa_2)_\gamma}{\kappa_2-\kappa_3} \ = &  -\displaystyle\frac{e^{-2P}}{\sin\varphi}(e^P\sin\varphi)_z,\\[3mm]
\displaystyle \frac{(\kappa_3)_\alpha}{\kappa_3-\kappa_1} \ = &  \displaystyle\frac{(e^{-P})_x}{\cos\varphi}, &  
& \displaystyle   \frac{(\kappa_3)_\beta}{\kappa_3-\kappa_2} \ = &  \displaystyle\frac{(e^{-P})_y}{\sin\varphi}.
\end{array}\right. \hfill{(2.2.4)}$ \\[3mm]
\\
Furthermore, the following equations also hold: 

\hspace*{12mm}$   
\begin{array}{lll}
 \varphi_x=(\kappa_1)_x \ e^P \cos^2 \varphi, \ \ \ & 
\varphi_y=(\kappa_2)_y \ e^P\sin^2 \varphi, \ \ \ &  
\varphi_z=-(\kappa_3)_z \ e^P
\end{array} . \hfill{(2.2.5)}$ \\
Note that we have 
$$\nabla'_{\partial/\partial x}X_{\beta}=\frac{(e^P\cos\varphi)_y}{e^P\sin\varphi}X_{\alpha}, \ \ \ \ 
\nabla'_{\partial/\partial y}X_{\alpha}=\frac{(e^P\sin\varphi)_x}{e^P\cos\varphi}X_{\beta}  \eqno{(2.2.6)}$$ 
by (2.2.3) and (2.2.4). 
The sectional curvatures of the metric $e^{2P}g$ are given by 
$$
\left\{
\begin{array}{l}
K(X_{\alpha}\wedge X_{\beta})=\displaystyle{\frac{e^{-P}}{\cos^2\varphi}
\left((e^{-P})_{xx}+\frac{\varphi_x}{\sin\varphi\cos\varphi}(e^{-P})_x\right)
+\frac{e^{-P}}{\sin^2\varphi}\left((e^{-P})_{yy}-\frac{\varphi_y}{\sin\varphi\cos\varphi}(e^{-P})_y\right)}  \\
\ \ \ \ \ -\displaystyle{\frac{[(e^{-P})_x]^2}{\cos^2\varphi}-\frac{[(e^{-P})_y]^2}{\sin^2\varphi}
-\frac{e^{-2P}}{\sin\varphi\cos\varphi}(\varphi_{xx}-\varphi_{yy}) }
 -\displaystyle\frac{e^{-4P}}{\sin\varphi\cos\varphi}(e^P\cos\varphi)_z(e^P\sin\varphi)_z,\\[4mm]
K(X_{\beta}\wedge X_{\gamma})  
 = \ \displaystyle{\frac{e^{-P}\varphi_x}{\sin\varphi\cos\varphi}(e^{-P})_x+\frac{e^{-P}}{\sin^2\varphi}
 \left((e^{-P})_{yy}-\varphi_y\cot\varphi(e^{-P})_y \right) }     
\hspace{\fill}{(2.2.7)}\\
 +\displaystyle{e^{-P}\left((e^{-P})_{zz}+\varphi_z\cot\varphi(e^{-P})_z \right)
 -\frac{[(e^{-P})_x]^2}{\cos^2\varphi}-\frac{[(e^{-P})_y]^2}{\sin^2\varphi}-[(e^{-P})_z]^2 }
+e^{-2P}\left(\varphi_z^2-\varphi_{zz}\cot\varphi\right),\\[4mm]
K(X_{\gamma}\wedge X_{\alpha}) 
= \ \displaystyle{\frac{e^{-P}}{\cos^2\varphi}\left((e^{-P})_{xx}+\varphi_x\tan\varphi(e^{-P})_x \right)
-\frac{e^{-P}\varphi_y}{\sin\varphi\cos\varphi}(e^{-P})_y }\\
+\displaystyle{e^{-P}\left((e^{-P})_{zz}-\varphi_z\tan\varphi(e^{-P})_z\right)
  -\frac{[(e^{-P})_x]^2}{\cos^2\varphi}-\frac{[(e^{-P})_y]^2}{\sin^2\varphi}-[(e^{-P})_z]^2}
+e^{-2P}\left(\varphi_z^2+\varphi_{zz}\tan\varphi\right).\\
\end{array} 
\right. 
$$
Note that, as far as $g$ is expressed as (2.1.1), these formulae on the sectional curvatures are satisfied 
for any metric $e^{2P}g$ even if $g$ is not conformally flat.

By the Gauss equation, we have 
$$K(X_{\alpha}\wedge X_{\beta})=\kappa_1\kappa_2, \ \ \ 
K(X_{\beta}\wedge X_{\gamma})=\kappa_2\kappa_3, \ \ \ 
K(X_{\gamma}\wedge X_{\alpha})=\kappa_3\kappa_1.  \eqno{(2.2.8)}$$
Furthermore, $\kappa_i$ satisfy the following equations: 
$$\left\{
\begin{array}{l}
(\kappa_2-\kappa_3)(\kappa_1)_x+(\kappa_1-\kappa_3)(\kappa_2)_x+(\kappa_2-\kappa_1)(\kappa_3)_x=0,\\[2mm]
(\kappa_3-\kappa_1)(\kappa_2)_y+(\kappa_2-\kappa_1)(\kappa_3)_y+(\kappa_3-\kappa_2)(\kappa_1)_y=0,\\[2mm]
(\kappa_1-\kappa_2)(\kappa_3)_z+(\kappa_3-\kappa_2)(\kappa_1)_z+(\kappa_1-\kappa_3)(\kappa_2)_z=0.
\end{array} 
\right. \eqno{(2.2.9)}
$$ \\

\vspace{2mm}

\noindent
{\bf 2.3. Determination of fundamental forms for hypersurfaces from Guichard nets}

In this subsection, we review the method of determining 
the first and the second fundamental forms for a generic conformally flat hypersurface realized in $R^4$ from a conformally flat metric $g$ with the Guichard condition (\cite{hs3}, \S2 and \S3). 
Through this subsection, let $g\in CFM$ and $\varphi(x,y,z)$ be a function on $U$ arising from $g$ by (2.1.1). 

Suppose that $f(x,y,z)$ is a conformally flat hypersurface in $R^4$ with the Guichard net $g$. 
Let $I_f=e^{2P}g$ be the first fundamental form and $\kappa_i \ (i=1,2,3)$ be the principal curvatures of $f$, as in \S2.2. 
Firstly, note that the problem of determining the first and the second fundamental forms of $f$ 
is reduced that of determining $e^P$ and $\kappa_3$ from $\varphi$, by (2.2.1) and (2.2.2).
The following proposition 2.2 gives further relation among these $\varphi,$ $e^P$ and $\kappa_3$ (cf. \cite{hs3}, Proposition 1). \\

\vspace{2mm}
{\pro 2.2}. \ \ {\it Let $g\in CFM$, and $\varphi(x,y,z)$ be a function arising from $g$ by (2.1.1). 
Let $f$ be a generic conformally flat hypersurface in $R^4$ with the Guichard net $g$. 
Let $e^P$ and $\kappa_i \ (i=1,2,3)$ be defined as above. Then, $\kappa_3$ and $e^P$ satisfy the following equations:  

$${\rm (i)} \ (\kappa_3)_x=-(e^{-P})_x\tan \varphi, \ \ \ \    
{\rm (ii)} \ (\kappa_3)_y=(e^{-P})_y\cot \varphi, \ \ \ \ 
{\rm (iii)} \ (\kappa_3)_z=-e^{-P}\varphi_z,    $$
$$
\begin{array}{l} 
 \kappa_3=\left(\tan\varphi(e^{-P})_{xx}-\varphi_x \displaystyle\frac{\cos 2\varphi}{\cos ^2\varphi} (e^{-P})_x\right)
 -\left(\cot\varphi(e^{-P})_{yy}-\varphi_y\displaystyle\frac{\cos 2\varphi}{\sin^2 \varphi}(e^{-P})_y\right) \\[2mm]
 \hspace{4cm} +\left(e^{-P}\varphi_{zz}-(e^{-P})_z\varphi_z\right),
\end{array} 
 \leqno{ \ \ \ \ ({\rm iv})}
$$
$$
\begin{array}{l} 
e^{-P}=-\left((e^{-P})_{xx}+2\varphi_x\tan\varphi(e^{-P})_x\right)-\left((e^{-P})_{yy}-2\varphi_y\cot\varphi(e^{-P})_y\right)  
\\[2mm]
 \hspace{3cm} +(e^{-P})_{zz}+2\psi_{zz}e^{-P},
\end{array} 
 \leqno{ \ \ \ \ ({\rm v})}
$$
where \ $\psi_{zz}=(\varphi_{xx}-\varphi_{yy}
 -\varphi_{zz}\cos2\varphi)/\sin2\varphi.$ }\\

\vspace*{2mm}
Now, let $e^{-P}$ be an arbitrary (positive) function on $U$. For $e^{-P}$, we define functions $I^{x,y}, \ I^{x,z} \ {\rm and } \ I^{y,z}$ on $U$ by 
$$
\left\{
\begin{array}{l}
I^{x,y}:=(e^{-P})_{xy}-(e^{-P})_y\varphi_x\cot\varphi+(e^{-P})_x\varphi_y\tan\varphi, \\

\vspace*{2mm}
I^{x,z}:=(e^{-P})_{zx}+(e^{-P})_x\varphi_z\tan\varphi-e^{-P}\varphi_{zx}\cot\varphi,\\

\vspace*{2mm}
I^{y,z}:=(e^{-P})_{zy}-(e^{-P})_y\varphi_z\cot\varphi+e^{-P}\varphi_{zy}\tan\varphi,
\end{array} \right. \eqno{(2.3.1)}
$$
respectively. Then, by (i), (ii), (iii) of Proposition 2.2 and the integrability conditions on $\kappa_3$:\ $(\kappa_3)_{xy}=(\kappa_3)_{yx}$, \ $(\kappa_3)_{xz}=(\kappa_3)_{zx}$ \ and \ $(\kappa_3)_{yz}=(\kappa_3)_{zy}$, \ 
we have the following three equations:  
$$I^{x,y}(x,y,z)=I^{x,z}(x,y,z)=I^{y,z}(x,y,z)\equiv 0, \eqno{(2.3.2)}$$
which are linear differential equations for $e^{-P}$ defined by $\varphi$. 
These equations are equivalent to the condition that the curvature tensor $R$ of a metric $e^{2P}g$ is diagonal: \ $R_{\alpha\gamma\beta\gamma}=R_{\alpha\beta\gamma\beta}=R_{\beta\alpha\gamma\alpha}=0$. 
Note that, as far as $g$ is expressed as (2.1.1), this equivalence is satisfied for any metric $e^{2P}g$ even if $g$ is not conformally flat.

The following lemma was given in (\cite{hs3}, Lemma 1). 
However, we still give a proof, since all equations in the proof are also useful for our argument.  \\

{\lem 2.1}. \ \ {\it Let $g\in CFM$, and $\varphi(x,y,z)$ be a function arising from $g$ by (2.1.1). 
Let $e^{-P}$ be an arbitrary (positive) function on $U$ satisfying Proposition 2.2-(v). 
Let $\kappa_3$ be a function defined by Proposition 2.2-(iv) from $e^{-P}$. 
Then, we have the following three equations: 
$$(\kappa_3)_x+(e^{-P})_x\tan\varphi=-\frac{1}{\sin\varphi\cos\varphi}(I^{x,y})_y+\tan\varphi(I^{x,z})_z
+\frac{\varphi_y}{\sin^2\varphi\cos^2\varphi}I^{x,y}-\frac{\varphi_z}{\cos^2\varphi}I^{x,z} ,$$
$$(\kappa_3)_y-(e^{-P})_y\cot\varphi=\frac{1}{\sin\varphi\cos\varphi}(I^{x,y})_x
-\cot\varphi(I^{y,z})_z+\frac{\varphi_x}{\sin^2\varphi\cos^2\varphi}I^{x,y}-\frac{\varphi_z}{\sin^2\varphi}I^{y,z} ,$$
$$(\kappa_3)_z+\varphi_ze^{-P}=\tan\varphi(I^{x,z})_x-\cot\varphi(I^{y,z})_y-\varphi_x\frac{\cos2\varphi}{\cos^2\varphi}I^{x,z}+
\varphi_y\frac{\cos2\varphi}{\sin^2\varphi}I^{y,z}.$$
In particular, if $e^{-P}$ satisfies (2.3.2), then (i), (ii) and (iii) in Proposition 2.2 hold. } \\

{\proof} \ Let $\kappa_3$ and $e^{-P}$ satisfy (iv) and (v) of Proposition 2.2, respectively. Then, we have 
the following equations by direct calculation: 
$$\begin{array}{l}
\displaystyle{(\kappa_3)_x+(e^{-P})_x\tan\varphi=-\frac{1}{\sin\varphi\cos\varphi}(I^{x,y})_y+\tan\varphi(I^{x,z})_z
+\frac{\varphi_y}{\sin^2\varphi\cos^2\varphi}I^{x,y}-\frac{\varphi_z}{\cos^2\varphi}I^{x,z} } \\[3mm]
\hspace*{3.5cm} \displaystyle{+2\tan\varphi\left(\psi_{xzz}+\varphi_{xzz}\cot\varphi-\frac{\varphi_z\varphi_{xz}}{\sin^2\varphi}  \right)e^{-P},}
\end{array} \eqno{(2.3.3)}$$
$$\begin{array}{l}
\displaystyle{(\kappa_3)_y-(e^{-P})_y\cot\varphi=\frac{1}{\sin\varphi\cos\varphi}(I^{x,y})_x
-\cot\varphi(I^{y,z})_z+\frac{\varphi_x}{\sin^2\varphi\cos^2\varphi}I^{x,y}-\frac{\varphi_z}{\sin^2\varphi}I^{y,z}}   \\[3mm]
\hspace*{3.5cm} \displaystyle{-2\cot\varphi \left(\psi_{yzz}-\varphi_{yzz}\tan\varphi-\frac{\varphi_z\varphi_{yz}}{\cos^2\varphi}  \right)e^{-P},}
\end{array} \eqno{(2.3.4)}$$
$$\begin{array}{l}
\displaystyle{(\kappa_3)_z+\varphi_ze^{-P}=\tan\varphi(I^{x,z})_x-\cot\varphi(I^{y,z})_y-\varphi_x\frac{\cos2\varphi}{\cos^2\varphi}I^{x,z}+
\varphi_y\frac{\cos2\varphi}{\sin^2\varphi}I^{y,z} } \\[3mm]
\displaystyle{ +2\left(\frac{(\varphi_{xx}+\varphi_{yy}+\varphi_{zz})_z}{2}+
\varphi_z\psi_{zz}-\varphi_x\varphi_{xz}\cot\varphi+\varphi_y\varphi_{yz}\tan\varphi \right)e^{-P}.}
\end{array} \eqno{(2.3.5)}$$
Note that, as far as $\varphi$, $\kappa_3$ and $e^{-P}$ satisfy (iv) and (v) of Proposition 2.2, the equations (2.3.3)$\sim$(2.3.5) hold even if a metric $g$ given by (2.1.1) from $\varphi$ is not conformally flat. 
(Then, we interpret $\psi_{zz}$ as a function defined by $\varphi$ in Proposition 2.2-(v)).

Now, in our case, all coefficients of $e^{-P}$-term in the right hand sides of (2.3.3)$\sim$(2.3.5) vanish by $g\in CFM$, that is, we have 
$$\sin2\varphi\left(\psi_{xzz}+\varphi_{xzz}\cot\varphi-\frac{\varphi_z\varphi_{xz}}{\sin^2\varphi} \right)= 
2({\rm left \ side \ of \ (ii) \ in \ Propsition \ 2.1-(2)})=0,$$
$$\sin2\varphi\left(\psi_{yzz}-\varphi_{yzz}\tan\varphi-\frac{\varphi_z\varphi_{yz}}{\cos^2\varphi} \right)= 
2({\rm left \ side \ of \ (iii) \ in \ Propsition \ 2.1-(2)})=0,$$
$${\rm coefficient \ of} \ e^{-P} \ {\rm in \ (2.3.5)}= 
2({\rm left \ side \ of \ (iv) \ in \ Propsition \ 2.1-(2)})=0.$$
In consequence, all equations in Lemma have been obtained. The last statement follows directly from these equations.  
\hspace{\fill}$\Box$\\

\vspace*{2mm}
Now, let $f$ be a generic conformally flat hypersurface in $R^4$ with the Guichard net $g$, and $e^{2P}$ be the conformal element of $I_f$. 
Then, we can rewrite the equation in Proposition 2.2-(v) as a linear evolution equation in $z$ for $e^{-P}$:   
$$(e^{-P})_{zz}=\left((e^{-P})_{xx}+2\varphi_x\tan\varphi(e^{-P})_x\right)
+\left((e^{-P})_{yy}-2\varphi_y\cot\varphi(e^{-P})_y\right)+(1-2\psi_{zz})e^{-P},   \eqno{(2.3.6)}$$
where \ $\psi_{zz}=(\varphi_{xx}-\varphi_{yy}
 -\varphi_{zz}\cos2\varphi)/\sin2\varphi.$ \ The evolution equation (2.3.6) for $e^{-P}$ will be the starting point of our argument in the next section.

Next, let us take an arbitrary (positive) function $e^{-P}$ satisfying (2.3.2) and (2.3.6) defined by $\varphi$ arising from $g\in CFM$. 
Let $\kappa_3$ be a function on $U$ given by Proposition 2.2-(iv) from $e^{-P}$, and $\kappa_1$ and $\kappa_2$ be functions given by (2.2.2).
Then, there is a function $\chi=\chi(x,y,z)$ on $U$ such that the sectional curvatures of a metric $e^{2P}g$ are expressed as  
$$K(X_{\alpha}\wedge X_{\beta})=\kappa_1\kappa_2 +\chi,\hspace{1cm}K(X_{\beta}\wedge X_{\gamma})=\kappa_2\kappa_3+\chi,
\hspace{1cm}K(X_{\gamma}\wedge X_{\alpha})=\kappa_3\kappa_1+\chi \eqno{(2.3.7)}$$ (cf. \cite{hs3}, Corollary 1 and Lemma 2).
Hence, if we choose a function $e^{-P}$ satisfying (2.3.2), (2.3.6) and \ $\chi\equiv 0$, 
then $e^{2P}g$ and $\kappa_i$ satisfy all equations of Proposition 2.2 and the Gauss condition for a hypersurface in $R^4$, as $e^{2P}g$ (resp. $\kappa_i$) is the first fundamental form (resp. principal curvatures). 

In order to investigate the condition $\chi\equiv 0$, we define a new function $\zeta$ on $U$ from $e^{-P}$ by 
$$\zeta:= 
e^{-P}\times
\left[
\begin{array}{l}
( \ (e^{-P})_{xx}+2\varphi_x\tan\varphi(e^{-P})_x \ )+( \ (e^{-P})_{yy}-2\varphi_y\cot\varphi(e^{-P})_y \ )   \\[2mm]
+( \ (e^{-P})_{zz}+\varphi_z^2e^{-P} \ )
\end{array}
\right] \eqno{(2.3.8)}$$
$$-\left( \ \displaystyle\frac{[(e^{-P})_x]^2}{\cos^2\varphi}+\frac{[(e^{-P})_y]^2}{\sin^2\varphi}+[(e^{-P})_z]^2 \ \right). $$ 
Note that $\zeta$ satisfies \ $\zeta=K(X_{\beta}\wedge X_{\gamma})\sin^2\varphi +\ K(X_{\gamma}\wedge 
X_{\alpha})\cos^2\varphi$ \ 
by (2.2.7). Then, we have that $\chi\equiv 0$ holds for a metric $e^{2P}g$ 
if and only if $e^{-P}$ satisfies \ $\zeta=\kappa_3^2$ \ (cf. \cite{hs3}, Lemma 3). 
 
\vspace*{2mm}
Now, we have the following proposition by (\cite{hs3}, Theorem 1), because we can replace the condition on $\kappa_3$ there with that on $e^{-P}$, by Proposition 2.2 and Lemma 2.1. \\ 

\vspace{2mm}
{\pro 2.3}. \ \ {\it Let $g\in CFM$, and $\varphi(x,y,z)$ be a function on $U$ arising from $g$ by (2.1.1). Let 
$e^{-P}$ be a (positive) function on $U$ satisfying (2.3.2) and (2.3.6). 
Let $\kappa_3$ and $\zeta$ be functions on $U$ defined from $\varphi$ and $e^{-P}$ by the Proposition 2.2-(iv) and (2.3.8), respectively. 
Suppose that \ $\zeta=\kappa_3^2$ \ is satisfied on $U$. 
Then, there is a hypersurface $f$ in $R^4$ with the Guichard net $g$ such that $f$ has the first fundamental form $I_f=e^{2P}g$ and 
the principal curvatures $\kappa_1, \ \kappa_2$ and $\kappa_3$, where  $\kappa_1, \ \kappa_2$ are given by (2.2.2).

More precisely, the metric $e^{2P}g$ and the three functions $\kappa_i$ satisfy the Gauss and the Codazzi conditions for a generic conformally flat hypersurface $f$ in $R^4$ as the first fundamental form and the principal curvatures of $f$, respectively. }\\

\vspace*{2mm}
The existence theorem of generic conformally flat hypersurfaces with Guichard nets $g$ ensures the existence of functions $e^{-P}$ satisfying all assumptions in Proposition 2.3.
 
At the end of this subsection, we fix some notations:  

\vspace{2mm}
{\sc Definition and Notation 2}. \ 
Let $g\in CFM$, and $\varphi(x,y,z)$ be a function on $U$ arising from $g$ by (2.1.1). Let $e^{-P}$ be a (positive) function on $U$. Let $\kappa_3$ and $\zeta$ be functions on $U$ determined from $\varphi$ and $e^{-P}$ by the Proposition 2.2-(iv) and (2.3.8), respectively. 
Then, for $\varphi$ and $e^{-P}$, we define a function $J=J(x,y,z)$ on $U$ by \ $J:=\zeta-\kappa_3^2$. 
Furthermore, for $g\in CFM$, we define a set\ $\Xi_{g}$\ of functions $e^{-P}(x,y,z)$ on $U$ by
$$\Xi_{g}:=\{ e^{-P} \ | \ e^{-P} \ {\rm satisfies \ equations} \ (2.3.2), \ (2.3.6) \ \ {\rm and} \ \ 
J\equiv 0 \ {\rm defined \ by \ \varphi} \}.   $$
\hspace{\fill}$\Box$\\

\vspace{2mm}
The set $\Xi_g$ is 5-dimensional. In fact, all non-isometric conformal transformations (homotheties and inversions) acting on $R^4$ generate a 5-dimensional set and the Guichard nets of hypersurfaces are preserved by these actions. Here, for any $q\in R^4$, an inversion $\iota_q$ with $q$ is defined by \ $\iota_q(x)=(x-q)/|x-q|^2$ \ for $x\in R^4\setminus \{q\}$, where $|x-q|$ is the Euclidean norm of $(x-q)$. \\

\newpage
\vspace*{2mm}
\noindent
{\large\bf 3. Determination of fundamental forms for hypersurfaces from 2-metrics.}

Let $g\in CFM$ and $e^{-P}(x,y,z)\in \Xi_g$. 
Then, there is a generic conformally flat hypersurface $f$ in $R^4$ which has the Guichard net $g$ and the first fundamental form $I_f=e^{2P}g$, 
by Proposition 2.3. 
In this section, we verify that, for any 2-metric ${\hat g}\in Met^0$, all functions $e^{-P}(x,y,z)\in \Xi_g$ with \ $g={\hat g}\in (Met^0\subset)Met$ are obtained as solutions to the evolution equation (2.3.6) in $z$ under initial conditions on $V$ (or on $V\times \{0\}$) determined by ${\hat g}$.

\vspace*{2mm}
Firstly, let $g$ be an arbitrary metric of $CFM$. Let $\varphi=\varphi(x,y,z)$ be a function on $U$ arising from $g$ by (2.1.1). 
Let $e^{-P}=e^{-P}(x,y,z)$ be a (positive) solution on $U$ to (2.3.6) defined by $\varphi$. 
Let $I^{x,y}, \ I^{x,z}, \ I^{y,z}$ be functions on $U$ defined from $\varphi$ and $e^{-P}$ by (2.3.1).   
Let $\kappa_3$ and $\zeta$ be functions on $U$ defined from $\varphi$ and $e^{-P}$ by Proposition 2.2-(iv) and (2.3.8), respectively. Then, as in Definition and Notation 2 and (2.2.2), functions \ $J:= \zeta-\kappa_3^2$, \ $\kappa_1:=e^{-P}\tan\varphi+\kappa_3$ \ and \ $\kappa_2:=-e^{-P}\cot\varphi+\kappa_3$ \ on $U$ are also determined.  
 
Then, we have the following two Propositions 3.1 and 3.2:\\ 

\vspace*{2mm} 
{\pro 3.1}. \ \ {\it Let $g\in CFM$, and $\varphi(x,y,z)$ be a function on $U$ arising from $g$ by (2.1.1). Let $e^{-P}(x,y,z)$ be a (positive) solution on $U$ to (2.3.6) defined by $\varphi$. 
Let $\kappa_3$ be a function defined from $e^{-P}$ by Proposition 2.2-(iv). 
Then, we have the following equations (1), (2), (3) and (4) on $U$:
$$(I^{x,y})_z=(I^{x,z})_y-\varphi_z\tan\varphi I^{x,y}+\varphi_y\tan\varphi I^{x,z}-\varphi_x\cot\varphi I^{y,z}. \leqno{(1)}$$
$$-\kappa_2\tan\varphi(I^{x,z})_z=\frac{1}{2}J_x-\frac{\kappa_3}{\sin\varphi\cos\varphi}(I^{x,y})_y
+\left(\frac{\kappa_3\varphi_y}{\cos^2\varphi}+(e^{-P})_y\right)\frac{1}{\sin^2\varphi} I^{x,y}  \leqno{(2)}$$
$$+\left(-\frac{\kappa_3\varphi_z}{\cos^2\varphi}+e^{-P}\varphi_z\tan\varphi+(e^{-P})_z\right)I^{x,z}.$$
$$\kappa_1\cot\varphi(I^{y,z})_z=\frac{1}{2}J_y+\frac{\kappa_3}{\sin\varphi\cos\varphi}(I^{x,y})_x
+\left(\frac{\kappa_3\varphi_x}{\sin^2\varphi}+(e^{-P})_x\right)\frac{1}{\cos^2\varphi}I^{x,y}  \leqno{(3)}$$
$$-\left(\frac{\kappa_3\varphi_z}{\sin^2\varphi}+e^{-P}\varphi_z\cot\varphi-(e^{-P})_z\right)I^{y,z}.$$
$$\frac{1}{2}J_z=-\kappa_2\tan\varphi(I^{x,z})_x+\kappa_1\cot\varphi(I^{y,z})_y  \leqno{(4)}$$
$$+\left(2e^{-P}\varphi_x\tan\varphi+\kappa_3\varphi_x\frac{\cos2\varphi}{\cos^2\varphi}
-\frac{(e^{-P})_x}{\cos^2\varphi}\right)I^{x,z}
-\left(2e^{-P}\varphi_y\cot\varphi+\kappa_3\varphi_y\frac{\cos2\varphi}{\sin^2\varphi}
+\frac{(e^{-P})_y}{\sin^2\varphi}  \right)I^{y,z}. $$ }\\

{\proof} \ \ These equations are obtained by direct calculation under the assumption for $g$ to be conformally flat. 
For the first equation, we have 
$$\begin{array}{l}
(I^{x,y})_z=(I^{x,z})_y-\varphi_z\tan\varphi I^{x,y}
+\varphi_y\tan\varphi I^{x,z}-\varphi_x\cot\varphi I^{y,z}   \\[3mm]
+\cot\varphi(\varphi_{xyz}+\varphi_x\varphi_{yz}\tan\varphi-\varphi_y\varphi_{xz}\cot\varphi)e^{-P}.
\end{array} \eqno{(3.1)}$$
Here, the coefficient of $e^{-P}$ vanishes by (i) of Proposition 2.1-(2).

In order to obtain the other equations, we firstly provide the following equalities: 
$$[ (e^{-P})_{xx}+2\varphi_x\tan\varphi (e^{-P})_x  ]+[ (e^{-P})_{yy}-2\varphi_y\cot\varphi (e^{-P})_y  ]   
+[  (e^{-P})_{zz}+\varphi_z^2e^{-P}  ]  \eqno{(3.2)}$$
$$=2(e^{-P})_{zz}+(2\psi_{zz}+\varphi_z^2-1)e^{-P}  \eqno{(3.3)}$$
$$=2[(e^{-P})_{xx}+2\varphi_x\tan\varphi (e^{-P})_x]+2[(e^{-P})_{yy}
-2\varphi_y\cot\varphi (e^{-P})_y]+(1+\varphi_z^2-2\psi_{zz})e^{-P},  \eqno{(3.4)}$$
where $e^{-P}\times (3.2)$ is the first term of $\zeta$ in (2.3.8). These equalities follow from (2.3.6). 
Next, we apply (3.2)$\sim$(3.4) as follows: 
$$[e^{-P}\times (3.2)]_x=(e^{-P})_x\times (3.4)+e^{-P}\times [(3.3)]_x,$$  
$$[e^{-P}\times (3.2)]_y=(e^{-P})_y\times (3.4)+e^{-P}\times [(3.3)]_y,$$ 
$$[e^{-P}\times (3.2)]_z=(e^{-P})_z\times (3.3)+e^{-P}\times [(3.4)]_z.$$ 

Now, for the derivative $\zeta_x$ of $\zeta$, by Proposition 2.2-(iv) and (2.3.1) we have 
$$\begin{array}{l}
\displaystyle{\zeta_x=(\kappa_3^2)_x-2\kappa_3((\kappa_3)_x+(e^{-P})_x\tan\varphi)+2e^{-P}(I^{x,z})_z
-\frac{2(e^{-P})_y}{\sin^2\varphi}I^{x,y}}   \\[3mm]
\displaystyle{-2((e^{-P})_z+e^{-P}\varphi_z\tan\varphi)I^{x,z}
+2\left(\psi_{xzz}+\varphi_{xzz}\cot\varphi-\frac{\varphi_z\varphi_{xz}}{\sin^2\varphi}\right)e^{-2P}.}
\end{array}\eqno{(3.5)}$$
The coefficient of $e^{-2P}$ in (3.5) vanishes by the assumption for $g$ to be conformally flat, as in Proof of Lemma 2.1. 
Hence, we obtain the second equation by (2.3.3), (3.5) and 
$e^{-P}-\kappa_3\tan\varphi=-\kappa_2\tan\varphi$. 

Similarly, by the Proposition 2.2-(iv) and (2.3.1), we have 
$$\begin{array}{l}
\displaystyle{\zeta_y=(\kappa_3^2)_y-2\kappa_3((\kappa_3)_y-(e^{-P})_y\cot\varphi)+2e^{-P}(I^{y,z})_z-\frac{2}{\cos^2\varphi}(e^{-P})_xI^{x,y} }  \\[3mm]
\displaystyle{-2[(e^{-P})_z-e^{-P}\varphi_z\cot\varphi]I^{y,z}+2\left(\psi_{yzz}-\varphi_{yzz}\tan\varphi-\frac{\varphi_z\varphi_{yz}}
{\cos^2\varphi}\right)e^{-2P}. }
\end{array}\eqno{(3.6)}$$
Since the coefficient of $e^{-2P}$ in (3.6) vanishes, we obtain the third equation by (2.3.4), (3.6) and 
$e^{-P}+\kappa_3\cot\varphi=\kappa_1\cot\varphi$.

We also have 
$$
\hspace*{3cm}\displaystyle{\zeta_z=(\kappa_3^2)_z-2\kappa_3((\kappa_3)_z+e^{-P}\varphi_z)+2e^{-P}((I^{x,z})_x+(I^{y,z})_y)}  $$
$$\displaystyle{+2\left(2e^{-P}\varphi_x\tan\varphi-\frac{(e^{-P})_x}{\cos^2\varphi}\right)I^{x,z}
-2\left(2e^{-P}\varphi_y\cot\varphi+\frac{(e^{-P})_y}{\sin^2\varphi}\right)I^{y,z}} \eqno{(3.7)}$$
$$\displaystyle{+\frac{4\cos2\varphi}{\sin2\varphi}
\left(\frac{(\varphi_{xx}+\varphi_{yy}+\varphi_{zz})_z}{2}+\frac{(\varphi_{xx}-\varphi_{yy})-\varphi_{zz}\cos\varphi}{\sin2\varphi}
\varphi_z-\varphi_x\varphi_{xz}\cot\varphi+\varphi_y\varphi_{yz}\tan\varphi \right)e^{-2P}.}$$
Since the coefficient of $e^{-2P}$ in (3.7) vanishes, we obtain the last equation by (2.3.5), (3.7), \ 
$e^{-P}-\kappa_3\tan\varphi=-\kappa_2\tan\varphi$ \ and \ $e^{-P}+\kappa_3\cot\varphi=\kappa_1\cot\varphi$. 

Note that, as far as $\varphi$, $e^{-P}$ and $\kappa_3$ satisfy (2.3.6) and Proposition 2.2-(iv), the equations (3.1) and (3.5)$\sim$(3.7) hold even if a metric $g$ determined by (2.1.1) from $\varphi$ is not conformally flat.  
\hspace{\fill}$\Box$ \\

\vspace*{2mm}
{\sc Remark 2.} \ We give another proof of Theorem due to Canevari and Tojeiro (\cite{ct}), in relation to our argument. In the following statement, a metric $g$ is defined by (2.1.1) from a function $\varphi(x,y,z)$, but we do not assume that $g$ is conformally flat:\\

{\it ``A generic hypersurface $f$ in $R^4$ is conformally flat, if and only if 
$f$ satisfies the following conditions (1) and (2): 

(1) There is a principal curvature line coordinate system $(x,y,z)$ of $f$ such that the first fundamental form $I_f$ is expressed as \ $I_f=e^{2P}g:=e^{2P}(\cos^2\varphi dx^2+\sin^2\varphi dy^2+dz^2)$.\ 

(2) The principal curvatures $\kappa_i \ (i=1,2,3)$ corresponding to  the coordinate lines $x$, $y$ and $z$, respectively, satisfy the equations \ 
$\kappa_1=e^{-P}\tan\varphi+\kappa_3,$ \ and \ $\kappa_2=-e^{-P}\cot\varphi+\kappa_3.$" }\\

Suppose that a generic hypersurface $f$ in $R^4$ is conformally flat. Then, $f$ satisfies (1) and (2), as mentioned in \S2.2. Now, we verify the converse. Let a generic hypersurface $f(x,y,z)$ in $R^4$ satisfy (1) and (2). 
Then, three sectional curvatures $K(X_{\alpha}\wedge X_{\beta})$, $K(X_{\beta}\wedge X_{\gamma})$ and $K(X_{\gamma}\wedge X_{\alpha})$ of the metric $e^{2P}g$ are given by (2.2.7), and \ $K(X_{\alpha}\wedge X_{\beta})=\kappa_1\kappa_2$, \ $K(X_{\beta}\wedge X_{\gamma})=\kappa_2\kappa_3$ \ and 
$K(X_{\gamma}\wedge X_{\alpha})=\kappa_3\kappa_1$ \ are satisfied, since $(x,y,z)$ is a principal curvature line coordinate system. By these facts and the condition (2), $\kappa_3$ and $e^{-P}$ satisfy (iv) and (v) (or (2.3.6)) of Proposition 2.2, respectively, of which proofs are same as in the proposition.

Since \ $\kappa_3=\kappa_1\cos^2\varphi+\kappa_2\sin^2\varphi$ \ holds by (2), we have \ $J=\zeta-\kappa_3^2= 0$ on $U$. \ 
Furthermore, we have \ $I^{x,y}(x,y,z)=I^{x,z}(x,y,z)=I^{y,z}(x,y,z)= 0$ \ on $U$, since the curvature tensor of $I_f=e^{2P}g$ is diagonal with respect to the principal curvature line coordinate system.

Thus, by (2.3.3)$\sim$(2.3.5), (3.1) and (3.5)$\sim$(3.7), we have whether $g$ is conformally flat or one of the following three equations is satisfied: $$e^{-P}=2\kappa_3\tan\varphi, \ \ \ e^{-P}= -2\kappa_3\cot\varphi, \ \ \ e^{-P}=\kappa_3\tan2\varphi.$$ 
However, all these three equations are not conformally invariant. In consequence, we have verified that $g$ is conformally flat. 
\hspace{\fill}$\Box$\\

\vspace*{2mm} 
By Proposition 3.1, we have the following proposition: \\

\vspace*{2mm}
{\pro 3.2}. \ \ {\it  In Proposition 3.1, let $g$ be an analytic metric of $CFM$ and $e^{-P}(x,y,z)$ be an analytic solution to (2.3.6).   
Suppose that $e^{-P}(x,y,z)$ satisfies \ $(\kappa_1\kappa_2)(x,y,0)\neq 0$ \ and \ $I^{x,y}(x,y,0)=I^{x,z}(x,y,0)=I^{y,z}(x,y,0)=J(x,y,0)= 0$ \ for $(x,y)\in V$.  
Then, we have $e^{-P}\in \Xi_g$, by choosing a suitable subdomain $V'$ of $V$ if necessary. That is, the equations \ $I^{x,y}(x,y,z)=I^{x,z}(x,y,z)=I^{y,z}(x,y,z)=J(x,y,z)= 0 $ \ hold on $U$. 
In particular, \ 
$(\kappa_3)_x=-(e^{-P})_x\tan\varphi$, \ $(\kappa_3)_y=(e^{-P})_y\cot\varphi$ \ and \ $(\kappa_3)_z=-e^{-P}\varphi_3$ \ also hold on $U$. } \\

\vspace*{2mm}
{\proof} \ Since the proof relies on the Cauchy-Kovalevskaya theorem for analytic evolution equations, we have assumed that $g$ and $e^{-P}$ are analytic. 

Now, by choosing a suitable subdomain $V'$ of $V$ if necessary, we can assume that \ $(\kappa_1\kappa_2)(x,y,z)\neq 0$ \ holds on $U=V\times I$ from the assumption \ $(\kappa_1\kappa_2)(x,y,0)\neq 0$ \ on $V$. 
Then, the equations for $(I^{x,y}, I^{x,z},I^{y,z},J)$ in Proposition 3.1 are a system of evolution equations in $z$. 
Hence, under the initial condition $(I^{x,y}=0, \ I^{x,z}=0, \ I^{y,z}=0, \ J=0)$ on $V\times \{0\}$, 
we have $(I^{x,y}\equiv 0, \ I^{x,z}\equiv 0, \ I^{y,z}\equiv 0, \ J\equiv 0)$ on $U=V\times I$ by the uniqueness of solutions. 

The last three equations follow from Lemma 2.1.
\hspace{\fill}$\Box$ \\

Now, let \ ${\hat g}\in Met^0$. \ Let $(\varphi(x,y),\varphi_z(x,y))$ be functions on $V$ arising from ${\hat g}$ and $\varphi(x,y,z)$ be a function on $U$ arising from $g={\hat g}\in (Met^0\subset)Met$ by (2.1.1). 
Then, by Definition and Notation 1 in \S2.1, the pair $(\varphi(x,y),\varphi_z(x,y))$ implies a class $(\varphi(x,y),\varphi_z^t(x,y),\varphi_{zz}^t(x,y),\psi^t(x,y),\psi_z^t(x,y),\psi_{zz}^t(x,y))$ 
of functions on $V$ determined by ${\hat g}$ and $t\neq 0$, and $\varphi(x,y,z)$ implies a pair $(\varphi^t(x,y,z),\psi^t(x,y,z))$ of functions on $U$ obtained as evolution from $(\varphi(x,y),\psi^t(x,y),\varphi_z^t(x,y),\psi_z^t(x,y))$ with the same $t$. 
In particular, note that all functions above are analytic, since $\hat g$ is analytic by the definition. 

The assumption: \ $\kappa_1\kappa_2\neq 0$ \ and \ $I^{x,y}=I^{x,z}=I^{y,z}=J=0$ \ hold on $V\times\{0\}$, in Proposition 3.2 imposes constraints on two functions $e^{-P}(x,y,0)$ and $(e^{-P})_z(x,y,0)$ on $V$. We clarify the fact by the following Definition and Notation 3: \\

\vspace*{2mm}
{\sc Definition and Notation 3}. \ Let ${\hat g}\in Met^0$, and $(\varphi(x,y),\varphi_z(x,y))$ be functions on $V$ 
arising from ${\hat g}$. 
Then, the following functions $\bar P$, $\bar P_z$, $\bar\kappa_i$, $\bar\zeta$ on $V$ are restrictions to $V\times\{0\}$ of functions $P$, $P_z$, $\kappa_i$, $\zeta$ on $U$, respectively: 

\vspace*{2mm}
(1) Let $\bar P(x,y)$ and $\bar P_z(x,y)$ be arbitrary analytic functions on $V$ independent of each other. 
Then, we denote $(e^{-\bar P})_z(x,y):=-(e^{-\bar P}\bar P_z)(x,y)$.\\[-3mm] 

(2) For $e^{-\bar P}(x,y)$ and $(e^{-\bar P})_z(x,y)$, \ let $\bar{\kappa}_3(x,y)$ and $(e^{-\bar P})_{zz}(x,y)$ be functions on $V$ defined by 
$$\bar{\kappa}_3:=\left(\tan\varphi(e^{-\bar P})_{xx}-\varphi_x\frac{\cos2\varphi}{\cos^2\varphi}(e^{-\bar P})_x\right)
-\left(\cot\varphi (e^{-\bar P})_{yy}-\varphi_y\frac{\cos2\varphi}{\sin^2\varphi} (e^{-\bar P})_y\right)$$
$$+ \ [e^{-\bar P}\varphi_{zz}-(e^{-\bar P})_z\varphi_z],$$
$$ (e^{-\bar P})_{zz}:=((e^{-\bar P})_{xx}+2\varphi_x\tan\varphi~(e^{-\bar P})_x)+((e^{-\bar P})_{yy}
-2\varphi_y\cot\varphi~(e^{-\bar P})_y)+(1-2\psi_{zz})e^{-\bar P},$$
where \ $\psi_{zz}(x,y):=[(\varphi_{xx}-\varphi_{yy}-\varphi_{zz}\cos2\varphi)/\sin2\varphi](x,y)$. \\[-3mm]

(3) For $e^{-\bar P}(x,y)$ and $(e^{-\bar P})_z(x,y)$, \ let ${\bar \kappa}_1(x,y)$ and ${\bar\kappa}_2(x,y)$ be functions on $V$ defined 
by \ ${\bar \kappa}_1:=e^{-\bar P}\tan\varphi+{\bar \kappa}_3, \ \ 
{\bar\kappa}_2:=-e^{-\bar P}\cot\varphi+{\bar \kappa}_3.$ \\[-3mm]

(4) For $e^{-\bar P}(x,y)$ and $(e^{-\bar P})_z(x,y)$, \ let $\bar\zeta(x,y)$ be a function on $V$ defined by 
$$\bar\zeta:= 
e^{-\bar P}\times
\left[
\begin{array}{l}
( \ (e^{-\bar P})_{xx}+2\varphi_x\tan\varphi~ (e^{-\bar P})_x \ )+( \ (e^{-\bar P})_{yy}-2\varphi_y\cot\varphi~(e^{-\bar P})_y \ )   \\[2mm]
+( \ (e^{-\bar P})_{zz}+\varphi_z^2 e^{-\bar P} \ )
\end{array}
\right] $$
$$-\left( \ \displaystyle\frac{((e^{-\bar P})_x)^2}{\cos^2\varphi}+\frac{((e^{-\bar P})_y)^2}{\sin^2\varphi}+((e^{-\bar P})_z)^2 \ \right). $$
\hspace{\fill}$\Box$\\ 

\vspace*{2mm}
Now, let ${\hat g}\in Met^0$, and $(\varphi(x,y),\varphi_z(x,y))$ be functions arising from $\hat g$. For analytic functions $\bar P(x,y)$ and $\bar P_z(x,y)$ on $V$, we consider the following four equations:   
$$ \ \ (e^{-\bar P})_{xy}-(e^{-\bar P})_y\varphi_x\cot\varphi+(e^{-\bar P})_x\varphi_y\tan\varphi=0,  \eqno{(3.8)}$$
$$((e^{-\bar P})_z)_x+(e^{-\bar P})_x\varphi_z\tan\varphi- e^{-\bar P}\varphi_{zx}\cot\varphi=0, \eqno{(3.9)}$$
$$((e^{-\bar P})_z)_y-(e^{-\bar P})_y\varphi_z\cot\varphi+e^{-\bar P}\varphi_{zy}\tan\varphi=0,  \eqno{(3.10)}$$
$$\bar{\kappa}_3^2(x,y)=\bar\zeta(x,y), \eqno{(3.11)}$$ 
where \ $\varphi_{zx}(x,y):=(\varphi_z)_x(x,y)$ and $\varphi_{zy}(x,y):=(\varphi_z)_y(x,y)$. Note that (3.8) is the integrability condition on $(e^{-\bar P})_z$: \ $((e^{-\bar P})_z)_{xy}=((e^{-\bar P})_z)_{yx}$, \ in this case. 

Then, we have the following Theorem. \\

\vspace*{2mm}
{\thm 1}. \ \ {\it Let $\hat{g}\in Met^0$, and $(\varphi(x,y),\varphi_z(x,y))$ be functions on $V$ arising from ${\hat g}$. 
Suppose that analytic functions $\bar P(x,y)$ and $\bar P_z(x,y)$ on $V$ satisfy \ $(\bar\kappa_1\bar\kappa_2)(x,y)\neq 0$ \ and (3.8)$\sim$(3.11) defined by $(\varphi(x,y),\varphi_z(x,y))$. 
Let $e^{-P}(x,y,z)$ be a solution on $U=V\times I$ to (2.3.6) defined by $\varphi(x,y,z)$ arising from $g={\hat g}\in (Met^0\subset) Met$, 
under the initial condition $e^{-P}(x,y,0)=e^{-\bar P}(x,y)$ and $(e^{-P})_z(x,y,0)=(e^{-\bar P})_z(x,y)$. 
Let $\kappa_3(x,y,z)$ (resp. $\kappa_i(x,y,z), \ i=1,2$) be a function defined from $e^{-P}$ by Proposition 2.2-(iv) (resp. by (2.2.2)). 
Then, there is a generic conformally flat hypersurface $f(x,y,z)$ of $U$ into $R^4$ such that $f$ has the Guichard net $g$, the first fundamental form $I_f=e^{2P}g$ and principal curvatures $\kappa_i$, by replacing $V$ with a suitable subdomain $V'$ if necessary. 

In particular, \ $(\kappa_3)_x=-(e^{-P})_x\tan\varphi$, \ $(\kappa_3)_y=(e^{-P})_y\cot\varphi$, \ $(\kappa_3)_z=-e^{-P}\varphi_3$ \ hold on $U$.}\\

\vspace*{2mm}
{\sc Remark 3}. \ (1) In Theorem 1, we choose the subdomain $V'$ of $V$ such that $\kappa_1\kappa_2\neq 0$ holds on $V'\times I$.

(2) Let $f(x,y,z)$ be a generic conformally flat hypersurface in $R^4$ determined by Theorem 1, and $N$ be a unit normal vector field of $f$. Then, \ $\phi(x,y):=N(x,y,0)$ \ for $(x,y)\in V$ \ is a surface in $S^3$, further 
$\phi^z(x,y):=N(x,y,z)$, $z\in I$, is an evolution of surfaces in $S^3$ issuing from $\phi$. In the next section, we relate $(\bar P(x,y),\bar P_z(x,y))$ to both $\phi$ and $\phi^z$ more directly. 
\hspace{\fill}$\Box$\\

\vspace*{2mm}
{\sc Proofs of Theorem 1 and Explanation of Remark 3}. \ Let \ $\hat g\in Met^0$ \ and \ $g={\hat g}\in (Met^0\subset) Met$. \ 
Let $(\varphi(x,y),\varphi_z(x,y))$, $(\bar P(x,y),\bar P_z(x,y))$ and $\varphi(x,y,z)$, $e^{-P}(x,y,z)$, $\kappa_i(x,y,z)$ be functions given in the theorem. 
Let $\zeta(x,y,z)$ be a function defined by (2.3.8) from $\varphi(x,y,z)$ and $e^{-P}(x,y,z)$. 
Then, we have \ $\varphi(x,y,0)=\varphi(x,y)$, \ $\varphi_z(x,y,0)=\varphi_z(x,y)$ \ and \ $\varphi_{zz}(x,y,0)=\varphi_{zz}(x,y)$ \ by Definition and Notation 1, and have \ $e^{-P}(x,y,0)=e^{-\bar P}(x,y)$, \  $(e^{-P})_z(x,y,0)=(e^{-\bar P})_z(x,y)$, \  $(e^{-P})_{zz}(x,y,0)=(e^{-\bar P})_{zz}(x,y)$, \ $\kappa_i(x,y,0)=\bar\kappa_i(x,y)$ \ and \ $\zeta(x,y,0)=\bar\zeta(x,y)$ \ by Definition and Notation 3. 

Therefore, for $I^{x,y}(x,y,z)$, $I^{x,z}(x,y,z)$, $I^{y,z}(x,y,z)$ defined by (2.3.1) and $J(x,y,z)$ in Definition and Notation 2, we have \ $I^{x,y}=I^{x,z}=I^{y,z}=J=0$ \ on $V\times\{0\}$, by the assumptions (3.8)$\sim$(3.11) for $({\bar P}, {\bar P}_z)$. 
Furthermore, we have \ $(\kappa_1\kappa_2)(x,y,z)\neq 0$ \ on $V\times I$ from the assumption \ $(\bar\kappa_1\bar\kappa_2)(x,y)\neq 0$ \ on $V$, by choosing a suitable subdomain $V'$ of $V$ if necessary. Hence, $e^{-P}(x,y,z)$ belongs to $\Xi_g$ \ and the last three equations also hold by Proposition 3.2.

In consequence, Theorem 1 is verified by Proposition 2.3. 
For Remark 3-(2), each map $\phi^z(x,y)$ with $z\in I$ defines a surface in $S^3$, since $N_x=-\kappa_1f_x$, $N_y=-\kappa_2f_y$ and $\kappa_1\kappa_2\neq 0$ hold on $U$.   \hspace{\fill}$\Box$\\

Next, we study the assumption $(\bar\kappa_1\bar\kappa_2)(x,y)\neq 0$ for $(\bar P(x,y),\bar P_z(x,y))$ in Theorem 1. The function $\bar\kappa_1\bar\kappa_2$ is analytic on $V$. Hence, if it does not vanish identically, then we can assume $\bar\kappa_1\bar\kappa_2\neq 0$ on $V$ by choosing a suitable subdomain $V'$ of $V$ if necessary. \\

\vspace*{2mm}
{\pro 3.3}. \ \ {\it Let $(\varphi(x,y),\varphi_z(x,y))$ be functions on $V$ arising from $\hat g\in Met^0$. Suppose that a pair $({\bar P}(x,y),{\bar P}_z(x,y))$ on $V$ is a solution to (3.8)$\sim$(3.11) defined by $(\varphi(x,y),\varphi_z(x,y))$ and satisfies $(\bar\kappa_1\bar\kappa_2)(x,y)\equiv 0$ on $V$. Then, $\hat g$ is a 2-metric given in (\cite{bhs}, Theorem 5). However, for any Guichard net $g$ arising from such a $\hat g$, there is no \ $e^{-P}(x,y,z)\in \Xi_g$ \ such that \ $e^{-P}(x,y,0)=e^{-\bar P}(x,y)$ \ and \ $(e^{-P})_z(x,y,0)=(e^{-\bar P})_z(x,y)$ \ hold.}

\vspace*{2mm}
{\sc Proof}. \ Any 2-metric $\hat g$ in (\cite{bhs}, Theorem 5) belongs to $Met^0$, and it is characterized by the condition that $\varphi(x,y)$ arising from $\hat g$ is a one-variable function. 
Hence, for $(\varphi(x,y),\varphi_z(x,y))$ arising from $\hat g\in Met^0$, we show the following fact: if a solution $({\bar P}(x,y),{\bar P}_z(x,y))$ to (3.8)$\sim$(3.11) defined by $(\varphi(x,y),\varphi_z(x,y))$ satisfies either $\bar\kappa_1\equiv 0$ or $\bar\kappa_2\equiv 0$ on $V$, then $\varphi(x,y)$ is a one-variable function.

Let $\hat g$, $(\varphi(x,y),\varphi_z(x,y))$ and $(\bar P(x,y),\bar P_z(x,y))$ be as above. 
Let $g=\hat g\in(Met^0\subset)Met$, and $\varphi(x,y,z)$ be a function arising from $g$ by (2.1.1). Let $e^{-P}(x,y,z)$ be a solution to (2.3.6) defined by $\varphi(x,y,z)$ under the initial data $e^{-\bar P}(x,y)$ and $(e^{-\bar P})_z(x,y)$ on $V\times\{0\}$. Let $\kappa_3(x,y,z)$ (resp. $\kappa_i(x,y,z), \ i=1,2$) be a function defined by Proposition 2.2-(iv) (resp. by (2.2.2)) from $e^{-P}$. Let $I^{x,y}(x,y,z)$, $I^{x,z}(x,y,z)$, $I^{y,z}(x,y,z)$ be functions defined by (2.3.1) and $J(x,y,z)$ be a function in Definition and Notation 2. Then, we have $\kappa_i(x,y,0)=\bar\kappa_i(x,y)$ for $i=1,2,3$, by Definition and Notation 3. The equations (3.8)$\sim$(3.11) for $(e^{-\bar P}(x,y),(e^{-\bar P})_z(x,y))$ are equivalent to \ 
$I^{x,y}(x,y,0)=I^{x,z}(x,y,0)=I^{y,z}(x,y,0)=J(x,y,0)\equiv 0$ \ on $V$. 
Hence, we have \ $(I^{x,y})_z(x,y,0)=[\kappa_2(I^{x,z})_z](x,y,0)=[\kappa_1(I^{y,z})_z](x,y,0)=J_z(x,y,0)\equiv 0$ \ on $V$ by Proposition 3.1. 

Now, suppose that \ $\kappa_1(x,y,0)(=[e^{-P}\tan\varphi+\kappa_3](x,y,0))\equiv 0$ \ holds on $V$. 
Then, we have \ $(I^{x,z})_z(x,y,0)\equiv 0$ \ by $\kappa_2(x,y,0)\neq 0$ on $V$, hence we firstly have \ $[(\kappa_3)_x+(e^{-P})_x\tan\varphi](x,y,0)\equiv 0$ \ on $V$ by Lemma 2.1. 
On the other hand, since $\kappa_1(x,y,0)\equiv 0$, we directly have \ $(\kappa_3)_x(x,y,0)=-[(e^{-P})_x\tan\varphi+e^{-P}\varphi_x/\cos^2\varphi](x,y,0)$ \ on $V$. From these two equations and $e^{-\bar P}\neq 0$, we have $\varphi_x(x,y,0)\equiv 0$ on $V$, which shows that $\varphi(x,y)$ is a one-variable function for $y$. In the case $\kappa_2(x,y,0)\equiv 0$, we have that $\varphi(x,y)$ is a one-variable function for $x$, in the same way. 

For the second statement, we shall verify it in Example 1 below. 
\hspace{\fill}$\Box$\\

\vspace*{2mm}
Let $\hat g$ be any 2-metric in (\cite{bhs}, Theorem 5) and $(\varphi(x,y),\varphi_z(x,y))$ be functions arising from $\hat g$. 
Then, we have that all pairs $({\bar P}(x,y),{\bar P}_z(x,y))$ on $V$ satisfying $\bar\kappa_1\bar\kappa_2\neq 0$ and (3.8)$\sim$(3.11) defined by $(\varphi(x,y),\varphi_z(x,y))$ generate a 5-dimensional set, by Theorem 1 and Proposition 3.3. Furthermore, for $(\varphi(x,y),\varphi_z(x,y))$ arising from the other $\hat g\in Met^0$, all pairs $({\bar P}(x,y),{\bar P}_z(x,y))$ on $V$ satisfying (3.8)$\sim$(3.11) defined by $(\varphi(x,y),\varphi_z(x,y))$ generate a 5-dimensional set, since such pairs necessarily satisfy $\bar\kappa_1\bar\kappa_2\neq 0$ on $V$ by choosing a suitable subdomain $V'$ of $V$ if necessary. \\

\vspace{2mm}
{\bf Examples.} \ Through examples, we study a method of solving the equations (3.8)$\sim$(3.11) defined by $(\varphi(x,y),\varphi_z(x,y))$ arising from $\hat g\in Met^0$. We firstly give some remarks. 

When we solve (3.8)$\sim$(3.11) for $({\bar P}(x,y), {\bar P}_z(x,y))$, we can replace (3.8) with two equations \ $(\bar\kappa_3)_x=-(e^{-\bar P})_x\tan\varphi$ \ and \ $(\bar\kappa_3)_y=(e^{-\bar P})_y\cot\varphi$. \ 
In fact, (3.8) follows from these two equations as the integrability condition on $\bar\kappa_3(x,y)$. 
On the other hand, suppose that $({\bar P}(x,y), {\bar P}_z(x,y))$ further satisfies $\bar\kappa_1\bar\kappa_2\neq 0$ in addition to (3.8)$\sim$(3.11). Then, by Theorem 1, a generic conformally flat hypersurface $f$ in $R^4$ is determined from $(\varphi(x,y),\varphi_z(x,y))$ and $({\bar P}(x,y), {\bar P}_z(x,y))$, and these two equations hold. 
This additional condition $\bar\kappa_1\bar\kappa_2\neq 0$ is necessary only for 2-metrics $\hat g$ in Example 1, where we treat every 2-metric $\hat g$ in (\cite{bhs}, Theorem 5). 
In Example 1, we shall show that there is not any $(\bar P(x,y),\bar P_z(x,y))$ satisfying $\bar\kappa_1\bar\kappa_2\equiv 0$ in solutions to \ $(\bar\kappa_3)_x=-(e^{-\bar P})_x\tan\varphi$, \ $(\bar\kappa_3)_y=(e^{-\bar P})_y\cot\varphi$ \ and  (3.9)$\sim$(3.11) defined by $\hat g$. This implies the following fact by Proposition 2.2: if $(\bar P(x,y),\bar P_z(x,y))$ is a solution to (3.8)$\sim$(3.11) defined by $\hat g$ and satisfies $\bar\kappa_1\bar\kappa_2\equiv 0$, then there is no $e^{-P}(x,y,z)\in \Xi_g$ such that $e^{-P}(x,y,0)=e^{-\bar P}(x,y)$ and $(e^{-P})_z(x,y,0)=(e^{-\bar P})_z(x,y)$ hold, where $g=\hat g\in (Met^0\subset)Met$.

For the sake of simplicity, we study only $(\varphi(x,y),\varphi_z(x,y)):=(\varphi(x,y),\varphi_z^1(x,y))$, that is, the case of $t=1$ for \ $(\varphi(x,y),\varphi^t_z(x,y))$ in Definition and Notation 1-(2) and (3) of \S2.1. 

For functions $\varphi_z(x,y)$ and $(\bar P(x,y),{\bar P}_z(x,y))$ on $V$, we denote their derivatives by $\varphi_{zx}(x,y):=(\varphi_z)_x(x,y)$, \ 
$(e^{-\bar P})_{zy}(x,y):=-(e^{-\bar P}\bar P_{z})_y(x,y)$ \ and so on, as usual. 

Although $e^{-\bar P}(x,y)$ is a positive function as its notation shows, the sign of $e^{-\bar P}$ is not important, since we can replace $e^{-\bar P}$ with $-e^{-\bar P}$.

\vspace*{2mm}
{\sc Example 1.} \ Let $\hat g\in Met^0$ be a 2-metric given in (\cite{bhs}, Theorem 5). 
Such a 2-metric $\hat g$ gives rise to functions 
$\varphi(x,y)$, $\varphi_z(x,y)$ and $L\psi(x,y):=(\psi_{xx}-\psi_{yy})(x,y)$ 
such that \ $\cos\varphi(x,y)=1/\sqrt{1+e^{2\rho(y)}},$ \ $\sin\varphi(x,y)=e^{\rho(y)}/\sqrt{1+e^{2\rho(y)}}$ \ and 
$$\varphi_z(x,y)=\sigma(x)\sin\varphi, \ \ \ 
L\psi(x,y)=(1/2)[\sigma^2(x)-\varphi_y^2/\cos^2\varphi]-\varphi_{yy}\tan\varphi,$$  
respectively, where $\rho(y)$ and $\sigma(x)(>0)$ are (strongly) increasing (or decreasing) functions of one-variable. For the sake of simplicity, we assume that the both domains of $\rho(y)$ and $\sigma(x)$ include $0\in R$.

\vspace*{2mm} 
Let $X_1(x)$ be a solution to the linear ordinary differential equation  
$$\sigma X_1^{'''}-\sigma' X_1^{''}+\sigma(1+\sigma^2) X_1'-\sigma' X_1=0. \eqno{(3.12)}$$
The equation (3.12) is rewritten as 
$$X_1^{''}+(1+\sigma^2)X_1-\sigma(\int_0^x\sigma' X_1dx+c_0)=0  \eqno{(3.13)}$$
with a constant $c_0$. Let us denote \ $X_2(x):=\int_0^x\sigma' X_1dx+c_0$, \ which appeared in (3.13). 
Let $Y_1(y)$ and $Y_2(y)$  be functions of one-variable satisfying the equations  
$$Y_2=-e^{-\rho}Y_1^{''}+\rho'e^{-\rho}\cos2\varphi~ Y_1'-[(\rho^{''}+\rho'^2\cos^2\varphi)\sin\varphi\cos\varphi-e^{\rho}]Y_1, \eqno{(3.14)}$$
$$Y_2'=(e^{\rho}Y_1)'+e^{-\rho}Y_1'.  \eqno{(3.15)}$$
Note that $Y_1(y)$ is a solution to a linear ordinary differential equation of third order, and 
the function $Y_2(y)$ is expressed as \ $Y_2(y)=(e^{\rho}Y_1)(y)+\int_0^y e^{-\rho}Y_1'dy+c_1$ \ with a constant $c_1$.

Then, the functions $e^{-\bar P}(x,y),$ $\bar\kappa_3(x,y)$ and $(e^{-\bar P})_z(x,y)$ are determined as follows: 
$$\begin{array}{l}
e^{-\bar P}:=\cos\varphi~X_1+Y_1, \ \ \ \ \bar\kappa_3:=-\sin\varphi~X_1-e^{\rho}Y_1+Y_2, \\[3mm]
\ \ \ \ (e^{-\bar P})_z:=-\sigma\sin^2\varphi~X_1+X_2+\sigma\cos\varphi~Y_1
\end{array} \eqno{(3.16)}
$$
from \ $(\bar\kappa_3)_x=-(e^{-\bar P})_x\tan\varphi$, \ $(\bar\kappa_3)_y=(e^{-\bar P})_y\cot\varphi$, \ (3.9), (3.10) and the definition of $\bar\kappa_3$, where $(X_1(x),Y_1(y))$ is an arbitrary pair of solutions to the above equations except for $X_1=Y_1\equiv 0$. 

For functions in (3.16), we have \ $\bar\kappa_3^2-\bar\zeta=G(x)+H(y),$ \ 
where
$$G(x):=\left[(\sigma X_1-X_2)^2+X_1^2+(X_1')^2\right](x),   \eqno{(3.17)}$$
$$H(y):=\left[Y_2^2+\left(e^{-\rho}\sqrt{1+e^{2\rho}}~Y_1'+\frac{\rho'e^{\rho}}{\sqrt{1+e^{2\rho}}}Y_1\right)^2-(1+e^{2\rho})Y_1^2\right](y)
  \eqno{(3.18)}$$
and both $G(x)(\geq 0)$ and $H(y)$ are constant functions. 

In particular, for a given $(\sigma(x),\rho(y))$, all pairs $(X_1(x),Y_1(y))$ 
satisfying \ $(\bar\kappa_3^2-\bar\zeta)(x,y)\equiv 0$ \ 
generate a 5-dimensional set. Furthermore, for any pair $(X_1(x),Y_1(y))\neq (0,0)$ satisfying \ $(\bar\kappa_3^2-\bar\zeta)(x,y)\equiv 0$, \ $\bar\kappa_1$ and $\bar\kappa_2$ do not vanish identically. \\

\vspace*{2mm}
{\sc Proof.} \ We firstly summarize equations following from $\varphi(x,y)$ and $\varphi_z(x,y)$:
$$\varphi_y=\rho'\sin\varphi\cos\varphi, \ \ \varphi_{yy}=(\rho^{''}+\rho'^2\cos2\varphi)\sin\varphi\cos\varphi, \ \ 
\varphi_{zx}=\sigma'\sin\varphi, \ \ \varphi_{zy}=\sigma\rho'\sin\varphi\cos^2\varphi, $$ 
$$(\cos\varphi)_y=-\rho'\sin^2\varphi\cos\varphi, \ \ 
(\cos\varphi)_{yy}=-[\rho^{''}+\rho'^2(\cos2\varphi+\cos^2\varphi)]\sin^2\varphi\cos\varphi,$$
$$L\psi=(1/2)[\sigma^2-2\rho^{''}\sin^2\varphi-\rho'^2\sin^2\varphi(2\cos^2\varphi+\cos2\varphi)],$$
$$\varphi_{zz}=(\sigma^2-\rho^{''}-\rho'^2\cos^2\varphi)\sin\varphi\cos\varphi,$$
$$2\psi_{zz}=-\sigma^2\cos2\varphi-2\rho~{''}\sin^2\varphi-\rho'^2\sin^2\varphi\cos2\varphi,$$
$$1+\varphi_z^2-2\psi_{zz}=1+\sigma^2\cos^2\varphi+2\rho^{''}\sin^2\varphi+\rho'^2\sin^2\varphi\cos2\varphi.$$

Now, we study \ $(\bar \kappa_3)_x=-(e^{-\bar P})_x\tan\varphi$ \ and \ $(\bar \kappa_3)_y=(e^{-\bar P})_y\cot\varphi$ \ to obtain $e^{-\bar P}(x,y)$ and $\bar\kappa_3(x,y)$. By the first equation, we have \ 
$\bar\kappa_3(x,y)=-(e^{\rho}e^{-\bar P})(x,y)+ Y_2(y)$ \ with a one-variable function $ Y_2(y)$. By the second equation we have \ $Y_2'=(e^{-\bar P})_y(e^{\rho}+e^{-\rho})+\rho'e^{\rho}e^{-\bar P}$. \ 
Differentiating this equation by $x$, we have \ $(e^{-\bar P})_{xy}/(e^{-\bar P})_x+\rho'e^{2\rho}/(1+e^{2\rho})=0.$ \
Then, we have \ $e^{-\bar P}(x,y)=\cos\varphi(y)X_1(x) +Y_1(y)$, \ $\bar\kappa_3(x,y)=-\sin\varphi(y) X_1(x)-(e^{\rho}Y_1)(y)+Y_2(y)$ \ and \ $Y_2'=(e^{\rho}Y_1)'+e^{-\rho}Y_1'.$

Next, we study (3.9) and (3.10) to obtain $(e^{-\bar P})_z(x,y)$. Let us define a one-variable function $X_2(x)$ by $X_2(x):=\int^x \sigma'X_1dx$. 
By (3.9), we have \ $(e^{-\bar P})_{zx}=[-\sigma \sin^2\varphi~X_1 +X_2+\sigma\cos\varphi~ Y_1 ]_x$. In the same way, we have \ 
$(e^{-\bar P})_{zy}=[-\sigma \sin^2\varphi~ X_1+\sigma\cos\varphi~ Y_1]_y$ \ by (3.10). 
Hence, we have \ $(e^{-\bar P})_{z}=-\sigma \sin^2\varphi~ X_1+X_2+\sigma\cos\varphi~ Y_1$. By the argument above, $e^{-\bar P}(x,y)$, $(e^{-\bar P})_z(x,y)$ and $\bar\kappa_3(x,y)$ have been determined by (3.16) from $X_1(x)$ and $Y_1(y)$, and further $X_2(x)$ and $Y_2(y)$ also have been determined in the desired forms.

Now, in order to determine $X_1(x)$ and $Y_1(y)$, we study the equation of $\bar \kappa_3$ in Definition and Notation 3-(2). Then, we have (3.13) and (3.14). 
In fact, we have 
$$-\sin\varphi~X_1-e^{\rho}Y_1+Y_2=(X_1^{''}+\sigma^2X_1-\sigma X_2)\sin\varphi$$
$$-[e^{-\rho}~Y_1^{''}-\rho'e^{-\rho}\cos2\varphi~ Y_1'+(\rho^{''}+\rho'^2\cos^2\varphi)\sin\varphi\cos\varphi~ Y_1]$$
by substituting $e^{-\bar P}$, $(e^{-\bar P})_z$ and $\bar\kappa_3$ in (3.16) for those in the equation of $\bar \kappa_3$. 
In this equation, $\sigma$ and $X_i$ are one-variable functions for $x$ and $\rho$, $\varphi$ and $Y_i$ are one-variable functions for $y$, hence we firstly have the following two equations:
$$X_1''+(1+\sigma^2)X_1-\sigma X_2=c,$$
$$e^{\rho}Y_1-Y_2-[e^{-\rho}Y_1''-\rho'e^{-\rho}\cos2\varphi Y_1'+(\rho''+\rho'^2\cos^2\varphi)\sin\varphi\cos\varphi Y_1]=-c\sin\varphi$$
with a constant $c$. Next, in these two equations we show that it is sufficient to study the case $c=0$, that is, the case of (3.13) and (3.14). Suppose $c\neq 0$. Then, with $X_2(x)=c\sigma(x)$, \ $X_1(x)\equiv c$ \ is a solution to the first equation. With $Y_2(y)\equiv 0$, \ $Y_1(y)=-c\cos\varphi$ \ is a solution to the second equation and (3.15). These two equations are linear, and $e^{-\bar P}$, $\bar\kappa_3$ and $(e^{-\bar P})_z$ in (3.16) do not depend on all $(X_1,Y_1)=(c,-c\cos\varphi)$ with $c\neq 0$. Hence, the assertion has been shown.

It still remains to be shown that $X_1(x)$ is a solution to (3.12).  
The equation (3.13) is rewritten as \ $[X_1^{''}(x)+(1+\sigma^2)(x)X_1(x)-\sigma(x)\int^x_0\sigma'X_1dx]/\sigma(x)=c_0$.\  We differentiate 
this equation by $x$, then we have  
$$\left[\frac{X_1^{''}}{\sigma}+\frac{(1+\sigma^2)}{\sigma}X_1-\int^x_0\sigma'X_1dx \right]_x=
\frac{1}{\sigma^2}[\sigma X_1^{'''}-\sigma'X_1^{''}+\sigma(1+\sigma^2)X_1'-\sigma'X_1]=0.$$
Hence, we have showed that (3.12) and (3.13) are equivalent to each other.
In the same way, we can show that $Y_1(y)$ satisfies a linear differential equation of third order, by (3.14) and (3.15).

Next, we study the condition $\bar\kappa_3^2-\bar\zeta\equiv 0$. We firstly 
have 
$$(e^{-\bar P})_x=\cos\varphi~X_1^{'}, \ \ \ \  (e^{-\bar P})_{xx}=\cos\varphi[-(1+\sigma^2)X_1+\sigma X_2], \ \ \ \  (e^{-\bar P})_y=-\rho'\sin^2\varphi\cos\varphi~X_1+Y_1^{'},$$
$$(e^{-\bar P})_{yy}=-\{\rho{''}+\rho'^2(\cos2\varphi+\cos^2\varphi)\}\sin^2\varphi\cos\varphi~X_1-e^{\rho}Y_2+\rho{'}\cos2\varphi~Y_1{'}$$
$$-\{(\rho{''}+\rho'^2\cos^2\varphi)\sin^2\varphi-e^{2\rho}\}Y_1 $$ 
by (3.13) and (3.14).
Then, by replacing $(e^{-\bar P})_{zz}$ in $\bar\zeta$ with it in Definition and Notation 3-(2), we have \ $\bar\kappa_3^2-\bar\zeta=G(x)+H(y).$ \ 
Next, the derivative of $G(x)$ (resp. $H(y)$) vanishes by (3.13) and the definition of $X_2(x)$ (resp. by (3.14) and (3.15)), hence $G(x)$ and $H(y)$) are constant functions.

Now, we show that all pairs $(X_1(x),Y_1(y))$ such that \ $(\bar\kappa_3^2-\bar\zeta)(x,y)\equiv 0$ \ generate a 5-dimensional set. 
Let $c_2$ be a constant determined by
$$c_2:=e^{-\rho(0)}Y_1'(0)+\frac{(\rho'(0)e^{\rho(0)})}{(1+e^{2\rho(0)})}Y_1(0).\eqno{(3.19)}$$
Then, we have 
$$-H(0)+(c_1^2+c_2^2)(1+e^{2\rho(0)})=(Y_1(0)-c_1e^{\rho(0)})^2 \eqno{(3.20)}$$
by (3.18) and $Y_2(0)=e^{\rho(0)}Y_1(0)+c_1$. Hence, for any $(c_1,c_2)$ and non-negative $-H(0)$, there is $Y_1(0)$ such that (3.20) holds. 

By the argument above, we firstly choose a solution $X_1(x)$ arbitrarily. 
For $X_1(x)$, $G(0)=-H(0)$ is determined. Then, for any pair $(c_1,c_2)$, $Y_1(0)$ is determined by (3.20). Next, for their $c_2$ and $Y_1(0)$, $Y_1'(0)$ is determined by (3.19). Note that, for $(Y_1(0),Y_1'(0),c_1)$, $Y_1(x)$ is uniquely determined, since $Y_1''(0)$ is also determined by (3.14).
Since the space of solutions $X_1(x)$ is 3-dimensional and pairs $(c_1,c_2)$ are 2-dimensional, we have showed that all pairs $(X_1(x),Y_1(y))$ satisfying \ $(\bar\kappa_3^2-\bar\zeta)(x,y)\equiv 0$ \ generate a 5-dimensional set.

Finally, we show that, for any pair $(X_1(x),Y_1(y))\neq (0,0)$ satisfying \ $(\bar\kappa_3^2-\bar\zeta)(x,y)\equiv 0$, \ each function $\bar\kappa_i \ (i=1,2)$ does not vanish identically. We firstly have $\bar \kappa_1=Y_2$ by (3.16). Suppose $Y_2\equiv 0$. Then, by $Y_2'\equiv 0$ we have $Y_1=c_3\cos\varphi$ with a constant $c_3$. Furthermore, we have $c_3\sin\varphi\equiv 0$ by (3.14), that is, $c_3=0$. Hence, we have $H(y)\equiv 0$ by $Y_1=Y_2\equiv 0$ and (3.18), then we have $X_1\equiv 0$ by $G(x)\equiv 0$ and (3.17). In consequence, if $\bar\kappa_1\equiv 0$, then we have $e^{-\bar P}\equiv 0$, which is not the case that we want. 
For $\bar\kappa_2$, we have 
$\bar\kappa_2(x,y)=-(1/\sin\varphi)(y)X_1(x)-(e^{\rho}+e^{-\rho})(y)Y_1(y)+Y_2(y)$. If $\bar\kappa_2\equiv 0$, then we also have \ $X_1=Y_1=Y_2=e^{-\bar P}\equiv 0$ \ in the same way.  
\hspace{\fill}$\Box$\\

\vspace*{2mm}
The following lemma, which is directly verified, is useful for Example 2: \\

{\sc Lemma 3.1.} \ {\it Let $\hat g\in Met^0$, and $(\varphi(x,y),\varphi_z(x,y))$ be functions arising from $\hat g$. Suppose that $\bar \kappa_3(x,y)$ and $e^{-\bar P}(x,y)$ satisfy \ $(\bar \kappa_3)_x=-(e^{-\bar P})_x\tan\varphi$ \ and \ $(\bar \kappa_3)_y=(e^{-\bar P})_y\cot\varphi$. \ Then, (3.9) and (3.10) for $e^{-\bar P}(x,y)$ and $(e^{-\bar P})_z(x,y)$, respectively, are equivalent to the following equations: 
$$[(e^{-\bar P})_z-\bar\kappa_3\varphi_z]_x=\varphi_{zx}(e^{-\bar P}\cot\varphi-\bar\kappa_3) \ \ (=-\bar\kappa_2\varphi_{zx}),$$
$$[(e^{-\bar P})_z-\bar\kappa_3\varphi_z]_y=-\varphi_{zy}(e^{-\bar P}\tan\varphi+\bar\kappa_3) \ \ (=-\bar\kappa_1\varphi_{zy}).$$
}

\vspace*{2mm}
{\sc Example 2.} \ The hyperbolic 2-metric $\hat g=[(dx)^2+(dy)^2]/y^2$ on the upper half plane belongs to $Met^0$ (\cite{bhs}, \S2.2, Example 1).  
$\hat g$ gives rise to functions 
$\varphi(x,y)$, $\varphi_z(x,y)$ and $\psi(x,y)$ 
such that \ $\cos\varphi(x,y)=(x^2-y^2)/(x^2+y^2),$ \ $\sin\varphi(x,y)=2xy/(x^2+y^2)$ \ and 
$$\varphi_z(x,y)=y/(x^2+y^2), \ \ \ \ \psi(x,y)=\log(x^2+y^2)-(9/8)\log |x|.$$  

\vspace*{2mm}
For arbitrarily fixed two constants $c_0$ and $c_1$, let $X_1=X_1(x)$ and $Y=Y(y)$ be solutions to linear ordinary differential equations   
$$xX_1 ^{'''}-X_1 ^{''}+(x+\frac{9}{4x})X_1^{'}-X_1=c_0x^2+c_1, \ \ \ \ \ \ Y^{''}+Y=c_0y^2+c_1, \eqno{(3.21)}$$
respectively. 
Note that all solutions $X_1(x)$ and $Y(y)$ to each equation are defined on whole real number $R$: the point $x=0$ for the equation of $X_1(x)$ is a regular singularity and \ $X_1(x):=c_0(x^2+5/2)-c_1$ \ and \ $Y(y):=c_0(y^2-2)+c_1$ \ are particular solutions to each equations.

We define a one-variable function $X_2(x)$ by \ 
$(X_2/x)^{'}=(9/(4x^3))X_1^{'}$, \ then the first equation in (3.21) is equivalent to the following equation (3.22):
$$X_1^{''}+X_1+X_2=c_0x^2-c_1. \eqno{(3.22)}$$
Let \ $X_2(x)=(9x/4)\int_1^x(1/x^3) X_1^{'}(x)dx+c_2x$ \ with a constant $c_2$.

From \ $(\bar\kappa_3)_x=-(e^{-\bar P})_x\tan\varphi$, \ $(\bar\kappa_3)_y=(e^{-\bar P})_y\cot\varphi$, \ (3.9), (3.10) and the definition of $\bar\kappa_3$, the functions $e^{-\bar P}(x,y),$ $\bar\kappa_3(x,y)$ and $(e^{-\bar P})_z(x,y))$ are determined as follows:  
$$\begin{array}{l}
e^{-\bar P}=X_1^{'}-\displaystyle\frac{2x}{x^2+y^2}A, \ \ \ \ \ \bar\kappa_3=-Y^{'}+\displaystyle\frac{2y}{x^2+y^2}A, \\[3mm]
(e^{-\bar P})_z=\displaystyle\frac{1}{x^2+y^2}\left(xX_1^{'}-A+\displaystyle\frac{2y^2}{x^2+y^2}A\right)
\end{array}  \eqno{(3.23)}$$
with any solution $(X_1(x),Y(y))$ to two equations in (3.21), where \ $A:=\bar X +\bar Y$, \ $\bar X:=xX_1^{'}-X_1$ \ and \ $\bar Y:=yY^{'}-Y$. \ In particular, three functions in (3.23) do not depend on $c_1$. 

For functions in (3.23), we have \ $(\bar\kappa_3^2-\bar\zeta)(x,y)=G(x)+H(y)$, \ where 
$$G(x):=(X_1+c_1+X_2-c_0x^2)^2+4c_0(X_1+c_1)+(1+\frac{9}{4x^2})(X_1^{'})^2+(\frac{2}{x}X_2-4c_0x)X_1^{'},$$
$$H(y):=(Y^{'}-2c_0y)^2+(Y-c_1-c_0y^2+2c_0)^2-4c_0^2$$
and both $G(x)$ and $H(y)$ are constant functions. 
Furthermore, all pairs $(X_1(x),Y(y))$ of solutions satisfying $\bar\kappa_3^2-\bar\zeta\equiv 0$ generate a 5-dimensional set. \\

\vspace*{2mm}
{\sc Proof.} \ We firstly summarize equations following from $\varphi(x,y)$ and $\varphi_z(x,y)$:
$$\varphi_x=-\frac{2y}{x^2+y^2}, \ \ \varphi_y=\frac{2x}{x^2+y^2}, \ \ \psi_z=\frac{-x}{x^2+y^2}, \ \ \varphi_{zx}=\frac{-2xy}{(x^2+y^2)^2}, \ \ \varphi_{zy}=\frac{x^2-y^2}{(x^2+y^2)^2}, \ \ \varphi_{xx}=\frac{4xy}{(x^2+y^2)^2},$$
$$\varphi_{xy}=-2\frac{x^2-y^2}{(x^2+y^2)^2}, \ \ \varphi_{yy}=\frac{-4xy}{(x^2+y^2)^2}, \ \ L\varphi=\frac{8xy}{(x^2+y^2)^2}, \ \ \psi_x=\frac{2x}{x^2+y^2}-\frac{9}{8x}, \ \ \psi_y=\frac{2y}{x^2+y^2},$$
$$\psi_{xx}=-2\frac{x^2-y^2}{(x^2+y^2)^2}+\frac{9}{8x^2}, \ \ 
\psi_{yy}=2\frac{x^2-y^2}{(x^2+y^2)^2}, \ \ 
L\psi=-4\frac{x^2-y^2}{(x^2+y^2)^2}+\frac{9}{8x^2}, $$
$$\varphi_{zz}=-\frac{y}{2x(x^2+y^2)^2}(7x^2+9y^2), \ \ \psi_{zz}=\frac{1}{8x^2(x^2+y^2)^2}(23x^4+22x^2y^2-9y^4),$$
$$1+\varphi_z^2-2\psi_{zz}=1-\frac{23x^4+18x^2y^2-9y^4}{4x^2(x^2+y^2)^2}.$$

Now, suppose that $(e^{-\bar P}(x,y),(e^{-\bar P})_z(x,y))$ satisfy \ $(\bar \kappa_3)_x=-(e^{-\bar P})_x\tan\varphi$, \ $(\bar \kappa_3)_y=(e^{-\bar P})_y\cot\varphi$, \ (3.9) and (3.10). Let \ $f(x,y):=[(e^{-\bar P})_z-\bar\kappa_3\varphi_z](x,y)$. \ Then, by Lemma 3.1 we have \ $-f_x=e^{-\bar P}\varphi_{zy}+\bar\kappa_3\varphi_{zx}$ \ and \ $-f_y=-e^{-\bar P}\varphi_{zx}+\bar\kappa_3\varphi_{zy}$. \ Hence, we have 
$$e^{-\bar P}=-(x^2-y^2)f_x-2xyf_y, \ \ \ \ \bar\kappa_3=2xyf_x-(x^2-y^2)f_y. 
\eqno{(3.24)}$$

Next, by (3.24), each equation of \ $(\bar\kappa_3)_x=-(e^{-\bar P})_x\tan\varphi$ \ and \ $(\bar\kappa_3)_y=(e^{-\bar P})_y\cot\varphi$ \ is equivalent to \ $2yf_x+2xf_y+(x^2+y^2)f_{xy}=0.$ \ Hence, there are two functions $X_1=X_1(x)$ and $Y=Y(y)$ of one-variable such that \ $f(x,y)=(X_1+Y)/(x^2+y^2)$ \ holds. 
Then, we obtain $e^{-\bar P}$, $(e^{-\bar P})_z$ and $\bar\kappa_3$ in (3.23) from (3.24). 

Next, we substitute $e^{-\bar P}$, $(e^{-\bar P})_z$ and $\bar\kappa_3$ in (3.23) for those in the definition of $\bar\kappa_3$. Then, we have 
$$-2xX_1^{'''}+2X_1^{''}-(2x+\frac{9}{2x})X_1^{'}+2X_1+\frac{(x^2-y^2)}{y}(Y^{'''}+Y^{'})+2(Y^{''}+Y)=0.$$
In this equation, the term $[(x^2-y^2)/y](Y^{'''}+Y^{'})$ is expressed as sum of two one-variable functions for $x$ and $y$. Hence, we have \ $Y^{'''}+Y^{'}=2c_0y$ \ with a constant $c_0$. In consequence, we obtain the equations in (3.21) and (3.22). 
However, we can fix $c_1=0$ in the equations of (3.21) and (3.22). In fact, let $X_1(x)$ and $Y(y)$, respectively, be solutions of two equations in (3.21), as above. Then, \ $\hat X_1(x):=X_1(x)+c_1$ \ and \ $\hat Y(y):=Y(y)-c_1$, \ respectively, are solutions to these equations with $c_1=0$. 
Furthermore, $e^{-\bar P}$, $\bar\kappa_3$ and $(e^{-\bar P})_z$ in (3.23) are determined by $(\hat X_1(x),\hat Y(y))$ (not by $(X_1(x),Y(y))$).

The equation of $X_1(x)$ in (3.21) has a regular singularity at $x=0$. In fact, let $c_0=c_1=0$, and suppose that a solution $X_1(x)$ is expressed as $X_1(x)=x^r(1+\Sigma_{n=1}^{\infty}a_nx^n)$. Then, we have \ $r=0$ and \ $r=2\pm\sqrt{-5}/2$. \ For $r=2\pm\sqrt{-5}/2$, we have \ $x^{2\pm\sqrt{-5}/2}=x^2[\cos(\sqrt{5}/2\log |x|)\pm\sqrt{-1}\sin(\sqrt{5}/2\log |x|)]$. \ 
Furthermore, we have \ $a_{2n-1}=0$ \ for all natural numbers $n$ and  
$$(r-1)+(r+2)(r^2+5/4)a_2=0$$ 
$$(n+r+2)[(n+r)^2+5/4]a_{n+2}+(n+r-1)a_n=0 \hspace*{1cm}for \ \ n\geq 2.$$ 

Now, by (3.21), (3.22) and (3.23), we have 
$$
\begin{array}{l}
(e^{-\bar P})_x=\displaystyle\cos\varphi[(X_1+c_1+X_2+\frac{2}{x^2+y^2}A)-c_0x^2)],\\[3mm]
(e^{-\bar P})_y=\displaystyle\sin\varphi[(Y-c_1+\frac{2}{x^2+y^2}A)-c_0y^2],\\[3mm]
(e^{-\bar P})_{xx}=\displaystyle\frac{x^2-y^2}{x^2+y^2}(1+\frac{9}{4x^2})X_1^{'}-\frac{2x(x^2-3y^2)}{(x^2+y^2)^2}(X_1+c_1)-\frac{x^4-6x^2y^2+y^4}{x(x^2+y^2)^2}X_2 \\[3mm]
\hspace*{2cm} -\displaystyle\frac{4x(x^2-3y^2)}{(x^2+y^2)^3}A-2c_0\frac{xy^2(3x^2-y^2)}{(x^2+y^2)^2}
.\\[3mm]
(e^{-\bar P})_{yy}=\displaystyle\frac{2xy}{x^2+y^2}Y^{'}+\frac{2x(x^2-3y^2)}{(x^2+y^2)^2}(Y-c_1)+\frac{4x(x^2-3y^2)}{(x^2+y^2)^3}A -\displaystyle{2c_0\frac{xy^2(3x^2-y^2)}{(x^2+y^2)^2}}. 
\end{array} \eqno{(3.25)}$$
Then, we have $(\bar\kappa_3)^2-\bar\zeta=G(x)+H(y)$ by (3.23), (3.25) and the definition of $\bar\zeta$. In this calculation, we attain to the desired equation by firstly collecting terms of the same degree for $(x^2+y^2)^k, \ k=0,-1,-2,-3,-4$. Furthermore, since $G^{'}(x)\equiv 0$ and $H^{'}(y)\equiv 0$ by (3.21) and (3.22), both $G(x)$ and $H(y)$ are constant functions.

Finally, we show that all solutions $(X_1(x),Y(y))$ satisfying $G(x)+H(y)\equiv 0$ generate a 5-dimensional set. Let us fix $c_1=0$.
By $X_2(1)=c_2$, we have 
$$H(0)=(Y^{'}(0))^2+[(Y(0))^2+4c_0Y(0)],$$
$$-H(0)=G(1)=[(X_1(1))^2+2(c_0+c_2)X_1(1)+(-c_0+c_2)^2]+2X_1^{'}(1)[\frac{13}{8}X_1^{'}(1)+(c_2-2c_0)].$$
Let \ $2c_3:=Y^{'}(0)$ \ and \ $2c_4:=X_1^{'}(1)[\frac{13}{8}X_1^{'}(1)+(c_2-2c_0)]$. \ Then, by $X_1(1),\ X_1'(1)\in R$ we have \ 
$$13c_4\geq-(c_2-2c_0)^2, \hspace*{1cm} 4(c_0c_2-c_4)\geq H(0)\geq 4(-c_0^2+c_3^2)$$  
and \ $Y(0)=-2c_0\pm \sqrt{4(c_0^2-c_3^2)+H(0)},$ \  
$X_1(1)=-(c_0+c_2)\pm\sqrt{4(c_0c_2-c_4)-H(0)}.$ \ Hence, we have \ 
$$17c_0^2+9c_0c_2+c_2^2\geq 13c_3^2    \eqno{(3.26)}$$
by \ $13c_0(c_0+c_2)-13c_3^2\geq 13c_4\geq -(c_2-2c_0)^2$. \ Now, for any $(c_0,c_2,c_3)$ satisfying (3.26), we can choose $c_4$ and $H(0)$ in order such that they satisfy the inequalities above. In consequence, all pairs $(X_1(x),Y(y))$ such that $G(x)+H(y)\equiv 0$ generate a 5-dimensional set.
\hspace{\fill}$\Box$\\

\newpage
\vspace{2mm}
\noindent
{\bf \large  4. Condition on $(\bar P(x,y),\bar P_z(x,y))$ and surfaces in $S^3$.}

Let $\hat{g}\in Met^0$, and $(\varphi(x,y),\varphi_z(x,y))$ be (analytic) functions on $V$ arising from $\hat{g}$. 
In this section, we firstly give a geometrical interpretation of the assumption for a pair $({\bar P}(x,y),{\bar P}_z(x,y))$ in Theorem 1 of \S3: for such a pair $({\bar P}(x,y),{\bar P}_z(x,y))$, a surface $\phi(x,y)$ in $S^3$ is determined, and then (3.11) is the Gauss equation for $\phi(x,y)$. 
Next, we show that the surface $\phi(x,y)$ gives rise to an evolution $\phi^z(x,y)$ of surfaces in $S^3$ issuing from $\phi(x,y)$ and 
the evolution $\phi^z$ also determines an evolution $f^z(x,y)$ of surfaces in $R^4$ such that $f(x,y,z):=f^z(x,y)$ is a generic conformally flat hypersurface.

\vspace*{2mm}
Now, for any functions $\bar P(x,y)$ and $\bar P_z(x,y)$ on $V$, we have \ 
$${\bar \kappa}_1e^{\bar P}\cos\varphi\ =\sin\varphi+{\bar \kappa}_3e^{\bar P}\cos\varphi, \ \ \ \ 
{\bar \kappa}_2e^{\bar P}\sin\varphi\ =-\cos\varphi+{\bar \kappa}_3e^{\bar P}\sin\varphi   \eqno{(4.1)}$$
by the definitions of ${\bar \kappa}_i \ (i=1,2,3)$ in Definition and Notation 3 in \S3. 
We have denoted by $\nabla'$ the standard connection on $R^4$.\\ 

\vspace*{2mm}
{\pro 4.1}. \ \ {\it Let $\hat{g}\in Met^0$, and $(\varphi(x,y),\varphi_z(x,y))$ be functions on $V$ arising from $\hat{g}$. 
Let $\phi:V\ni(x,y)\mapsto \phi(x,y)\in S^3$ be a generic surface and $(x,y)$ be a principal curvature line coordinate system.  
Let $\bar P(x,y)$ and $\bar P_z(x,y)$ be functions on $V$. 
Suppose that $\phi(x,y)$ satisfies
$$\phi_x=-{\bar \kappa}_1e^{\bar P}\cos\varphi\ X^0_{\alpha}, \ \ \ \ 
\phi_y=-{\bar \kappa}_2e^{\bar P}\sin\varphi\ X^0_{\beta}  \eqno{(4.2)}$$
with respect to orthonormal vector fields $X^0_{\alpha}(x,y)$ and $X^0_{\beta}(x,y)$. 
Then, the following facts (1), (2) and (3) are equivalent to each other: 

(1) there is a surface $f^0:V\ni(x,y)\mapsto f^0(x,y)\in R^4$ such that  
$$f^0_x=e^{\bar P}\cos\varphi\ X^0_{\alpha}, \ \ f^0_y=e^{\bar P}\sin\varphi\ X^0_{\beta}.$$  

(2) $\nabla'_{\partial/\partial y}X^0_{\alpha}$ and $\nabla'_{\partial/\partial x}X^0_{\beta}$ satisfy 
$$\nabla'_{\partial/\partial y}X^0_{\alpha}=\frac{(e^{\bar P}\sin\varphi)_x}{e^{\bar P}\cos\varphi}~X_{\beta}^0, \ \ \ \ 
\nabla'_{\partial/\partial x}X^0_{\beta}=\frac{(e^{\bar P}\cos\varphi)_y}{e^{\bar P}\sin\varphi}~X_{\alpha}^0.  \eqno{(4.3)}$$ 

(3)\ $({\bar \kappa}_3)_x=-(e^{-\bar P})_x\tan\varphi$\ and ~$({\bar \kappa}_3)_y=(e^{-\bar P})_y\cot\varphi$\ 
are satisfied. In particular, $\bar P(x,y)$ satisfies (3.8) defined by $(\varphi(x,y),\varphi_z(x,y))$. }\\

\vspace*{2mm}

In Proposition 4.1, the assumption that $\phi$ is a surface implies that $\bar\kappa_1\bar\kappa_2\neq 0$ holds on $V$. 
The first fundamental form $I_{\phi}$ of $\phi$ is expressed as 
$$I_{\phi}:=\langle d\phi,d\phi\rangle=({\bar \kappa}_1e^{\bar P}\cos\varphi)^2dx^2+
({\bar \kappa}_2e^{\bar P}\sin\varphi)^2dy^2  \eqno{(4.4)}$$
by (4.2), where $\langle x,y\rangle$ is the standard inner product in $R^4$.
The first fundamental form $I_{f^0}$ of $f^0$ is expressed as 
$$I_{f^0}:=\langle df^0,df^0\rangle=e^{2\bar P}(\cos^2\varphi dx^2+\sin^2\varphi dy^2)  \eqno{(4.5)}$$  
by (1). The second statement in (3) follows directly from the integrability condition on $\bar\kappa_3$: \ 
$(\bar\kappa_3)_{xy}=(\bar\kappa_3)_{yx}$.

We break up the proof of Proposition 4.1 into three lemmata below. 
Let $D$ (resp. $D'$) be the Riemannian connection on $\phi$ (resp. on $S^3$) determined from the metric $I_{\phi}$ 
(resp. the canonical metric of $S^3$). 
Let $\xi$ be a unit normal vector field (in $S^3$) of $\phi$.

\vspace*{2mm}
{\lem 4.1}. \ \ {\it Let $\phi(x,y)$ be a generic surface in $S^3$ satisfying (4.2). Suppose that $(x,y)$ is a principal curvature line coordinate system. Then, there are functions $b_i(x,y)$ and $c_i(x,y)$ $(i=1,2)$ such that  
$$\nabla'_{\partial/\partial x}\xi=D'_{\partial/\partial x}\xi=b_1X^0_{\alpha}, \hspace*{1cm} 
\nabla'_{\partial/\partial y}\xi=D'_{\partial/\partial y}\xi=b_2X^0_{\beta},  \eqno{(4.6)}$$
$$\nabla'_{\partial/\partial x}X^0_{\beta}=D'_{\partial/\partial x}X^0_{\beta}=D_{\partial/\partial x}X^0_{\beta}=c_1 X^0_{\alpha} \ \ \ \ \ 
\nabla'_{\partial/\partial y}X^0_{\alpha}=D'_{\partial/\partial y}X^0_{\alpha}=D_{\partial/\partial y}X^0_{\alpha}=c_2 X^0_{\beta}.$$  }

\vspace*{2mm}
{\proof} \ Since $(x,y)$ is a principal curvature line coordinate system, we have \ $D'_{\partial/\partial x}\xi=b_1X^0_{\alpha}$ \ and \ $D'_{\partial/\partial y}\xi=b_2X^0_{\beta}$. \  
For (4.6), we have only to show that both $\nabla'_{\partial/\partial y}\xi$ and $\nabla'_{\partial/\partial x}\xi$ are perpendicular to $\phi$.
We have \  
$\langle \nabla'_{\partial/\partial y}\xi,\phi\rangle+\langle \xi,
\nabla'_{\partial/\partial y}\phi\rangle=0$ \ and \ 
$\langle \nabla'_{\partial/\partial x}\xi,\phi\rangle+\langle \xi,
\nabla'_{\partial/\partial x}\phi\rangle=0$ \ by $\langle \xi,\phi\rangle=0$. \ Then, since 
$\xi$ is a normal vector field of $\phi$, we have \ $\langle \xi,\phi_x\rangle=
\langle \xi,\phi_y\rangle=0$ and $\phi_x=\nabla'_{\partial/\partial x}\phi$, \ 
$\phi_y=\nabla'_{\partial/\partial y}\phi$, \ which shows (4.6).  

Next, we express the derivatives $\nabla'_{\partial/\partial x}X^0_{\beta}$  
and $\nabla'_{\partial/\partial y}X^0_{\alpha}$ as 
$$\nabla'_{\partial/\partial x}X^0_{\beta}=c_1X^0_{\alpha}+d_1\xi+e_1\phi, \ \ \ \ \nabla'_{\partial/\partial y}X^0_{\alpha}=c_2X^0_{\beta}+d_2\xi+e_2\phi.$$
We have to show $d_1=d_2=e_1=e_2\equiv 0.$  
By (4.2) and (4.6), we firstly have $d_1=d_2\equiv 0$. 
Next, since $\phi_x=\nabla'_{\partial/\partial x}\phi$, 
$\phi_y=\nabla'_{\partial/\partial y}\phi$ and (4.2), we have $e_1=e_2\equiv 0$. 
\hspace{\fill}$\Box$\\

\vspace*{2mm}
{\lem 4.2}. \ \ {\it Let $\phi(x,y)$ be a generic surface in $S^3$ satisfying (4.2). Suppose that $(x,y)$ is a principal curvature line coordinate system.     
Then, there is a surface $f^0:V\ni(x,y)\mapsto f^0(x,y)\in R^4$ satisfying \ 
$f^0_x=e^{\bar P}\cos\varphi X^0_{\alpha}$ \ and \ $f^0_y=e^{\bar P}\sin\varphi X^0_{\beta}$, \ if and only if (4.3) is satisfied. }

\vspace*{2mm}
{\proof} \ We have only to study the integrability condition on $f^0$. We have 
$$f^0_{xy}=(e^{\bar P}\cos\varphi)_y X^0_{\alpha}+e^{\bar P}\cos\varphi (X^0_{\alpha})_y, \ \ \ \ 
f^0_{yx}=(e^{\bar P}\sin\varphi)_x X^0_{\beta}+e^{\bar P}\sin\varphi (X^0_{\beta})_x.$$
Then, by $(X^0_{\alpha})_y=\nabla'_{\partial/\partial y}X^0_{\alpha}$, \ 
$(X^0_{\beta})_x=\nabla'_{\partial/\partial x}X^0_{\beta}$ \ and Lemma 4.1, \  
$f^0_{xy}=f^0_{yx}$ holds if and only if (4.3) is satisfied. 
\hspace{\fill}$\Box$ \\

\vspace*{2mm}
{\lem 4.3}. \ \ {\it Let $\phi(x,y)$ be a generic surface in $S^3$ satisfying (4.2). Suppose that $(x,y)$ is a principal curvature line coordinate system.   
Then, there is a surface $f^0:V\ni(x,y)\mapsto f^0(x,y)\in R^4$ satisfying \ 
$f^0_x=e^{\bar P}\cos\varphi X^0_{\alpha}$ \ and \ $f^0_y=e^{\bar P}\sin\varphi X^0_{\beta}$, \ 
if and only if \ $({\bar \kappa}_3)_x=-(e^{-\bar P})_x\tan\varphi$ \ and \ $({\bar \kappa}_3)_y=(e^{-\bar P})_y\cot\varphi$ \ 
are satisfied. }\\

\vspace*{2mm}
{\proof} \ Since the integrability condition on $\phi$ is given by 
$$-\phi_{xy}=(\bar\kappa_1e^{\bar P}\cos\varphi)_yX_{\alpha}^0+\bar\kappa_1e^{\bar P}\cos\varphi(X_{\alpha}^0)_y
=(\bar\kappa_2e^{\bar P}\sin\varphi)_xX_{\beta}^0+\bar\kappa_2e^{\bar P}\sin\varphi(X_{\beta}^0)_x=-\phi_{yx},$$
we have 
$$\nabla'_{\partial/\partial y}X_{\alpha}^0=\frac{(\bar\kappa_2e^{\bar P}\sin\varphi)_x}{\bar\kappa_1e^{\bar P}\cos\varphi}X_{\beta}^0, \ \ \ \ 
\nabla'_{\partial/\partial x}X_{\beta}^0=\frac{(\bar\kappa_1e^{\bar P}\cos\varphi)_y}{\bar\kappa_2e^{\bar P}\sin\varphi}X_{\alpha}^0  $$
by Lemma 4.1. By these equations and Lemma 4.2, 
there exists a surface 
$f^0:V\ni(x,y)\mapsto f^0(x,y)\in R^4$ satisfying \ $f^0_x=e^{\bar P}\cos\varphi X^0_{\alpha}, \ f^0_y=e^{\bar P}\sin\varphi X^0_{\beta}$, 
if and only if \ 
$(\bar\kappa_2e^{\bar P}\sin\varphi)_x=\bar\kappa_1(e^{\bar P}\sin\varphi)_x$ \ and \ 
$(\bar\kappa_1e^{\bar P}\cos\varphi)_y=\bar\kappa_2(e^{\bar P}\cos\varphi)_y$ \ hold. 
From the definitions of $\bar\kappa_1$, $\bar\kappa_2$ and (4.1), these equations, respectively, are equivalent to \ 
$({\bar \kappa}_3)_x=-(e^{-\bar P})_x\tan\varphi$ \ and \ $({\bar \kappa}_3)_y=(e^{-\bar P})_y\cot\varphi$.  
\hspace{\fill}$\Box$ \\

By Lemmata 4.1, 4.2 and 4.3, Proposition 4.1 has been verified.

\vspace*{2mm}
Now, let $\phi(x,y)$ be a surface in $S^3$. Then, we say that $\{X^0_{\alpha}(x,y),X^0_{\beta}(x,y),\xi(x,y)\}$ is an orthonormal frame field of $\phi$, if 
$X^0_{\alpha}$ and $X^0_{\beta}$ are orthonormal tangent vector fields and $\xi$ is a unit normal vector field (in $S^3$) of $\phi$. 
For functions $\bar P(x,y)$ and $\bar P_z(x,y)$ on $V$, we denote $(e^{\bar P})_z(x,y):=(e^{\bar P}\bar P_z)(x,y)$ and 
$$(e^{\bar P}\cos\varphi)_z(x,y):=((e^{\bar P})_z\cos\varphi-e^{\bar P}\varphi_z\sin\varphi)(x,y),$$   
$$(e^{\bar P}\sin\varphi)_z(x,y):=((e^{\bar P})_z\sin\varphi+e^{\bar P}\varphi_z\cos\varphi)(x,y).$$

\vspace*{2mm}
{\pro 4.2}. \ \ {\it Let $\hat{g}\in Met^0$, and $(\varphi(x,y),\varphi_z(x,y))$ be functions on $V$ arising from $\hat{g}$. 
Let $\phi:V\ni(x,y)\mapsto \phi(x,y)\in S^3$ be a surface. 
For functions $\bar P(x,y)$ and $\bar P_z(x,y)$ on $V$, 
suppose that $\phi$ satisfies the following conditions (1), (2) and (3):
$$\phi_x=-{\bar \kappa}_1e^{\bar P}\cos\varphi\ X^0_{\alpha} \ \ \ \ 
\phi_y=-{\bar \kappa}_2e^{\bar P}\sin\varphi\ X^0_{\beta}, \leqno{ \ \ \ (1)}$$ 
$$\nabla'_{\partial/\partial y}X^0_{\alpha}=\frac{(e^{\bar P}\sin\varphi)_x}{e^{\bar P}\cos\varphi}~X_{\beta}^0, \ \ \ \ 
\nabla'_{\partial/\partial x}X^0_{\beta}=\frac{(e^{\bar P}\cos\varphi)_y}{e^{\bar P}\sin\varphi}~X_{\alpha}^0, \leqno{ \ \ \ (2)}$$ 
$$\nabla'_{\partial/\partial x}\xi=e^{-\bar P}(e^{\bar P}\cos\varphi)_z\ X_{\alpha}^0, \ \ \ \ 
\nabla'_{\partial/\partial y}\xi=e^{-\bar P}(e^{\bar P}\sin\varphi)_z\ X_{\beta}^0  \leqno{ \ \ \ (3)}$$
with respect to an orthonormal frame field $\{X^0_{\alpha}(x,y),X^0_{\beta}(x,y),\xi(x,y)\}$. Then, the Codazzi condition on $\phi$ is that $\bar P(x,y)$ and $\bar P_z(x,y)$ satisfy (3.9) and (3.10) defined by $(\varphi(x,y),\varphi_z(x,y))$. }\\

{\sc Proof}. \ Note that (1) and (3) in Proposition 4.1 are also satisfied by the assumption. Let $\lambda_1$ and $\lambda_2$ be the principal curvatures of $\phi$. Then, we have \ 
$$\lambda_1=(e^{\bar P}\cos\varphi)_z/(\bar\kappa_1e^{2\bar P}\cos\varphi), \hspace*{1cm}
\lambda_2=(e^{\bar P}\sin\varphi)_z/(\bar\kappa_2e^{2\bar P}\sin\varphi).$$

Now, the Codazzi condition on $\phi$ is given by the equation \ 
$\nabla'_{\partial/\partial y}\nabla'_{\partial/\partial x}\xi=\nabla'_{\partial/\partial x}\nabla'_{\partial/\partial y}\xi$, \ which 
implies that $\bar P(x,y)$ and $\bar P_z(x,y)$ satisfy (3.9) and (3.10).   
\hspace{\fill}$\Box$\\

\vspace*{2mm}
The following theorem gives a geometrical interpretation of the equations (3.8)$\sim$(3.11) for ${\bar P}(x,y)$ and 
${\bar P}_z(x,y)$: \\

\vspace*{2mm}
{\thm 2}. \ \ {\it Let $\hat{g}\in Met^0$, and $(\varphi(x,y),\varphi_z(x,y))$ be functions on $V$ arising from $\hat{g}$. 
Then, analytic functions ${\bar P}(x,y)$ and ${\bar P}_z(x,y)$ on $V$ satisfy \ $\bar\kappa_1\bar\kappa_2\neq 0$ \ and (3.8)$\sim$(3.11) defined by $(\varphi(x,y),\varphi_z(x,y))$, if and only if 
there is an analytic surface $\phi:V\ni(x,y)\mapsto \phi(x,y)\in S^3$ such that $\phi$ satisfies the following conditions (1), (2) and (3) on $V$: 
$$\phi_x=-{\bar \kappa}_1e^{\bar P}\cos\varphi\ X^0_{\alpha}, \ \ \ \ 
\phi_y=-{\bar \kappa}_2e^{\bar P}\sin\varphi\ X^0_{\beta},  \leqno{ \ \ \ (1)}$$$$\nabla'_{\partial/\partial y}X^0_{\alpha}=\frac{(e^{\bar P}\sin\varphi)_x}{e^{\bar P}\cos\varphi}~X_{\beta}^0, \ \ \ \ 
\nabla'_{\partial/\partial x}X^0_{\beta}=\frac{(e^{\bar P}\cos\varphi)_y}{e^{\bar P}\sin\varphi}~X_{\alpha}^0,  \leqno{ \ \ \ (2)}$$
$$\nabla'_{\partial/\partial x}\xi=e^{-\bar P}(e^{\bar P}\cos\varphi)_z\ X_{\alpha}^0, \ \ \ \ 
\nabla'_{\partial/\partial y}\xi=e^{-\bar P}(e^{\bar P}\sin\varphi)_z\ X_{\beta}^0   \leqno{ \ \ \ (3)}$$ 
with respect to an orthonormal frame field $\{X^0_{\alpha}(x,y),X^0_{\beta}(x,y),\xi(x,y)\}$.} \\

\vspace*{2mm}
{\sc Remark 4}. \ (1) In (2) and (3) of Theorem 2, we can replace the connection $\nabla'$ on $R^4$ with 
the canonical connection $D'$ on $S^3$, by Lemma 4.1. 

(2) Suppose that $({\bar P}(x,y),{\bar P}_z(x,y))$ satisfies $\bar\kappa_1\bar\kappa_2(x,y)\neq 0$ and (3.8)$\sim$(3.11) defined by $(\varphi(x,y),\varphi_z(x,y))$ arising from $\hat g$. Then, there is a generic conformally flat hypersurface $f(x,y,z)$ in $R^4$ determined from $({\bar P}(x,y),{\bar P}_z(x,y))$ by Theorem 1. For an orthonormal frame field $(X_{\alpha},X_{\beta},X_{\gamma},N)$ along $f$ defined in \S2.2, let us take \ $\phi(x,y):=N(x,y,0)$ \ as in Remark 3-(2) and put \ $X^0_{\alpha}(x,y):=X_{\alpha}(x,y,0)$, \ $X^0_{\beta}(x,y):=X_{\beta}(x,y,0)$ \ and \ $\xi(x,y):=X_{\gamma}(x,y,0)$ \ for $(x,y)\in V$. Then, $\{X^0_{\alpha}(x,y),X^0_{\beta}(x,y),\xi(x,y),\phi(x,y)\}$ and $(\bar P(x,y),\bar P_z(x,y))$ satisfy (1), (2) and (3), by the results in \S2.2.

(3) If two surfaces $\phi$ and $\bar\phi$ are isometric, we regard $\phi$ and $\bar\phi$ as the same surface. 
Then, for a pair $(\varphi(x,y),\varphi_z(x,y))$ arising from 
$\hat g\in Met^0$, all surfaces $\phi(x,y)$ obtained by Theorem 2 generate a 4-dimensional set. 
In fact, if a pair $(\bar P(x,y),\bar P_z(x,y))$ satisfies $\bar\kappa_1\bar\kappa_2\neq 0$ and (3.8)$\sim$(3.11), then so do all pairs $(\bar P(x,y)+c,\bar P_z(x,y))$ with constant $c$. 
On the other hand, every pair $(\bar P(x,y)+c,\bar P_z(x,y))$ with $c$ leads to the same equations in (1), (2) and (3), hence only one $\phi$ is determined for these pairs.  
The set of pairs $({\bar P}(x,y),\bar P_z(x,y))$ satisfying $\bar\kappa_1\bar\kappa_2\neq 0$ and (3.8)$\sim$(3.11) is 5-dimensional as mentioned in \S3, hence all surfaces $\phi$ determined from a pair $(\varphi(x,y),\varphi_z(x,y))$ generate a 4-dimensional set. 

In consequence, all surfaces $\phi(x,y)$ determined for $\hat g\in Met^0$ generate a 5-dimensional set, since $\hat g$ gives rise to a one-parameter family $(\varphi(x,y),\varphi^t_z(x,y))$ with $t\in R\setminus\{0\}$.  
\hspace{\fill}$\Box$\\

\vspace*{2mm}
For Theorem 2, we have only to verify that (3.11) for $(\bar P(x,y),\bar P_z(x,y))$: \ $(\bar\zeta-\bar\kappa_3^2)(x,y)=0$, \ is induced from the existence of $\phi$ satisfying (1), (2) and (3), since the other facts have been verified by Propositions 4.1, 4.2 and Remark 4-(2).   
To this end, we firstly provide the following lemma. \\

\vspace*{2mm}
{\lem 4.4}. \ \ {\it Let $\hat{g}^0\in Met^0$, and $(\varphi(x,y),\varphi_z(x,y))$ be functions on $V$ arising from $\hat{g}$. 
Let $\phi:V\ni(x,y)\mapsto \phi(x,y)\in S^3$ be a generic surface and  $(x,y)$ be a principal curvature line coordinate system. 
Suppose that $\phi$ and $(\bar P(x,y),\bar P_z(x,y))$ satisfy (1) and (2) in Theorem 2. 
Then, the Gauss curvature $K^{\phi}$ of the metric $I_{\phi}$ satisfies the following equation: 
$$\bar\kappa_1\bar\kappa_2\ K^{\phi}=
\frac{e^{-\bar P}}{\cos^2\varphi}\left((e^{-\bar P})_{xx}+\frac{\varphi_x}{\sin\varphi\cos\varphi}(e^{-\bar P})_x\right)  
\eqno{(4.7)}$$ 
$$+\frac{e^{-\bar P}}{\sin^2\varphi}\left((e^{-\bar P})_{yy}-\frac{\varphi_y}{\sin\varphi\cos\varphi}(e^{-\bar P})_y\right)
-\left(\frac{((e^{-\bar P})_x)^2}{\cos^2\varphi}+\frac{((e^{-\bar P})_y)^2}{\sin^2\varphi}\right)
-\frac{e^{-2\bar P}}{\sin\varphi\cos\varphi}\left(\varphi_{xx}-\varphi_{yy}\right).  $$ }\\

\vspace*{2mm}
{\sc Proof}. \ There is a surface $f^0(x,y)$ in $R^4$ satisfying Proposition 4.1-(1), by the assumption.  
When we express metrics $I_{\phi}$ and $I_{f^0}$ given by (4.4) and (4.5), respectively, as 
$$I_{\phi}=Edx^2+Gdy^2 \ \ \ {\rm and} \ \ \ I_{f^0}=E^0dx^2+G^0dy^2,$$
we have 
$$\frac{1}{\sqrt{E}}(\sqrt{G})_x=\frac{1}{\sqrt{E^0}}(\sqrt{G^0})_x, \ \ \ \ \ \ 
\frac{1}{\sqrt{G}}(\sqrt{E})_y=\frac{1}{\sqrt{G^0}}(\sqrt{E^0})_y.$$
In fact, we have  $(\bar{\kappa}_3)_x=-(e^{-\bar P})_x\tan\varphi$ and $(\bar{\kappa}_3)_y=(e^{-\bar P})_y\cot\varphi$
by Proposition 4.1. Then, 
since \ $(\bar{\kappa}_2e^{\bar P}\sin\varphi)_x=\bar{\kappa}_1(e^{\bar P}\sin\varphi)_x$ \ and \ 
$(\bar{\kappa}_1e^{\bar P}\cos\varphi)_y=\bar{\kappa}_2(e^{\bar P}\cos\varphi)_y$ \ 
hold, we obtain 
$$\frac{1}{\sqrt{E}}(\sqrt{G})_x=\frac{1}{\bar{\kappa}_1e^{\bar P}\cos\varphi}(\bar{\kappa}_2e^{\bar P}\sin\varphi)_x
=\frac{(e^{\bar P}\sin\varphi)_x}{e^{\bar P}\cos\varphi}=\frac{1}{\sqrt{E^0}}(\sqrt{G^0})_x,$$
$$\frac{1}{\sqrt{G}}(\sqrt{E})_y
=\frac{1}{\bar{\kappa}_2e^{\bar P}\sin\varphi}(\bar{\kappa}_1e^{\bar P}\cos\varphi)_y
=\frac{(e^{\bar P}\cos\varphi)_y}{e^{\bar P}\sin\varphi}=\frac{1}{\sqrt{G^0}}(\sqrt{E^0})_y.$$
Hence, $K^{\phi}=K^{f^0}/(\bar{\kappa}_1\bar{\kappa}_2)$ holds, further the Gauss curvature $K^{f^0}$ of $I_{f^0}$ 
is given by the right hand side of (4.7). 
\\ \hspace{\fill}$\Box$\\

\vspace*{2mm}
{\sc Proof of Theorem 2}. \ 
By the assumptions (1) and (3), two principal curvatures of $\phi$ are given by 
$$\frac{(e^{\bar P}\cos\varphi)_z}{{\bar\kappa}_1e^{2\bar P}\cos\varphi}, \ \ \ \ \ \ 
\frac{(e^{\bar P}\sin\varphi)_z}{{\bar\kappa}_2e^{2\bar P}\sin\varphi}.$$
Hence, by the Gauss equation we have 
$$\bar{\kappa}_1\bar{\kappa}_2\ K^{\phi}=\bar{\kappa}_1\bar{\kappa}_2+\frac{e^{-4\bar P}}{\sin\varphi\cos\varphi}
(e^{\bar P}\cos\varphi)_z(e^{\bar P}\sin\varphi)_z.$$ 
We replace the left hand side of this equation with the right hand side of (4.7). (Then, we obtain 
the first equation of (2.2.7) on $z=0$ by putting $K(X_{\alpha}\wedge X_{\beta})=\bar\kappa_1\bar\kappa_2$.)   
Furthermore, by Definition and Notation 3 in \S3, we have \ 
$$\bar{\kappa}_1\bar{\kappa}_2={\bar\kappa}_3^2-e^{-\bar P}[e^{-\bar P}+(\cos2\varphi/\sin\varphi\cos\varphi)\bar\kappa_3], $$ 
where, as $e^{-\bar P}$ can be expressed in the form of Proposition 2.2-(v), 
$$e^{-\bar P}+\frac{\cos2\varphi}{\sin\varphi\cos\varphi}\bar\kappa_3
=-\tan^2\varphi(e^{-\bar P})_{xx}-\varphi_x \left[\frac{1}{\sin\varphi\cos^3\varphi}-2\tan\varphi \right](e^{-\bar P})_x
-\cot^2\varphi(e^{-\bar P})_{yy}  $$
$$+\varphi_y \left[\frac{1}{\sin^3\varphi\cos\varphi}-2\cot\varphi \right](e^{-\bar P})_y
+(e^{-\bar P})_{zz}+\frac{\varphi_{xx}-\varphi_{yy}}{\sin\varphi\cos\varphi}e^{-\bar P}
-\frac{\cos2\varphi}{\sin\varphi\cos\varphi}\varphi_z(e^{-\bar P})_z$$
holds. The equation $(\bar\zeta-\bar\kappa_3^2)(x,y)=0$ follows directly from these equations. 

In consequence, the proof of Theorem 2 has been completed.  
\hspace{\fill}$\Box$\\

\vspace*{2mm}
Let $\hat g\in Met^0$, and $(\varphi(x,y),\varphi_z(x,y))$ be functions on $V$ arising from $\hat g$, as in Theorems 1 and 2. 
Let $({\bar P}(x,y),{\bar P}_z(x,y))$ be analytic functions on $V$ satisfying $\bar\kappa_1\bar\kappa_2\neq 0$ and (3.8)$\sim$(3.11) defined by $(\varphi(x,y),\varphi_z(x,y))$.
Let $\phi(x,y)$ be a surface in $S^3$ determined from $({\bar P}(x,y),{\bar P}_z(x,y))$ by Theorem 2. 
Then, we have the following fact: an evolution $\phi^z(x,y)$, $z\in I$, of surfaces in $S^3$ issuing from $\phi^0=\phi$ is uniquely determined such that  
$\phi^z$ leads to an evolution $f^z(x,y)$ of surfaces in $R^4$ and 
$f(x,y,z):=f^z(x,y)$ is a generic conformally flat hypersurface $f$ in Theorem 1.

In fact, let $({\bar P}(x,y),{\bar P}_z(x,y))$ and $\phi$ be as above. Then, three elements: $g=\hat g\in (Met^0\subset)Met$, \ $\varphi(x,y,z)$ arising from $g$ and $e^{-P}(x,y,z)\in \Xi_g$ satisfying $e^{-P}(x,y,0)=e^{-\bar P}(x,y)$ and $(e^{-P})_z(x,y,0)=(e^{-\bar P})_z(x,y)$, are determined for $(x,y,z)\in V\times I$. 
Then, for each $z\in I$, a pair $(P(x,y,z),P_z(x,y,z))$ satisfies (3.8)$\sim$(3.11) defined by $(\varphi(x,y,z),\varphi_z(x,y,z))$. 
Let $\kappa_3(x,y,z)$ be a function defined from $e^{-P}$ by Proposition 2.2-(iv) and the other $\kappa_i(x,y,z) \ (i=1,2)$ be defined by \ $\kappa_1=e^{-P}\tan\varphi+\kappa_3$, \ $\kappa_2=-e^{-P}\cot\varphi+\kappa_3$. \  
Then, for each $z\in I$ a surface $\phi^z(x,y)$ in $S^3$ is determined such that it satisfies (1), (2) and (3) in Theorem 2. Because we can assume $\kappa_1\kappa_2\neq 0$ on $U=V\times I$ by the condition $\bar\kappa_1\bar\kappa_2\neq 0$ on $V$. The differential of $\phi^z$ with $z$ is expressed as
$$(d\phi^z)(x,y)=-(\kappa_1e^P\cos\varphi)(x,y,z)X_{\alpha}(x,y,z)dx-(\kappa_2e^P\sin\varphi)(x,y,z)X_{\beta}(x,y,z)dy$$
with respect to an orthonormal frame field $\{X_{\alpha}(x,y,z),X_{\beta}(x,y,z)\}$.

On the other hand, a generic conformally flat hypersurface $f(x,y,z)$ in $R^4$ is also determined from $\varphi(x,y,z)$ and $e^{-P}(x,y,z)$ by Theorem 1, which has the Guichard net $g$ and $I_f=e^{2P}g$. Hence, we have an evolution $\hat\phi^z(x,y):=N(x,y,z)$ of surface as in Remark 3-(2), where $N$ is a unit normal vector field of $f$. Then, for each $z$, $\phi^z$ and $\hat\phi^z$ are isometric, since $\phi^z$ and $\hat\phi^z$ satisfy the same conditions (1), (2) and (3) in Theorem 2 by (2.2.3) and (2.2.4).  
Therefore, from the other equations in (2.2.3) and (2.2.4), we determine $\phi^z(x,y)$ by the following system of evolution equations in $z$ under the initial condition $\phi^0(x,y)=\phi(x,y)$: 
$$ 
\begin{array}{l}
(\phi^z)_z=-{\kappa}_3e^P\xi^z, \ \ \ \ \ (X_{\alpha})_z=\displaystyle\frac{P_x}{\cos\varphi}\xi^z, \ \ \ \ \ \ 
(X_{\beta})_z=\displaystyle\frac{P_y}{\sin\varphi}\xi^z,\\[3mm]
(\xi^z)_z=\displaystyle-{\frac{P_x}{\cos\varphi}X_{\alpha}-\frac{P_y}{\sin\varphi}X_{\beta}}
+{\kappa}_3e^P\phi^z
\end{array}
\eqno{(4.8)}
$$
under the initial condition \ 
$$
\begin{array}{l}
\phi^0(x,y)=\phi(x,y), \ \ X_{\alpha}(x,y,0)=X_{\alpha}^0(x,y),\\[2mm] 
X_{\beta}(x,y,0)=X_{\beta}^0(x,y), \ \ \xi^0(x,y)=\xi(x,y), 
\end{array}
\eqno{(4.9)}
$$
where $\xi^z(x,y)$ (resp. $\xi(x,y)$) be a unit normal vector field (in $S^3$) of $\phi^z(x,y)$ (resp. $\phi$).
Note that the conditions (1), (2) and (3) in Theorem 2 for $\phi^z$ are compatible with (4.8), since $\{X_{\alpha},X_{\beta},\xi^z,\phi^z\}$ is the orthonormal frame field of $f$ given in \S2.2 as $X_{\gamma}(x,y,z):=\xi^z(x,y)$ and $N(x,y,z):=\phi^z(x,y)$ and the frame field satisfies (2.2.3) and (2.2.4). Therefore, for each $z\in I$, the solution $\{X_{\alpha},X_{\beta},\xi^z,\phi^z\}$ to (4.8) under the initial condition (4.9) satisfies (1), (2) and (3) in Theorem 2, by the uniqueness of solutions. 

Next, by Theorem 2 and Proposition 4.1, for all $z\in I$,  we define surfaces $f^z(x,y)$ in $R^4$ by
$$\begin{array}{l}
(df^z)(x,y)=e^{P(x,y,z)}[\cos\varphi(x,y,z)X_{\alpha}(x,y,z)dx+\sin\varphi(x,y,z)X_{\beta}(x,y,z)dy],   \\[3mm]
(f^z)_z(x_0,y_0)=e^P(x_0,y_0,z)\xi^z(x_0,y_0) \hspace*{0.5cm} and \hspace*{0.5cm} f(x_0,y_0,0)=0,
\end{array} \eqno{(4.10)}$$
where the second equation is an ordinary differential equation in $z$ and $(x_0,y_0)$ is a fixed point of $V$. 
Then, $f(x,y,z):=f^z(x,y)$ is a generic conformally flat hypersurface with the Guichard net $g$. In fact, by $e^{-P}(x,y,z)\in \Xi_g$ there is a generic conformally flat hypersurface, say $f$, such that $\{X_{\alpha},X_{\beta},\xi^z,\phi^z\}$ is the orthonormal frame field along $f$. Both equations in (4.10) imply that $f^z(x,y)$ has the orthonormal frame field $\{X_{\alpha},X_{\beta},\xi^z,\phi^z\}$ and the first fundamental form $e^{2P}g$, since \ $(f^z)_z(x,y)=e^P(x,y,z)\xi^z(x,y)$ \ for any $(x,y)\in V$ is also satisfied by the existence of $f$ and the uniqueness of solutions.

\vspace*{2mm}
By the argument above, we have verified the following theorem. \\

\vspace*{2mm}   
{\thm 3}. \ \ {\it Let $\hat g\in Met^0$, and $(\varphi(x,y),\varphi_z(x,y))$ be functions on $V$ arising from ${\hat g}$. 
Suppose that $\phi:V\ni(x,y)\mapsto \phi(x,y)\in S^3$ is an analytic surface and $(\bar P(x,y),\bar P_z(x,y))$ are 
analytic functions on $V$ 
such that $\phi$ and $({\bar P}(x,y),{\bar P}_z(x,y))$ satisfy the conditions (1), (2), (3) in Theorem 2. 
Then, a generic conformally flat hypersurface $f$ with the Guichard net $g=\hat g\in (Met^0\subset)Met$ is realized in $R^4$ 
via an evolution of surfaces in $S^3$ issuing from $\phi$. 

More precisely, we have the following facts (1), (2) and (3) under the assumption above:

(1) Let $\varphi(x,y,z)$ be a function arising from $g$ by (2.1.1) and 
$e^{-P}(x,y,z)\in \Xi_g$ be a solution to (2.3.6) satisfying $e^{-P}(x,y,0)=e^{-\bar P}(x,y)$ \ and \ 
$(e^{-P})_z(x,y,0)=(e^{-\bar P})_z(x,y)$, where $(x,y,z)\in U:=V\times I$. 
Let $\kappa_3(x,y,z)$ be a function defined from $e^{-P}(x,y,z)$ by Proposition 2.2-(iv).
Then, orthonormal $R^4$-valued functions $X_{\alpha}(x,y,z)$, $X_{\beta}(x,y,z)$, $\xi^z(x,y)$, $\phi^z(x,y)$ 
are obtained as solutions to (4.8) under the initial condition (4.9). 
In particular, $\phi^z(x,y)$ is an evolution of surfaces issuing from $\phi^0=\phi$.  

(2) For the orthonormal frame $(X_{\alpha}(x,y,z),X_{\beta}(x,y,z),\xi^z(x,y),\phi^z(x,y))$ of $(x,y,z)\in U$ to $R^4$, an evolution $f^z(x,y)$ of surfaces in $R^4$ is determined by (4.10).

(3) A generic conformally flat hypersurface $f$ in $R^4$ with the Guichard net $g$ is determined from the evolution $f^z(x,y)$ by $f(x,y,z):=f^z(x,y)$. Then, $f$ has the first fundamental form $I_f=e^{2P}g$, unit normal vector field $\phi^z(x,y)$ and principal curvatures $\kappa_i \ (i=1,2,3)$, where \ $\kappa_1=e^{-P}\tan\varphi +\kappa_3$ \ and \ $\kappa_2=-e^{-P}\cot\varphi+\kappa_3$.  
Furthermore, $X_{\alpha}$, $X_{\beta}$, $\xi^z$ are unit principal curvature vector fields corresponding to the coordinate lines $x$, $y$ and $z$, respectively.  }\\

\newpage
\vspace{2mm}
\noindent
{\bf \large  5. Surfaces in $S^3$ leading to generic conformally flat hypersurfaces.}

In this section, we firstly study the condition for surfaces in $S^3$ to give rise to generic conformally flat hypersurfaces in $R^4$. 
Next, we show that, if a surface $\phi$ in $S^3$ leads to a generic conformally flat hypersurface $f$, then the dual (generic conformally flat) hypersurface $f^*$ of $f$ is also induced from $\phi$.

Let $\phi:V\ni (x,y)\mapsto \phi(x,y)\in S^3$ be a generic analytic surface and $( x, y)$ be a principal curvature line coordinate system. Let $\lambda_i$ $(i=1,2)$ be 
the principal curvatures of $\phi$, and $X^0_{\alpha}$ and $X^0_{\beta}$ be orthonormal principal vector fields corresponding to $\lambda_1$ and $\lambda_2$, respectively. 
Let $\xi$ be a unit normal vector field (in $S^3$) of $\phi$.  
Then, we can express derivatives of $\phi$, $\xi$, $X^0_{\alpha}$ and $X^0_{\beta}$ as follows: \\[2mm]

$\left\{\begin{array}{ll}
d\phi=-( a_1d x)X^0_{\alpha}-( a_2d y)X^0_{\beta}, &
D'\xi=( b_1d x)X^0_{\alpha}+( b_2d y)X^0_{\beta}, \\[3mm]
\nabla'_{\partial/\partial  x}X^0_{\alpha}=- c_1X^0_{\beta}- b_1\xi+ a_1\phi, &
\nabla'_{\partial/\partial  y}X^0_{\beta}=- c_2X^0_{\alpha}- b_2\xi+ a_2\phi,   \\[3mm]
\nabla'_{\partial/\partial  x}\xi=D'_{\partial/\partial  x}\xi=
 b_1X^0_{\alpha}, &
\nabla'_{\partial/\partial  y}\xi=D'_{\partial/\partial  y}\xi
= b_2X^0_{\beta}, \\[3mm]
\nabla'_{\partial/\partial  x}X^0_{\beta}=D'_{\partial/\partial  x}X^0_{\beta}=D_{\partial/\partial  x}X^0_{\beta}=
 c_1X^0_{\alpha}, &
\nabla'_{\partial/\partial  y}X^0_{\alpha}=D'_{\partial/\partial  y}X^0_{\alpha}=D_{\partial/\partial  y}X^0_{\alpha}
= c_2X^0_{\beta}, 
\end{array}\right.  \hfill{(5.1)}$ \\[2mm]
where $ a_i, \  b_i, \  c_i$ are functions on $V$. 
In fact, the equations in the last two lines of (5.1) follow from Lemma 4.1. 
In particular, we have $\lambda_i= b_i/ a_i \ (i=1,2). $ 

The following lemma is fundamental for such a surface $\phi$:\\

\vspace*{2mm}
{\lem 5.1}. \ \ {\it We have \ 
$$c_1=(a_1)_{ y}/ a_2=( b_1)_{ y}/ b_2, \ \ \ \ 
 c_2=( a_2)_{ x}/ a_1=( b_2)_{ x}/ b_1.$$}

\vspace*{2mm}
{\proof} \ We have 
$$0=d^2\phi=\left\{[( a_1)_{ y}X^0_{\alpha}- a_2\nabla'_{\partial/\partial  x}X^0_{\beta}]
-[( a_2)_{ x}X^0_{\beta}- a_1\nabla'_{\partial/\partial  y}X^0_{\alpha}]\right\}d x\wedge d y$$ 
by $d=\nabla'$. Hence, we obtain 
$$\nabla'_{\partial/\partial  x}X^0_{\beta}=( a_1)_{ y}/ a_2~X^0_{\alpha}, \ \ \ 
\nabla'_{\partial/\partial  y}X^0_{\alpha}=( a_2)_{ x}/ a_1~X^0_{\beta}$$
by the last two equations of (5.1). 
In the same way, from (5.1) and $d^2\xi=0$, we obtain 
$$\nabla'_{\partial/\partial  x}X^0_{\beta}=( b_1)_{ y}/ b_2~X^0_{\alpha}, \ \ \ 
\nabla'_{\partial/\partial  y}X^0_{\alpha}=( b_2)_{ x}/ b_1~X^0_{\beta}.$$
\hspace{\fill}$\Box$\\

Two functions $( a_1)_{ x}$ and $( a_2)_{ y}$ might have some special property, if a generic analytic surface $\phi(x,y)$ gives rise to a generic conformally flat hypersurface in $R^4$. Because these functions do not appear in Lemma 5.1. 

Now, in the following lemma, let $f:V\times I\ni(x,y,z)\mapsto f(x,y,z)\in R^4$ be a generic conformally flat hypersurface with the Guichard net $g$; we denote by $I_f=e^{2P}g$, $N(x,y,z)$ and $\kappa_i(x,y,z)$ the first fundamental form, a unit normal vector field and the principal curvatures of $f$, respectively, and by $\varphi(x,y,z)$ a function arising from $g$ by (2.1.1). 
In particular, we have defined orthonormal tangent vector fields $X_{\alpha}$ and $X_{\beta}$ of $f$ by \ $X_{\alpha}=\frac{e^{-P}}{\cos\varphi}f_x$ \ and \ $X_{\beta}=\frac{e^{-P}}{\sin\varphi}f_y$ \ in \S2.2.
If $(\kappa_1\kappa_2)(x,y,0)\neq 0$ for $(x,y)\in V$, then \ $\phi(x,y):=N(x,y,0)$ \ is a generic surface in $S^3$. We denote \ $\bar P(x,y):=P(x,y,0)$, \ $\bar\kappa_i(x,y):=\kappa_i(x,y,0)$, \ $\varphi(x,y):=\varphi(x,y,0)$, \ $X^0_{\alpha}(x,y):=X_{\alpha}(x,y,0)$ and $X^0_{\beta}:=X_{\beta}(x,y,0)$, \ as in the previous sections. \\

\vspace*{2mm}
{\lem 5.2}. \ \ {\it Let $\phi:V\ni ( x, y)\mapsto \phi( x, y)\in S^3$ be a generic surface, as above. Then, $a_i(x,y)$, $\varphi( x, y)$ and $\bar P( x, y)$ satisfy the following conditions (1) and (2):\\ 

(1) \ $ a_1\sin\varphi-a_2\cos\varphi=1$. \hspace*{1cm}
(2) \ $( a_1)_{ x}=-\varphi_{ x} a_2+\bar P_{ x}a_1, \ \ \ \ 
( a_2)_{ y}=\varphi_{ y} a_1+\bar P_{ y}a_2$. } \\

\vspace*{2mm}
{\sc Proof}. \ Since $\kappa_1$ and $\kappa_2$, respectively, are the first and the second principal curvatures of $f$, we have
$$ a_1=\bar \kappa_1e^{\bar P}\cos\varphi=(\sin\varphi+\bar\kappa_3e^{\bar P}\cos\varphi), \ \ \ \ 
 a_2=\bar \kappa_2e^{\bar P}\sin\varphi=(-\cos\varphi+\bar\kappa_3e^{\bar P}\sin\varphi).$$
Hence, (1) is satisfied.

Next, $(\bar \kappa_3)_{ x}=-(e^{-\bar P})_{ x}\tan\varphi$ \ and \ $(\bar \kappa_3)_{ y}=(e^{-\bar P})_{ y}\cot\varphi$ \ hold 
by (i) and (ii) of Proposition 2.2. 
From these equations and (2.2.2), we have \ $(\bar \kappa_3e^{\bar P})_{ x}=\bar \kappa_1(e^{\bar P})_{ x}$ \ and \ 
$(\bar \kappa_3e^{\bar P})_{ y}=\bar \kappa_2(e^{\bar P})_{ y}$. \ 
Then, we obtain (2). 
\hspace{\fill}$\Box$\\

\vspace*{2mm}
In Lemma 5.2, we can not distinguish between $a_i$ and $-a_i$ by the surface $\phi$ itself. 
If we replace only the equation (1) with \ $ a_1\sin\varphi-a_2\cos\varphi=\pm 1,$ \ then these conditions (1) and (2) are independent of the choice of $\pm a_i$.  
For example, suppose that we chose $(a_1,-a_2)$ in place of $(a_1,a_2)$. Then, we can take $-\varphi$ as $\varphi$ in (2), then (1) changes to \ $a_1\sin(-\varphi)-(-a_2)\cos(-\varphi)=-1.$ This change $\varphi\rightarrow -\varphi$ corresponds to the change $N\rightarrow -N$ for $f$.

\vspace*{2mm}
Now, we have the following theorem. In the theorem, we denote by $\kappa_3$ the principal curvature in the direction of $z$ for a hypersurface in $R^4$, as usual. \\

\vspace*{2mm}
{\thm 4}. \ \ {\it Let $\phi:V\ni (x,y)\mapsto \phi(x,y)\in S^3$ be a generic analytic surface and $(x,y)$ be a principal curvature line coordinate system.  
Let $a_i$ be functions in (5.1). 
Suppose that, for $\phi$, there is a coordinate system $(x,y)$ and analytic functions $(\varphi(x,y),\bar P(x,y))$ on $V$ such that they satisfy the following conditions (1) and (2): 
$$(1) \ a_1\sin\varphi- a_2\cos\varphi=1 \ {\rm (} resp. \ -1{\rm )}, \ \ \ (2) \ (a_1)_{ x}=-\varphi_{ x} a_2+\bar P_{ x} a_1, \ \ 
 ( a_2)_{ y}=\varphi_{ y} a_1+\bar P_{ y} a_2.$$  
Then, an evolution $\phi^z(x,y)$, $z\in I$, of surfaces in $S^3$ issuing from $\phi(x,y)$ is uniquely determined such that $\phi^z(x,y)$ is the Gauss map of a generic conformally flat hypersurface $f(x,y,z)$ of $V\times I$ into $R^4$. 
In particular, $\kappa_3$ is middle among three principal curvatures of $f$. \\ }

\vspace*{2mm}
{\sc Remark 5}. \ (1) In Theorem 4, $\bar P( x, y)$ is only determined up to a constant term.    
This fact implies that, for $(\phi(x,y),\varphi(x,y),(d\bar P)(x,y))$, a generic conformally flat hypersurface $f$ in $R^4$ is determined uniquely up to homothety and parallel translation (see Remark 4-(3) in \S4).  

(2) By Theorem 4, we can judge by only the first fundamental form whether a surface $\phi$ leads to a generic conformally flat hypersurface or not.
Furthermore, it is possible that the Guichard net $g$ of $f$ is not the one arising from \ $\hat g\in Met^0$, as see in the following (3).   

(3) In almost all cases, $\bar P(x,y)$ would be really two-variables function. 
However, in some cases $\bar P(x,y)$ is a one-variable function or a constant function. In Example 1 of \S3, we can choose $X_1(x)\equiv 0$ and $Y_1(y)\neq 0$ such that $Y_1(y)$ satisfies $H(y)\equiv 0$. For these pairs $(X_1(x), Y_1(y))$, all three functions $(e^{-\bar P}, (e^{-\bar P})_z,\bar\kappa_3)$ are one-variable functions for $y$. 
As a case $g\notin Met$, let $g\in CFM$ be a cyclic Guichard net, of which $\varphi(x,y,z)$ satisfies \ $\varphi_{zx}=\varphi_{zy}\equiv 0$. All hypersurfaces $f$ with such Guichard nets are constructed from certain surfaces in either $R^3$, $S^3$ or the Hyperbolic 3-space $H^3$ (\cite{hs1}). Then, for surfaces $\phi$ arising from hypersurfaces constructed from surfaces in $H^3$, $\bar P(x,y)$ are really two-variables functions.
But, for all hypersurfaces $f$ given in (\cite{su3}, \S2.2) (which are constructed from surfaces in either $R^3$ or $S^3$), $P(x,y,z)$ are one-variable functions for $z$. Hence, $\bar P(x,y)$ for surfaces $\phi$ arising from these $f$ are constant functions. On the other hand, for surfaces $\phi$ arising from their conformal transformations $\iota_q(f)$ with $q\neq 0$, $\bar P(x,y)$ are really two-variables functions. \hspace{\fill}$\Box$\\

In the proof of Theorem 4, for the sake of simplicity, we assume that the condition (1) is given by \ $a_1\sin\varphi-a_2\cos\varphi=1$ \ with respect to the expression \ $d\phi=-(a_1dx)X^0_{\alpha}-(a_2dy)X^0_{\beta}$ \ in (5.1). We firstly provide the following lemmata 5.3 and 5.4 for the proof:\\

\vspace*{2mm}
{\lem 5.3.} \ \ {\it Let $\phi:V\ni ( x, y)\mapsto \phi(x,y)\in S^3$ be a generic analytic surface. Suppose that the coordinate system $(x,y)$ and functions $(\varphi( x, y),\bar P( x, y))$ satisfy the assumption in Theorem 4. Then, we have the following facts (1), (2) and (3):\\

(1) There is a surface $f^0:V\ni ( x, y)\mapsto f^0( x, y)\in R^4$ satisfying \ $df^0=(e^{\bar P}\cos\varphi~d x)X^0_{\alpha}+(e^{\bar P}\sin\varphi~d y)X^0_{\beta}.$ \ 
In particular, $c_i$ in (5.1) are given by 
$ c_1=(e^{\bar P}\cos\varphi)_{ y}/(e^{\bar P}\sin\varphi)$ \ and \ 
$ c_2=(e^{\bar P}\sin\varphi)_{ x}/(e^{\bar P}\cos\varphi)$. \\

(2) There is an analytic function $Q( x, y)$ on $V$ such that \
$ a_1=\sin\varphi+Q\cos\varphi$, \ $ a_2=-\cos\varphi+Q\sin\varphi$ \ and \ $Q_{ x}=\bar P_{ x}(\tan\varphi+Q),$ \ $Q_{ y}=\bar P_{ y}(-\cot\varphi+Q)$ \ hold.
In particular, $\varphi$ and $\bar P$ satisfy (3.8). \\

(3) Let us define $\bar \kappa_i(x, y) \ (i=1,2,3)$ by \ $ a_1=\bar \kappa_1e^{\bar P}\cos\varphi$, \ $a_2=\bar \kappa_2e^{\bar P}\sin\varphi$ \ and \ 
$\bar \kappa_3:=e^{-\bar P}Q$. \ Then, \ $\bar \kappa_1=e^{-\bar P}\tan\varphi+\bar\kappa_3$, \ 
$\bar \kappa_2=-e^{-\bar P}\cot\varphi+\bar\kappa_3$ \ and \ 
$(\bar\kappa_3)_{ x}=-(e^{-\bar P})_{ x}\tan\varphi,$ \ 
$(\bar\kappa_3)_{ y}=(e^{-\bar P})_{ y}\cot\varphi$ \ hold. }\\

\vspace*{2mm}
{\sc Proof.} \ (1) Firstly, we differentiate \ $ a_1\sin\varphi-a_2\cos\varphi=1$ \ by $ x$ and $ y$, respectively. 
Then, by the assumption (2) in Theorem 4 and \ $( a_2)_{ x}= a_1 c_2$, \ $(a_1)_{ y}= a_2 c_1$ \ in Lemma 5.1, we have 
$${ c_1}={\bar P_{ y}\cot\varphi-\varphi_{ y}}
={(e^{\bar P}\cos\varphi)_{ y}}/(e^{\bar P}\sin\varphi), \ \ \ \  
{ c_2}={\bar P_{ x}\tan\varphi+\varphi_{ x}}
={(e^{\bar P}\sin\varphi)_{x}}/(e^{\bar P}\cos\varphi).$$
The existence of $f^0(x,y)$ follows from these equations, by the definitions of $c_i$. 

(2) The condition \ $a_1\sin\varphi- a_2\cos\varphi=1$ \ implies the existence of $Q( x, y)$ such that \ 
$ a_1=\sin\varphi+Q\cos\varphi$ \ and \ $ a_2=-\cos\varphi+Q\sin\varphi$. \ Then, since \ 
$( a_1)_{ x}=-\varphi_{ x} a_2+Q_{ x}\cos\varphi,$ \ 
we have \ $Q_{ x}=(\bar P_{ x}/\cos\varphi) a_1=\bar P_{ x}(\tan\varphi+Q).$ \  In the same way, we have \ $Q_{ y}=(\bar P_{ y}/\sin\varphi) a_2
=\bar P_{ y}(-\cot\varphi+Q).$ 
Hence, if $\bar P(x,y)$ is a one-variable function, then so are $Q(x,y)$ and $\varphi(x,y)$ for the same variable, by \ $0=Q_{xy}=\bar P_x\varphi_y/\cos^2\varphi$\ and \ $0=Q_{yx}=\bar P_y\varphi_x/\sin^2\varphi$. \ If $\bar P(x,y)$ is constant, then so is $Q(x,y)$. 
The last statement follows from the integrability condition on $Q$. 

(3) Note that all $\bar \kappa_i( x, y) \ (i=1,2,3)$ are well-defined, since the other functions in each equation are determined, in particular, \ $Q=a_1\cos\varphi+a_2\sin\varphi$. \ Then, all equations follows from (2) directly. 
Note that, if $\bar P(x,y)$ is a one-variable (or constant) function, then so is $\bar\kappa_3(x,y)$. \\ 
\hspace{\fill}$\Box$\\

Now, under the assumption in Theorem 4, for $b_i$ in (5.1), we define functions $\varphi_z(x,y)$ and $\bar P_z(x,y)$ on $V$ by 
$$ b_1=\bar P_z\cos\varphi-\varphi_z\sin\varphi(=:e^{-\bar P}(e^{\bar P}\cos\varphi)_z), \ \ 
 b_2=\bar P_z\sin\varphi+\varphi_z\cos\varphi(=:e^{-\bar P}(e^{\bar P}\sin\varphi)_z). \eqno{(5.2)}$$
Then, by \ $d^2\xi=0$ \ and Lemma 5.3-(1) we have 
$$
\begin{array}{l}
[(e^{-\bar P})_z]_x+(e^{-\bar P})_x\varphi_z\tan\varphi-e^{-\bar P}\varphi_{zx}\cot\varphi=0,  \\[1mm]
 
[(e^{-\bar P})_z]_y-(e^{-\bar P})_y\varphi_z\cot\varphi+e^{-\bar P}\varphi_{zy}\tan\varphi=0,
\end{array}
\eqno{(5.3)}
$$
where \ $(e^{-\bar P})_z:=-e^{-\bar P}\bar P_z$, \ $\varphi_{zx}:=(\varphi_z)_x$ \ and \ $\varphi_{zy}:=(\varphi_z)_y$. \  The equations in (5.3) show that $(e^{-\bar P})_z(x,y)$ and $\varphi_z(x,y)$ satisfy (3.9) and (3.10). 

Next, let us set 
$$\Pi:=\{\varphi(x,y),\varphi_z(x,y),e^{-\bar P}(x,y),(e^{-\bar P})_z(x,y),\bar\kappa_3(x,y)\}.$$
We say that a class \ $\{\tilde\varphi(x,y,z),e^{-\tilde P}(x,y,z),\tilde\kappa_3(x,y,z)\}$ \ of analytic functions on \ $V\times I$ \ is an extension of $\Pi$ if they satisfy, on $V$, the following conditions: $$\tilde\varphi(x,y,0)=\varphi(x,y), \ \ \ \tilde\varphi_z(x,y,0)=\varphi_z(x,y), \ \ \ e^{-\tilde P}(x,y,0)=e^{-\bar P}(x,y),$$  
$$(e^{-\tilde P})_z(x,y,0)=(e^{-\bar P})_z(x,y), \ \ \ \tilde\kappa_3(x,y,0)=\bar\kappa_3(x,y), \ \ \ (\tilde\kappa_3)_z(x,y,0)=-(e^{-\bar P}\varphi_z)(x,y).$$
For extensions $\tilde\varphi$ and $e^{-\tilde P}$, we have \ $[(e^{-\bar P})_z]_x(x,y)=(e^{-\tilde P})_{zx}(x,y,0)$, \ $(\varphi_z)_x(x,y)=\tilde \varphi_{zx}(x,y,0)$ \ and so on. 
Hence, $\tilde\varphi$ and $e^{-\tilde P}$ satisfy two equations in (5.3) on $V\times\{0\}$ with respect to their real derivatives.
From \ $(\bar\kappa_3)_x=-(e^{-\bar P})_x\tan\varphi$ \ and \ $(\bar\kappa_3)_y=(e^{-\bar P})_y\cot\varphi$ \ in Lemma 5.3-(3) and the condition \ $(\tilde\kappa_3)_z(x,y,0)=-(e^{-\tilde P}\tilde\varphi_z)(x,y,0)$, \ we have
$$
\begin{array}{l}
 (\tilde\kappa_3)_x(x,y,z)=-(e^{-\tilde P})_x\tan\tilde\varphi(x,y,z)+O(z^2), \\[1mm]

 (\tilde\kappa_3)_y(x,y,z)=(e^{-\tilde P})_y\cot\tilde\varphi(x,y,z)+O(z^2) 
\end{array}
\eqno{(5.4)}
$$ 
by \ $(\tilde\kappa_3)_{zx}(x,y,0)=(\tilde\kappa_3)_{xz}(x,y,0)$, \ $(\tilde\kappa_3)_{zy}(x,y,0)=(\tilde\kappa_3)_{yz}(x,y,0)$ \ and (5.3).   
Here, for a function $h(x,y,z)$ on \ $V\times I$, we denote \ $h(x,y,z)=O(z^k)$ \ if there is a constant $C$ such that \ $h(x,y,z)/z^k<C$ \ holds as $z\rightarrow 0$. 

In the proof of the following lemma, we use equations of Definition and Notation 3 in \S3, which are equations for functions on $V$. Hence, we can use them even if $(\varphi(x,y),\varphi_z(x,y))$ does not lead to any metric $\hat g\in Met^0$.  \\

{\sc Lemma 5.4}. \ {\it Let $\phi:V\ni ( x, y)\mapsto \phi(x,y)\in S^3$ be a generic analytic surface. Suppose that the coordinate system $(x,y)$ and functions $(\varphi( x, y),\bar P( x, y))$ satisfy the assumption in Theorem 4. 
Then, there is an extension \ $\{\tilde\varphi(x,y,z),e^{-\tilde P}(x,y,z),\tilde\kappa_3(x,y,z)\}$ \ of $\Pi$ such that $\tilde\kappa_3(x,y,z)$ and $e^{-\tilde P}(x,y,z)$ satisfy (iv) and (v) of Proposition 2.2, respectively, with respect to \ $\varphi:=\tilde\varphi$, \ $e^{-P}:=e^{-\tilde P}$, \ $\kappa_3:=\tilde\kappa_3$ \ and \ $\psi_{zz}:=\tilde\psi_{zz}$, \ where \ $\tilde\psi_{zz}(x,y,z)=[(\tilde\varphi_{xx}-\tilde\varphi_{yy}-\tilde\varphi_{zz}\cos2\tilde\varphi)/\sin2\tilde\varphi](x,y,z)$. \  

Furthermore, for the other such extensions $\hat\varphi(x,y,z)$ and $e^{-\hat P}(x,y,z)$, we have \ $\hat\varphi(x,y,z)=\tilde\varphi(x,y,z)+O(z^4)$ \ and \ $e^{-\hat P}(x,y,z)=e^{-\tilde P}(x,y,z)+O(z^4)$.   }\\

{\sc Proof}. \ Note that all functions of $\Pi$ are determined by Lemma 5.3 and (5.2). 
Firstly, we use Definition and Notation 3-(2). Then, we determine $\tilde\varphi(x,y,z)$ and $\tilde\psi_{zz}(x,y,z)$ from $\Pi$ in the following way. $\varphi_{zz}(x,y)$ is determined from the equation of $\bar\kappa_3(x,y)$. Hence, we choose an arbitrary $\tilde\varphi(x,y,z)$ such that it satisfies \ $\tilde\varphi(x,y,0)=\varphi(x,y),$ \ $\tilde\varphi_z(x,y,0)=\varphi_z(x,y)$ \ and \ $\tilde\varphi_{zz}(x,y,0)=\varphi_{zz}(x,y)$. \ For its $\tilde\varphi(x,y,z)$, we define $\tilde\psi_{zz}(x,y,z)$ by \ $\tilde\psi_{zz}(x,y,z):=[(\tilde\varphi_{xx}-\tilde\varphi_{yy}-\tilde\varphi_{zz}\cos2\tilde\varphi)/\sin2\tilde\varphi](x,y,z)$. \ Note that $(e^{-\bar P})_{zz}(x,y)$ is determined from $\Pi$ and $\psi_{zz}(x,y):=\tilde\psi_{zz}(x,y,0)$, that is, $(e^{-\bar P})_{zz}(x,y)$ is independent of the choice of $\tilde \varphi(x,y,z)$.   

Next, since Proposition 2.2-(v) and (2.3.6) are the same equation, let $e^{-\tilde P}(x,y,z)$ be a solution to (2.3.6) defined by $\varphi:=\tilde\varphi$ and $\psi_{zz}:=\tilde\psi_{zz}$ under the initial condition \ $e^{-\tilde P}(x,y,0)=e^{-\bar P}(x,y)$ \ and \ 
$(e^{-\tilde P})_z(x,y,0)=(e^{-\bar P})_z(x,y)$. Note that we have \ $(e^{-\tilde P})_{zz}(x,y,0)=(e^{-\bar P})_{zz}(x,y)$. 
Then, we determine $\tilde\kappa_3(x,y,z):=\kappa_3(x,y,z)$ by Proposition 2.2-(iv) from $\varphi(x,y,z):=\tilde \varphi(x,y,z)$ and $e^{-P}(x,y,z):=e^{-\tilde P}(x,y,z)$. Note that we have \ $\tilde \kappa_3(x,y,0)=\bar\kappa_3(x,y)$. \

Finally, we use the condition \ $(\tilde \kappa_3)_z(x,y,0)=-(e^{-\tilde P}\varphi_z)(x,y,0)$. \ Then, by Proposition 2.2-(iv), $\tilde \varphi_{zzz}(x,y,0)$ and $\tilde \psi_{zzz}(x,y,0)$ are necessarily determined from $\Pi$, $\tilde \varphi_{zz}(x,y,0)$ and $(e^{-\tilde P})_{zz}(x,y,0)$. Hence, it is necessary to replace the first $\tilde \varphi(x,y,z)$ with a new one satisfying its $\tilde \varphi_{zzz}(x,y,0)$. Then, $(e^{-\tilde P})_{zzz}(x,y,0)$ is also determined from $\Pi$ by (2.3.6).

In consequence, we have verified Lemma.  
\hspace{\fill}$\Box$\\

{\sc Remark 6}. \ In the proof above, $\varphi_{zz}(x,y)$ is determined from $\Pi$. For $\varphi_{zz}(x,y)$, we can define a function $(\psi_{xx}-\psi_{yy})(x,y)$ on $V$ by the following equation: 
$$\varphi_{zz}=(\varphi_{xx}-\varphi_{yy})\cos2\varphi+(\psi_{xx}-\psi_{yy})\sin2\varphi.$$
Hence, $(\psi_{xx}-\psi_{yy})(x,y)$ is also determined from $\Pi$, and $\psi_{zz}(x,y)$ is expressed as 
$$\psi_{zz}=(\varphi_{xx}-\varphi_{yy})\sin2\varphi-(\psi_{xx}-\psi_{yy})\cos2\varphi.$$
\hspace{\fill}$\Box$\\

\vspace*{2mm}
{\sc Proof of Theorem 4.} \ Let $\phi( x, y)$ be a generic analytic surface in $S^3$. Suppose that the coordinate system $(x,y)$ and functions $(\varphi( x, y), \bar P(x, y))$ on $V$ satisfy the assumption. 
Then, we have the function $\bar\kappa_3( x, y):=e^{-\bar P}Q$ on $V$ by Lemma 5.3-(3) and functions $\varphi_z(x,y)$ and $(e^{-\bar P})_z(x,y)$ on $V$ by (5.2). Hence, all functions of $\Pi$ are determined. Then, $\bar\kappa_3( x, y)$ satisfies \ $(\bar\kappa_3)_x=-(e^{-\bar P})_x\tan\varphi$ \ and \ $(\bar\kappa_3)_y=(e^{-\bar P})_y\cot\varphi$ \ by Lemma 5.3-(3). Furthermore, $\varphi_{zz}(x,y)$ and $\psi_{zz}(x,y)$ are determined from $\Pi$, and then   
$\bar\kappa_3( x, y)$ and $(e^{-\bar P})_{zz}( x, y)$ satisfy two equations of Definition and Notation 3-(2), as verified in Lemma 5.4. 

Now, let $\{\tilde\varphi(x,y,z),e^{-\tilde P}(x,y,z),\tilde\kappa_3(x,y,z)\}$ be an extension of $\Pi$ in Lemma 5.4. 
We prove the theorem step by step. 

{\sc Step 1}. We have the following facts (A) and (B): 

(A) For \ $\bar \kappa_1:=e^{-\bar P}\tan\varphi+\bar\kappa_3$, \ $\bar\kappa_2:=-e^{-\bar P}\cot\varphi+\bar\kappa_3$, \ we have \ $\phi_x=-\bar\kappa_1e^{\bar P}\cos\varphi~X^0_{\alpha}, \ \ \phi_y=-\bar\kappa_2e^{\bar P}\sin\varphi~X^0_{\beta}$, \ by Lemma 5.3-(3).  
We have 
$$\nabla'_{\partial/\partial y}X^0_{\alpha}=\frac{(e^{\bar P}\sin\varphi)_x}{e^{\bar P}\cos\varphi}X^0_{\beta}, \ \ \ \ 
\nabla'_{\partial/\partial x}X^0_{\beta}=\frac{(e^{\bar P}\cos\varphi)_y}{e^{\bar P}\sin\varphi}X^0_{\alpha}$$
by Lemma 5.3-(1). 
Since $\varphi_z(x,y)$ and $\bar P_z(x,y)$ satisfy (5.2), we have 
$$\nabla'_{\partial/\partial x}\xi=e^{-\bar P}(e^{\bar P}\cos\varphi)_zX^0_{\alpha}, \ \ \ \ 
\nabla'_{\partial/\partial y}\xi=e^{-\bar P}(e^{\bar P}\sin\varphi)_zX^0_{\beta}.  $$

(B) Let $\tilde\zeta(x,y,z):=\zeta(x,y,z)$ be a function defined by (2.3.8) from $\varphi:=\tilde\varphi$ and $e^{-P}:=e^{-\tilde P}$. With respect to $\varphi:=\tilde\varphi$, $e^{-P}:=e^{-\tilde P}$ and $\kappa_3:=\tilde\kappa_3$, we define $I^{x,y}(x,y,z)$, $I^{x,z}(x,y,z)$, $I^{y,z}(x,y,z)$ and \ $J(x,y,z):=(\tilde\zeta-\tilde\kappa_3^2)(x,y,z)$ by (2.3.1) and Definition and Notation 2 in \S2.3, respectively. Then, we have
$$I^{x,y}(x,y,0)=I^{x,z}(x,y,0)=I^{y,z}(x,y,0)=J(x,y,0)=0. \eqno{(5.5)}$$
In fact, functions of $\Pi$ satisfy the equations (3.8)$\sim$(3.10), that is, the first three equations, by Lemma 5.3-(2) and (5.3). Furthermore, for $\bar\zeta(x,y)$ defined from $\Pi$ by Definition and Notation 3-(4), we have \ $J(x,y,0)=(\bar\zeta-\bar\kappa_3^2)(x,y)$. \ Then, \ $\bar\zeta-\bar\kappa_3^2\equiv 0$ \ is the Gauss equation for $\phi$, as verified in Theorem 2.

\vspace*{2mm}
By Theorem 3 and (A), if $(\varphi(x,y),\varphi_z(x,y))$ determines a 2-metric $\hat g\in Met^0$, then $\phi(x,y)$ leads to a generic conformally flat hypersurface. However, in this proof we directly show that the extension \\ $\{\tilde\varphi(x,y,z),e^{-\tilde P}(x,y,z),\tilde\kappa_3(x,y,z)\}$ determines an analytic metric $g\in  CFM$ ($g\notin Met$ is possible).

\vspace*{1mm}
{\sc Step 2}. Each analytic metric $g\in CFM$ corresponds to a pair $(\psi(x,y,z),\varphi(x,y,z))$ of functions on $V\times I$, where $(\psi(x,y,z),\varphi(x,y,z))$ is a solution to the system of evolution equations (iii) and (iv) in Proposition 2.1-(3). In order to obtain the solution, we have to find out other suitable initial data $(\psi(x,y),\psi_z(x,y))$ in addition to $(\varphi(x,y),\varphi_z(x,y))$ of $\Pi$.  
In the following step 3, we study the extension $\{\tilde\varphi(x,y,z),e^{-\tilde P}(x,y,z),\tilde\kappa_3(x,y,z)\}$ in detail, particularly, around $V\times\{0\}$. From the result, we determine these initial data in Step 4. 

Before Step 3, we remember equations (2.3.3)$\sim$(2.3.5), (3.1) and (3.5)$\sim$(3.7) obtained in the proofs of Lemma 2.1 and Proposition 3.1:  
with respect to $\varphi:=\tilde\varphi$, $e^{-P}:=e^{-\tilde P}$, $\kappa_3:=\tilde\kappa_3$ and $\psi_{zz}:=\tilde\psi_{zz}$, all these equations are satisfied, since $\tilde \kappa_3(x,y,z)$ and $e^{-\tilde P}(x,y,z)$ satisfy (iv) and (v) of Proposition 2.2, respectively.  

\vspace*{1mm}
{\sc Step 3}. In this step, let $I^{x,y}$, $I^{x,z}$, $I^{y,z}$, $J$ be functions defined in Step 1-(B). 
All $z$-derivatives $(I^{x,y})_z(x,y,0)$, $(I^{x,z})_z(x,y,0)$, $(I^{y,z})_z(x,y,0)$ and $J_z(x,y,0)$ are independent of the choice of extension $(\tilde \varphi(x,y,z),e^{-\tilde P}(x,y,z))$, by Lemma 5.4. 

Then, we have the following equations: 
$$(I^{x,y})_z(x,y,0)=(I^{x,z})_z(x,y,0)=(I^{y,z})_z(x,y,0)=J_z(x,y,0)=0 \hspace*{0.5cm}{\rm on} \ \ V\times\{0\}.  \eqno{(5.6)}$$
In fact, we have \ $(e^{-\tilde P})_{zx}=-(e^{-\tilde P})_x\tilde \varphi_z\tan\tilde \varphi+e^{-\tilde P}\tilde \varphi_{zx}\cot\tilde \varphi$ \ and \ $(e^{-\tilde P})_{zy}=(e^{-\tilde P})_y\tilde \varphi_z\cot\tilde \varphi-e^{-\tilde P}\tilde \varphi_{zy}\tan\tilde \varphi$ \ on $V\times\{0\}$, \ by (5.3). Then, by \ $(e^{-\tilde P})_{zyx}(x,y,0)=(e^{-\tilde P})_{zyx}(x,y,0)$ \ and \ $I^{x,y}(x,y,0)=0$ in (5.5), we have 
$$\tilde \varphi_{xyz}+\tilde \varphi_x\tilde \varphi_{yz}\tan\tilde \varphi-\tilde \varphi_y\tilde \varphi_{xz}\cot\tilde \varphi=0 \ \ \ \ {\rm on} \ \ V\times\{0\}. \eqno{(5.7)}$$
Thus, we obtain \ $(I^{x,y})_z(x,y,0)=0$ \ by (3.1), (5.5) and (5.7), and further obtain  
$$\frac{(\tilde \varphi_{xx}+\tilde \varphi_{yy}+\tilde \varphi_{zz})_z}{2}+\tilde \varphi_z\tilde \psi_{zz}-\tilde \varphi_x\tilde \varphi_{xz}\cot\tilde \varphi+\tilde \varphi_y\tilde \varphi_{yz}\tan\tilde \varphi =0  \ \ \ \ {\rm on} \ \ V\times\{0\} \eqno{(5.8)}$$
by (2.3.5), \ $((\tilde\kappa_3)_z+e^{-\tilde P}\tilde\varphi_z)(x,y,0)=0$ \ and (5.5). \ Hence, we also have \ $J_z(x,y,0)=0$ \ by (3.7). 

Next, by (2.3.3), (2.3.4), (5.4) and (5.5) we have 
$$(I^{x,z})_z+2\left(\tilde \psi_{xzz}+\tilde \varphi_{xzz}\cot\tilde \varphi-\frac{\tilde \varphi_z\tilde \varphi_{xz}}{\sin^2\tilde \varphi} \right)e^{-\tilde P}=
(I^{y,z})_z+2\left(\tilde \psi_{yzz}-\tilde \varphi_{yzz}\tan\tilde \varphi-\frac{\tilde \varphi_z\tilde \varphi_{yz}}{\cos^2\tilde \varphi} \right)e^{-\tilde P}=0$$
on $V\times \{0\}$. On the other hand, by (3.5), (3.6) and (5.5) we have 
$$(I^{x,z})_z+\left(\tilde \psi_{xzz}+\tilde \varphi_{xzz}\cot\tilde \varphi-\frac{\tilde \varphi_z\tilde \varphi_{xz}}{\sin^2\tilde \varphi} \right)e^{-\tilde P}=
(I^{y,z})_z+\left(\tilde \psi_{yzz}-\tilde \varphi_{yzz}\tan\tilde \varphi-\frac{\tilde \varphi_z\tilde \varphi_{yz}}{\cos^2\tilde \varphi} \right)e^{-\tilde P}=0$$
on $V\times \{0\}$. 
Therefore, we have \ $(I^{x,z})_z=(I^{y,z})_z=0$ \ and 
$$\left(\tilde \psi_{xzz}+\tilde \varphi_{xzz}\cot\tilde \varphi-\frac{\tilde \varphi_z\tilde \varphi_{xz}}{\sin^2\tilde \varphi} \right)=\left(\tilde \psi_{yzz}-\tilde \varphi_{yzz}\tan\tilde \varphi-\frac{\tilde \varphi_z\tilde \varphi_{yz}}{\cos^2\tilde \varphi} \right)=0 \eqno{(5.9)}$$
on $V\times\{0\}$. 

In consequence, we have verified (5.6) and further obtained (5.7)$\sim$(5.9). Equations (5.7)$\sim$(5.9) imply that $\tilde \varphi(x,y,z)$ satisfies all equations in Proposition 2.1-(2) on $V\times\{0\}$, as mentioned in the Proof of Lemma 2.1.

\vspace*{1mm}
{\sc Step 4}. $\tilde \varphi(x,y,z)$ satisfies (5.7)$\sim$(5.9) on $V\times\{0\}$. Then, there is an analytic function $\tilde \psi_z(x,y,z)$; $\tilde \psi_z(x_0,y_0,0)=0$ for a fixed $(x_0,y_0)\in V$, such that it satisfies 
$$
\begin{array}{l}
\tilde\psi_{zx}(x,y,z)=-\{\tilde\varphi_{zx}\cot\tilde\varphi\}(x,y,z)+O(z^2), 
\hspace{0.5cm} \tilde\psi_{zy}(x,y,z)=\{\tilde\varphi_{zy}\tan\tilde\varphi\}(x,y,z)+O(z^2),\\[3mm]
\hspace*{3cm}  \tilde\psi_{zz}(x,y,z)=[(\tilde\varphi_{xx}-\tilde\varphi_{yy}-\tilde\varphi_{zz}\cos2\tilde\varphi)/\sin2\tilde\varphi](x,y,z),
\end{array} \eqno{(5.10)}
$$
where $\tilde \psi_{zz}$ is given originally. Because we can interpret (5.7) and (5.9) as the integrability condition for $\tilde\psi_z(x,y,z)$ on $V\times\{0\}$. Furthermore, considering Remark 6, we have 
$$(\tilde\psi_{xx}-\tilde\psi_{yy})(x,y,z)=-\left[\{(\tilde\varphi_{xx}-\tilde\varphi_{yy})\cos2\tilde\varphi-\tilde\varphi_{zz}\}/\sin2\tilde\varphi\right](x,y,z)+O(z^2).  \eqno{(5.11)}$$
In fact, $((\tilde\psi_{zx})_x(x,y,0),(\tilde\psi_{zy})_y(x,y,0))$ and $(\tilde\psi_{xx}-\tilde\psi_{yy})_z(x,y,0)$, respectively, are determined from (5.10) and (5.11), and then \ $((\tilde\psi_{zx})_x-(\tilde\psi_{zy})_y)(x,y,0)=(\tilde\psi_{xx}-\tilde\psi_{yy})_z(x,y,0)$ \ holds by (5.8). 
Furthermore, by (5.11) and the last equation in (5.10), we have
$$
\begin{array}{l}
\tilde\psi_{zz}(x,y,z)=[(\tilde\varphi_{xx}-\tilde\varphi_{yy})\sin2\tilde\varphi-(\tilde\psi_{xx}-\tilde\psi_{yy})\cos2\tilde\varphi](x,y,0)+O(z^2),  \\[2mm]
\tilde\varphi_{zz}(x,y,z)=[(\tilde\varphi_{xx}-\tilde\varphi_{yy})\cos2\tilde\varphi+(\tilde\psi_{xx}-\tilde\psi_{yy})\sin2\tilde\varphi](x,y,z)+O(z^2). 
\end{array}
$$

Thus, $(\tilde\psi,\tilde\varphi)$ satisfies all equations in Proposition 2.1-(3) on $V\times\{0\}$. Furthermore, we have 
$$[\tilde\psi_{zx}+\tilde\varphi_{zx}\cot\tilde\varphi]_z=
[\tilde\psi_{zy}-\tilde\varphi_{zy}\tan\tilde\varphi]_z=0 \hspace*{1cm}
on \ \ V\times\{0\}$$
by the first two equations in (5.10). Hence, we can assume that $\tilde\psi(x,y,0)$ satisfies (2.1.2) on $V\times\{0\}$ (cf. \cite{bhs}, Proposition 3.1) and does not have any linear term. Then, $\tilde\psi(x,y,0)$ is uniquely determined for $(\tilde\varphi(x,y,0),\tilde\varphi_z(x,y,0))$ (cf. \cite{bhs}, Proposition 4.1). In consequence, an initial data $\tilde\psi(x,y,0)$, $\tilde\psi_z(x,y,0)$, $\tilde\varphi(x,y,0)$ and $\tilde\varphi_z(x,y,0)$ on $V\times\{0\}$ have been determined for the system of evolution equations (iii) and (iv) in Proposition 2.1-(3).

Now, let $(\psi(x,y,z),\varphi(x,y,z))$ be a solution to the system of evolution equations under the initial data above.
Then, $(\psi(x,y,z),\varphi(x,y,z))$ also satisfy (i) and (ii) of Proposition 2.1-(3) on $V\times I$ by (cf. \cite{bhs}, Proposition 4.2). Hence, we have obtained $\varphi(x,y,z)$, which determines a conformally flat metric $g$ by (2.1.1).

\vspace*{1mm}
{\sc Step 5}. We have obtained $\varphi(x,y,z)$ such that it leads to a metric $g\in CFM$, in Step 4. Let $e^{-P}(x,y,z)$ be a solution to (2.3.6) defined by $\varphi(x,y,z)$ under the initial condition \ $e^{-P}(x,y,0)=e^{-\tilde P}(x,y,0)$ \ and \ $(e^{-P})_z(x,y,0)=(e^{-\tilde P})_z(x,y,0)$. \ Let $\kappa_3(x,y,z)$ be a function defined by Proposition 2.2-(iv) from $\varphi(x,y,z)$ and $e^{-P}(x,y,z)$. 

With respect to these $\varphi(x,y,z)$, $e^{-P}(x,y,z)$ and $\kappa_3(x,y,z)$, we again define $I^{x,y}(x,y,z)$, $I^{x,z}(x,y,z)$, $I^{y,z}(x,y,z)$ and $J(x,y,z)$. Then, we also have (5.5) and \ $(\kappa_1\kappa_2)(x,y,0)\neq 0$ \ in this case, since \\
$(\varphi(x,y,z),e^{-P}(x,y,z),\kappa_3(x,y,z))$ is an extension of $\Pi$ by Proposition 2.2 and their initial conditions. 
Therefore, $e^{-P}(x,y,z)$ satisfies
$$I^{x,y}(x,y,z)=I^{x,z}(x,y,z)=I^{y,z}(x,y,z)=J(x,y,z)\equiv 0 \hspace*{0.5cm}for \ \ (x,y,z)\in V\times I$$
by Proposition 3.2. Thus, there is a generic conformally flat hypersurface $f$ in $R^4$ having the Guichard net $g$ and the first fundamental form $I_f=e^{2P}g$ by Proposition 2.3.  

In consequence, Theorem is verified in the same way as in Theorem 3. 
\hspace{\fill}$\Box$\\

\vspace*{2mm}
Suppose that a generic analytic surface $\phi( x, y)$ in $S^3$ satisfies the assumption of Theorem 4. 
Then, $\phi$ leads to not only a generic conformally flat hypersurface $f$ in the theorem but also the dual $f^*$ of $f$. We verify this fact by starting from reviews of dual generic conformally flat hypersurfaces (cf. \cite{hs3}, \cite{hsuy}).  

Let $f(x,y,z)$ be a generic conformally flat hypersurface in $R^4$, 
where $(x,y,z)$ is a coordinate system determined by the Guichard net of $f$. 
Let $S$ be the Schouten (1,1)-tensor of $f(x,y,z)$.
The tensor $S$ is diagonal with respect to the coordinate system, and its eigenvalues $\sigma_i \ (i=1,2,3)$ corresponding to 
the coordinate lines $x$, $y$ and $z$ are given by
$$\sigma_1=(\kappa_1\kappa_2-\kappa_2\kappa_3+\kappa_3\kappa_1)/2, \ \ \ \ 
\sigma_2=(\kappa_1\kappa_2+\kappa_2\kappa_3-\kappa_3\kappa_1)/2,$$
$$\sigma_3=(-\kappa_1\kappa_2+\kappa_2\kappa_3+\kappa_3\kappa_1)/2=(e^{-2P}+\kappa_3^2)/2 \ (>0),$$
respectively. 
Then, the dual $f^*(x,y,z)$ of $f(x,y,z)$ is defined by $df^*=df\circ S$,  
that is, $f^*_x=\sigma_1f_x$, $f^*_y=\sigma_2f_y$ and $f^*_z=\sigma_3f_z$ hold, and the coordinate system $(x,y,z)$ for $f^*$ is also a system determined by the Guichard net of $f^*$.  

For a function or a vector field $h$ on $f(x,y,z)$, 
we denote by $h^*$ those on $f^*(x,y,z)$ corresponding to $h$. The first fundamental form $I_{f^*}$ of $f^*$ is expressed as
$$I_{f^*}=e^{2P^*}(\cos^2\varphi^*(dx)^2+\sin^2\varphi^*(dy)^2+(dz)^2).$$ 
Note that $I_{f^*}$ is expressed in the same form as $I_f$. 
Then, we have 
$$e^{P^*}=\sigma_3e^P, \ \ \ \cos\varphi^*=(\sigma_1/\sigma_3)\cos\varphi, \ \ \ 
\sin\varphi^*=(\sigma_2/\sigma_3)\sin\varphi, \ \ \ \kappa_i^*=\kappa_i/\sigma_i,$$
$$X^*_{\alpha}=X_{\alpha}, \ \ \ X^*_{\beta}=X_{\beta}, \ \ \ X^*_{\gamma}=X_{\gamma}, \ \ \ N^*=N.$$
Furthermore, we have \ $\kappa^*_1=-e^{-P^*}\tan\varphi^*+\kappa^*_3$ \ and \ $\kappa^*_2=e^{-P^*}\cot\varphi^*+\kappa^*_3$ \ by \ $\kappa_1=e^{-P}\tan\varphi+\kappa_3$, \ $\kappa_2=-e^{-P}\cot\varphi+\kappa_3$ \ and \ $\kappa_i\sigma_j-\kappa_j\sigma_i=(\kappa_j-\kappa_i)\sigma_k$ \ for any permutation $\{i,j,k\}$ of $\{1,2,3\}$  (for these expressions of $\kappa^*_i$, see Remark 1 in \S2.2).   

Now, let $\phi(x,y)$ be a surface in $S^3$ satisfying the assumption of Theorem 4, and $f(x,y,z)$ be the hypersurface in $R^4$ arising from $\phi$. 
Then, the evolution $\phi^z(x,y)$, $z\in I$, of surfaces issuing from $\phi(x,y)$ has been given by $\phi^z(x,y):=N(x,y,z)$. Since $N(x,y,z)=N^*(x,y,z)$, \ $\phi^z(x,y)$ is the common Gauss map of $f$ and $f^*$.

\vspace*{2mm}
Now, we have the following three corollaries.  \\

\vspace*{2mm}
{\sc Corollary 5.1}. \ \ {\it Let $\phi(x, y)$ be a generic analytic surface in $S^3$ and satisfy the assumption in Theorem 4. 
Then, $\phi(x, y)$ also gives rise to the dual $f^*(x,y,z)$ of $f(x,y,z)$. For the dual $f^*$, we have \ 
(1) $a_1\sin\varphi^*- a_2\cos\varphi^*=-1$ (resp. 1) \ and \ (2) 
$( a_1)_{ x}=-\varphi^*_{ x} a_2+\bar P^*_{ x} a_1, \ \  
 ( a_2)_{ y}=\varphi^*_{ y} a_1+\bar P^*_{ y} a_2,$ \ 
where the sign $\mp 1$ in (1) corresponds to the one in Theorem 4-(1).} 

\vspace*{2mm}
{\sc Proof}. \ We have verified that $\phi$ also gives rise to the dual $f^*$ of $f$, in the argument above. Hence, we only show that the equations are satisfied.  
Let $f(x,y,z)$ and $f^*(x,y,z)$, respectively, be a generic conformally flat hypersurface and its dual arising from $\phi(x,y)$.

Then, we have \ $(\sigma_3)_x=\kappa_2(\kappa_3)_x$ \ and \ $(\sigma_3)_y=\kappa_1(\kappa_3)_y$ \ by (2.2.9). Hence, we have 
$$P^*_x=P_x[1+(\kappa_2/\sigma_3)e^{-P}\tan\varphi]=(\sigma_1/\sigma_3)P_x, \ \ \ 
P^*_y=P_y[1-(\kappa_1/\sigma_3)e^{-P}\cot\varphi]=(\sigma_2/\sigma_3)P_y \eqno{(5.12)}$$
by \ $(\kappa_3)_x=-(e^{-P})_x\tan\varphi$ \ and \ $(\kappa_3)_y=(e^{-P})_y\cot\varphi$. 

Next, we have
$$\varphi^*_x=\varphi_x-(\kappa_1/\sigma_3)(e^{-P})_x, \ \ \ \varphi^*_y=\varphi_y-(\kappa_2/\sigma_3)(e^{-P})_y. \eqno{(5.13)}$$
In fact, for $f$, \ $\varphi_x=(\kappa_1)_xe^P\cos^2\varphi$ \ and \ $\varphi_y=(\kappa_2)_ye^P\sin^2\varphi$ \ are satisfied as in (2.2.5): for example, this first equation is obtained from the first line of (2.2.9) and the first two equations of the second and the third lines in (2.2.4).  
For $f^*$, we have \ 
$\varphi^*_x=-(\kappa^*_1)_xe^{P^*}\cos^2\varphi^*$ \ and \ $\varphi^*_y=-(\kappa^*_2)_ye^{P^*}\sin^2\varphi^*$ \ from \ $\kappa^*_1=-e^{-P^*}\tan\varphi^*+\kappa^*_3$ \ and \ $\kappa^*_2=e^{-P^*}\cot\varphi^*+\kappa^*_3$ \ in the same way to obtain (2.2.5), since (2.2.4) and (2.2.9) are also satisfied for $f^*$.  
Then, since 
$$(\kappa_1)_x\sigma_1-\kappa_1(\sigma_1)_x=-(\kappa_1)_x\sigma_3+(\kappa_3)_x\kappa_1(\kappa_2-\kappa_1),$$ 
$$(\kappa_2)_y\sigma_2-\kappa_2(\sigma_2)_y=-(\kappa_2)_y\sigma_3-(\kappa_3)_y\kappa_2(\kappa_2-\kappa_1)$$ 
are satisfied by (2.2.9), we have (5.13) by \ $\kappa^*_i=\kappa_i/\sigma_i$. 

Now, since \ $\phi(x,y)=N(x,y,0)$, \ we have \ $\varphi(x,y)=\varphi(x,y,0), \ \ \bar P(x,y)=P(x,y,0), \ \ \bar\kappa_i(x,y)=\kappa_i(x,y,0)$, \ 
$ a_1=\bar\kappa_1e^{\bar P}\cos\varphi \ \ {\rm and} \ \ 
 a_2=\bar\kappa_2e^{\bar P}\sin\varphi \ \ {\rm for} \ \ (x,y)\in V.$ \  
Then, \ $ a_1\sin\varphi^*- a_2\cos\varphi^*=\mp 1$ \ follows from  \ $\kappa_1\sigma_2-\kappa_2\sigma_1=(\kappa_2-\kappa_1)\sigma_3$ \ and \ $ a_1\sin\varphi- a_2\cos\varphi=\pm 1$. \ The two equations in (2) follow directly from (5.12), (5.13) and Theorem 4-(2).
\hspace{\fill}$\Box$\\

\vspace*{2mm}
{\sc Corollary 5.2}. \ \ {\it We have  
$$P^*_z=-[\varphi_z\kappa_3e^{-P}+P_z(\sigma_3-\kappa_3^2)]/\sigma_3, \hspace{0.8cm} \varphi^*_z=[P_z\kappa_3e^{-P}-\varphi_z(\sigma_3-\kappa_3^2)]/\sigma_3.$$
In particular, let $\phi(x, y)$ be a generic analytic surface in $S^3$ and satisfy the assumption in Theorem 4. Then, we have 
$$ b_1=\bar P^*_z\cos\varphi^*-\varphi^*_z\sin\varphi^*, \ \ \ 
 b_2=\bar P^*_z\sin\varphi^*+\varphi^*_z\cos\varphi^*, $$
where \ $\bar P^*_z(x,y):=P^*_z(x,y,0)$ \ and \ $\varphi^*_z(x,y):=\varphi^*_z(x,y,0)$.
}\\

\vspace*{2mm}
{\sc Proof}. \ The first equation follows from \ $\sigma_3=(e^{-2P}+\kappa_3^2)/2$ \ and \ $(\kappa_3)_z=-\varphi_z e^{-P}$ \ in Proposition 2.2. For the second equation, since \ $(\sigma_2)_z=(\kappa_2)_z\kappa_1$ \ and \ $(\sigma_3)_z=(\kappa_2)_z(\kappa_3-\kappa_1)+(\kappa_3)_z\kappa_2$ \ are satisfied by (2.2.9), we firstly have 
$$\varphi^*_z=[-(e^{-P})_z\{\kappa_1\sigma_3+(\kappa_1-\kappa_3)\sigma_2\}+
\varphi_z\{\kappa_1(-\kappa_2+\kappa_3)\sigma_3+\kappa_3(\kappa_1-\kappa_3)\sigma_2+\sigma_2\sigma_3\}]/(\sigma_1\sigma_3)$$
by \ $\kappa_1=e^{-P}\tan\varphi+\kappa_3$, \ $\kappa_2=-e^{-P}\cot\varphi+\kappa_3$ \ and \ $(\kappa_3)_z=-\varphi_z e^{-P}$. \  
Next, we have \ $\kappa_1\sigma_3+(\kappa_1-\kappa_3)\sigma_2=\kappa_3\sigma_1$ \ by \ $\sigma_2+\sigma_3=\kappa_2\kappa_3$, \ and \ 
$\kappa_1(-\kappa_2+\kappa_3)\sigma_3+\kappa_3(\kappa_1-\kappa_3)\sigma_2+\sigma_2\sigma_3=\sigma_1(-\sigma_3+\kappa_3^2)$ \ by \ $\sigma_1+\sigma_2=\kappa_1\kappa_2$ \ and \ $\sigma_1+\sigma_3=\kappa_3\kappa_1$. \ Hence, the second equation is satisfied. 

For the last statement, we have the following equations on $U$ in general:
$$P_z\cos\varphi-\varphi_z\sin\varphi=P^*_z\cos\varphi^*-\varphi^*_z\sin\varphi^*, \hspace*{0.5cm}P_z\sin\varphi+\varphi_z\cos\varphi=P^*_z\sin\varphi^*+\varphi^*_z\cos\varphi^*.$$
Because \ $\kappa_3(\kappa_2-\kappa_3)\sigma_1+\sigma_2(\sigma_3-\kappa_3^2)=\kappa_3(\kappa_1-\kappa_3)\sigma_2+\sigma_1(\sigma_3-\kappa_3^2)=-\sigma_3^2$ \ hold by \ $\sigma_1+\sigma_3=\kappa_1\kappa_3$ \ and \ $\sigma_2+\sigma_3=\kappa_2\kappa_3$.
\hspace{\fill}$\Box$\\

We regard each coordinate element of the map \ $\phi:V\ni(x,y)\mapsto\phi(x,y)\in (S^3\subset) R^4$ \ as a function on the surface $\phi(x,y)$. 
Let $\Delta_{\phi}$ be the (positive) Laplacian acting on functions of the surface $\phi$. Then, it is known that the following equation is satisfied:
 $$\Delta_{\phi} \phi:=\frac{1}{a_1a_2}\left[\frac{\partial}{\partial x}(a_2X^0_{\alpha})+\frac{\partial}{\partial y}(a_1X^0_{\beta}) \right]=2(\phi-H\xi),$$
where \ $H:=(\lambda_1+\lambda_2)/2$ \ is the mean curvature of $\phi$.\\

\vspace*{2mm}
{\sc Corollary 5.3}. \ \ {\it Let $\phi(x, y)$ be a generic analytic surface in $S^3$ and satisfy the assumption in Theorem 4. 
Then, we have
$$-\frac{1}{{a_1}^2}\nabla'_{\partial/\partial x}\phi_x-\frac{1}{{a_2}^2}\nabla'_{\partial/\partial y}\phi_y$$ 
$$=2(\phi-H\xi)+\frac{{a_1}^2+{a_2}^2}{a_1a_2}\left[-\frac{1}{a_1}\frac{(\varphi+\varphi^*)_x}{2}X^0_{\alpha}+\frac{1}{a_2}\frac{(\varphi+\varphi^*)_y}{2}X^0_{\beta}\right].$$
} 

\vspace*{2mm}
{\sc Proof}. \ By Lemma 5.3-(1), we firstly have 
$$\frac{-1}{{a_1}^2}\nabla'_{\partial/\partial x}\phi_x=\frac{1}{{a_1}^2}\nabla'_{\partial/\partial x}(a_1X^0_{\alpha})$$
$$=\frac{1}{{a_1}^2}(-\varphi_xa_2+\bar P_xa_1)X^0_{\alpha}-\frac{1}{a_1}(\bar P_y\cot\varphi-\varphi_y)X^0_{\beta}-\frac{b_1}{a_1}\xi+\phi,$$
$$\frac{-1}{{a_2}^2}\nabla'_{\partial/\partial y}\phi_y=\frac{1}{{a_2}^2}\nabla'_{\partial/\partial y}(a_2X^0_{\beta})$$
$$=\frac{1}{{a_2}^2}(\varphi_ya_1+\bar P_ya_2)X^0_{\beta}-\frac{1}{a_2}(\bar P_x\tan\varphi+\varphi_x)X^0_{\alpha}-\frac{b_2}{a_2}\xi+\phi.$$
Next, by Lemma 5.3-(3) we have \ $(\bar P_x/\cos\varphi)a_1=\bar P_xe^{\bar P}\bar\kappa_1$, \ $(\bar P_y/\sin\varphi)a_2=\bar P_ye^{\bar P}\bar\kappa_2$, \ $a_1\sin\varphi-a_2\cos\varphi=1$ \ and \ ${a_1}^2+{a_2}^2=1+Q^2=e^{2\bar P}(e^{-2\bar P}+{\bar\kappa_3}^2)=2e^{2\bar P}\bar\sigma_3$, \ where \ $\bar\sigma_3(x,y):=\sigma_3(x,y,0)$. \ Hence, by (5.13) we have 
$$\frac{1}{{a_1}^2}\nabla'_{\partial/\partial x}(a_1X^0_{\alpha})+\frac{1}{{a_2}^2}\nabla'_{\partial/\partial y}(a_2X^0_{\beta})$$
$$=\frac{{a_1}^2+{a_2}^2}{a_1a_2}\left[-\frac{1}{a_1}(\varphi_x-\frac{(e^{-\bar P})_x\bar\kappa_1}{2\bar\sigma_3})X^0_{\alpha}+\frac{1}{a_2}(\varphi_y-\frac{(e^{-\bar P})_y\bar\kappa_2}{2\bar\sigma_3})X^0_{\beta}\right]+2(\phi-H\xi)$$
$$=\frac{{a_1}^2+{a_2}^2}{a_1a_2}\left[-\frac{1}{a_1}\frac{(\varphi+\varphi^*)_x}{2}X^0_{\alpha}+\frac{1}{a_2}\frac{(\varphi+\varphi^*)_y}{2}X^0_{\beta}\right]+2(\phi-H\xi).$$
\hspace{\fill}$\Box$\\

We remark on the elliptic differential equation in Corollary 5.3. For $\hat g\in Met^0$, all functions $a_i$, $H$ $\varphi$ and $\varphi^*$ on $V$ defining the equation are determined by the results in \S3. We could analyze such solutions $\phi$. \\

\newpage
{\small

\end{document}